\documentclass[final,3p,times]{elsarticle}
\usepackage{amsmath,color,amsthm,amsfonts,mathrsfs,dsfont,amssymb,amstext,float,amsopn}
\usepackage{graphicx,subfigure,texdraw}
\usepackage[all]{xy}
\usepackage{epstopdf}

\journal{ }

\newtheorem{theorem}{Theorem}[section]
\newtheorem{lemma}[theorem]{Lemma}

\newtheorem{proposition}[theorem]{Proposition}

\theoremstyle{remark}
\newtheorem{remark}[theorem]{Remark}

\newcommand{\ds}{\displaystyle}

\numberwithin{equation}{section}

\begin{document}

\begin{frontmatter}
\title{{\bf On Turing-Turing bifurcation of partial functional differential equations and its induced superposition patterns}}
\author{Xun Cao}
\author{Weihua Jiang\tnoteref{t1}}\tnotetext[t1]{Corresponding author. E-mail address: jiangwh@hit.edu.cn.}

\address{School of Mathematics, Harbin Institute of Technology, Harbin 150001, P.R.China}

\begin{abstract}
When two Turing modes interact, i.e., Turing-Turing bifurcation occurs, superposition patterns revealing complex dynamical phenomena appear.
In this paper, superposition patterns resulting from Turing-Turing bifurcation are investigated in theory.
Firstly, the third-order normal form locally topologically equivalent to original partial functional differential equations (PFDEs) is derived.
When selecting 1D domain and Neumann boundary conditions, three normal forms describing different spatial patterns are deduced from original third-order normal form. Also, formulas for computing coefficients of these normal forms are given, which are expressed in explicit form of original system parameters.
With the aid of three normal forms, spatial patterns of a diffusive predator-prey system with Crowley-Martin functional response near Turing-Turing singularity are investigated.
For one set of parameters, diffusive system supports \emph{the coexistence of four stable steady states} with different single characteristic wavelengths, which demonstrates our previous conjecture. For another set of parameters, \emph{superposition patterns, tri-stable patterns} that a pair of stable superposition steady states coexists with the stable coexistence equilibrium or another stable steady state, as well as \emph{quad-stable patterns} that a pair of stable superposition steady states and another pair of stable steady states coexist, arise.
Finally, numerical simulations are shown to support theory analysis.
\end{abstract}

\begin{keyword}
Partial functional differential equations; Turing-Turing bifurcation; Normal form; Superposition patterns; Quad-stability; Diffusive predator-prey system
\end{keyword}

\end{frontmatter}

\parskip=1mm

\section{Introduction}
In 1952, Turing proposed a striking idea of ''diffusion-driven instability'' in his seminal paper \cite{TAM}, which states that diffusion could destabilize an otherwise stable equilibrium of a reaction-diffusion system and induce new stable spatially inhomogeneous pattern. After a few years, experimental observations \cite{CV&DE&BJ, DP&CV&DE, LI&EIR} started to support Turing's theory.
This kind of spatial phenomenon is usually called Turing pattern, which is resulting from Turing instability. From the point of view of bifurcation, Turing patterns emerge as Turing bifurcation occurs.

Actually, most experimental studies about Turing patterns have explored structures with a single characteristic wavelength, while few experiments have concentrated on spatial patterns arising from interactions between different Turing modes.
Black-eye patterns which were firstly observed in a reaction-diffusion experiment involving the chlorite-iodide-malonic acid (CIMA) reaction in a thin layer gel reactor \cite{GGH&OQ&SHL}, is believed to be induced by two interacting Turing modes with different wavelengths. To reproduce the black-eye patterns in 2D domain, Yang et al. \cite{YL&DM&ZAM&EIR} firstly proposed a model and investigated pattern formations via numerical simulations. They also observed a variety of other spatial resonant patterns resulting from interactions between different Turing modes, including white-eye patterns and superposition patterns.
Also, we say that Turing-Turing bifurcation occurs, when two Turing modes with different wavelengths interact.

Lately, some spatial patterns arising from two interacting Turing modes with different wavelengths, have been theoretically investigated, see \cite{YR&SYL, WMH&WJH&GGH, GGH&LL&LBF&LJ, LSB&WJH&DYY, YWB, DG, GG&LMC&SM&NA13, GG&LMC&SM&PRE18}.
Yang and Song \cite{YR&SYL} investigated conditions of the occurrence of spatial resonance bifurcation and corresponding spatial patterns for a Gierer-Meinhardt system, utilizing linear stability analysis, center manifold theory and normal form method.
Meanwhile, applying Lyapunov-Schmidt technique and singularity theory method, Wei et al. \cite{WMH&WJH&GGH} investigated steady state bifurcation from a double eigenvalue for a two-species glycolysis model subject to Dirichlet boundary conditions. For an autocatalysis model with high order subject to Neumann boundary conditions, Guo et al. \cite{GGH&LL&LBF&LJ} also investigated steady state bifurcation from a double eigenvalue utilizing techniques of space decomposition and implicit function theorem. They found bifurcating steady state solutions in the form of $s(\omega)(\cos\omega\Psi_i+\sin\omega\Psi_j+W(\omega))$ for sufficiently small $|\omega-\omega_0|$, satisfying that $s(\omega_0)=0,W(\omega_0)=0$ and $\Psi_i,\Psi_j(j=2i\;\mathrm{or}\;i=2j)$ are eigenvectors corresponding to the double eigenvalue. They also provided numerical simulations to illustrate theory results.

In a two dimensional domain, Gambino et al. \cite{GG&LMC&SM&NA13} investigated the process of pattern formation resulting from regular, degenerate non-resonant and resonant Turing bifurcations for a reaction-diffusion system with cross-diffusion terms and the competitive Lotka-Volterra kinetics, using weakly nonlinear analysis. It was found that reaction-diffusion system supports spatial patterns, like squares, mixed-mode patterns, supersquares and hexagonal patterns.

These investigations on Turing-Turing bifurcation (or interactions between different Turing modes, or double eigenvalues, or degenerate Turing bifurcation) indicate that, there might appear new superposition patterns that reveal new dynamical phenomena, when two Turing modes interact.
However, these investigations only theoretically reveal a part of these spatial superposition patterns.
Then in this paper, based on center manifold theory and normal form method, we try to look for new superposition patterns induced by Turing-Turing bifurcation and to reveal spatial dynamics near Turing-Turing singularity as much as possible.

Exploring qualitative behaviors of partial differential equations, functional differential equations and other kinds of differential equations, which arise in engineering, biology and social sciences, is one significant branch of modern mathematics.
Normal form plays an important role in qualitative analysis of differential systems, like bifurcation analysis.
Also, normal form method is an efficient mathematical tool in investigating complex local dynamics of differential systems, when combined with center manifold theory. For ordinary differential systems with delays and without delay, center manifold reduction and normal form method have been established, see \cite{CJ, CSA&YY, CSN&LCZ&WD, FT&MLT&H, HBD&KND&WYH}.
Utilizing these mathematical tools, dynamical behaviors of functional differential systems, including periodic solutions, homoclinic orbits and heteroclinic orbits, have been investigated in \cite{JWH&YY, WCC&WJJ, WHB&JWH, YR&JWH&WY}.

For PFDEs, center manifold theory and normal form theory have also been developed \cite{FT&HWZ&WJH, LXD&SJWH&WJH, FT, FT01, WJH}. Applying normal form method, spatiotemporal dynamics resulting from Hopf bifurcation have been extensively investigated, see \cite{YFQ&WJJ&SJP, SY&WJJ&SJP, WJF, ZR&GSJ, SHB&RSG, CSS&SJP, LX&JWH&SJP}. Furthermore, based on center manifold theory and normal form method, Jiang et al.\cite{JWH&AQ&SJP} recently derived several concise formulas of computing coefficients of normal forms for partial functional differential equations at Turing-Hopf singularity.
\emph{And, these formulas make it easy for us to compute normal forms}.

Inspired by Jiang et al. \cite{JWH&AQ&SJP}, we discuss the calculation of normal forms of Turing-Turing bifurcation for some parameterized PFDEs, then investigate superposition patterns of a diffusive predator-prey system near Turing-Turing singularity by analyzing the obtained normal forms.
Firstly, we derive the third-order normal form which is locally topologically equivalent to the original parameterized PFDEs at Turing-Turing singularity, based on Faria's work \cite{FT, FT01} and center manifold theory \cite{LXD&SJWH&WJH, FT&HWZ&WJH}. Also, an explicit algorithm of computing the third-order normal form is provided.
And in a one dimensional domain and under Neumann boundary conditions, the third-order normal form have three different forms that could describe essentially different spatial phenomena, depending on spatial wavelengths of two interacting Turing modes.
Moreover, several concise formulas for computing coefficients of these three normal forms are also derived. It is remarkable that the process of computing coefficients of normal forms utilizing these concise formulas could be implemented by computer programs, and these formulas also apply to computing coefficients of normal forms of partial differential equations.
Then, spatial patterns of a diffusive predator-prey system with Crowley-Martin functional response near Turing-Turing singularity are investigated, with the aid of these three normal forms.
And, two of these three normal forms arise in predator-prey system for two different sets of parameters. Specifically, for one set of parameters, a pair of stable spatially inhomogeneous steady states with the shape of $\phi_1\cos2x-$like coexists with another pair of stable spatially inhomogeneous steady states with the shape of $\phi_2\cos3x-$like, which theoretically demonstrates our conjecture in \cite{CX&JWH}.
For another set of parameters, predator-prey system supports superposition steady states with the shape of $\phi_1h_1\cos x+\phi_2h_2\cos2x-$like, tri-stable patterns that a pair of stable superposition steady states coexists with the stable coexistence equilibrium or another stable steady state with the shape of $\phi_2\cos2x-$like, as well as quad-stable patterns that a pair of stable superposition steady states with the shape of $\phi_1h_1\cos x+\phi_2h_2\cos2x-$like and a pair of stable steady states with the shape of $\phi_2\cos2x-$like coexist.

This paper is organized as follows. In Section 2, the third-order normal form locally topologically equivalent to original PFDEs at Turing-Turing singularity, is derived, of which the proof is put in Section 4. Also, an explicit algorithm of computing the third-order normal form is provided.
Then, in 1D domain and under Neumann boundary conditions, the third-order normal form is simplified as three normal forms describing essentially different spatial phenomena. Moreover, several concise formulas of computing coefficients of these three normal forms are derived.
Next in Section 3, we investigate spatial patterns of a diffusive predator-prey system with Crowley-Martin functional response near Turing-Turing singularity, by analyzing the obtained normal forms.
At last, conclusions are shown in Section 5.

\section{Normal form of Turing-Turing bifurcation}
In this section, we investigate normal form of some parameterized PFDEs at Turing-Turing singularity.

Assume that $\Omega\subset \mathbb{R}^n(n\in\mathbb{N})$ is a bounded open set with smooth boundary, and $X$ is Hilbert space of functions defined on $\bar\Omega$ with inner product $ \left\langle\cdot,\cdot\right\rangle$. Actually, $X$ could be defined as the following real-valued Hilbert space,
\begin{equation*}
X\triangleq\left\{\begin{aligned}&u\in H^{2}(\Omega),\frac{\partial u}{\partial \nu}|_{x\in\partial\Omega}=0, & &\textrm{for the homogeneous Neumann boundary conditions},\\
&u\in H^{2}(\Omega)\cap H_0^1(\Omega), & &\textrm{for the homogeneous Dirichlet boundary conditions},
\end{aligned}\right.
\end{equation*}
with the inner product
\begin{equation*}
\langle u, v\rangle=\frac{1}{|\Omega|}\int_{\Omega}u(x)v(x) \mathrm{d}x, \quad \textrm{for}\; u,v\in X.
\end{equation*}

Let $\left\{\beta_{k}\right\}_{{k}\in\mathbb{N}_*}$ be eigenfunctions of $-\Delta$ on $\Omega$, with the corresponding eigenvalues $\{\mu_{k}\}_{{k}\in\mathbb{N}_*}$ satisfying $\mu_{k}\ge 0$ and $\mu_{k}\to+\infty$ as $k\to \infty$, where
\begin{equation*}
\mathbb{N}_*\triangleq\begin{cases} \mathbb{N}_0= \mathbb{N}\cup \{0\} &\textrm{for the homogeneous Neumann boundary conditions},\\
\mathbb{N} &\textrm{for the homogeneous Dirichlet boundary conditions}.
\end{cases}
\end{equation*}
Then, $\{\beta_{k}\}_{{k}\in\mathbb{N}_*}$ form an orthonormal basis of $X$.

Furthermore, let $\mathcal{C}=C([-r,0];X^m)\;(r>0, m\in\mathbb{N})$ be the Banach space of continuous maps from $[-r,0]$ to $X^m$ with the sup norm. Then in phase space $\mathcal{C}$, consider the following parameterized abstract PFDEs,
\begin{align}\label{e:2.1}
\frac{\mathrm{d}}{\mathrm{d}t}u(t)=D(\epsilon)\Delta u(t)+L(\epsilon)(u_t)+G(u_t,\epsilon),
\end{align}
where $u_t\in \mathcal{C}$ is defined as $u_t(\theta)=u(t+\theta),-r\leq\theta\leq0$. And, $D(\epsilon)=\mathrm{diag}\left(d_1(\epsilon), d_2(\epsilon), \cdots, d_m(\epsilon))\right)$ satisfies $d_j(\epsilon)>0$ for $\epsilon=(\epsilon_1,\epsilon_2)\in V_0\subset\mathbb{R}^2$ and $j=1,2,\cdots, m$, where $V_0$ is a neighborhood of the origin. Also, $L:V_0\to \mathcal{L}(\mathcal{C};X^m)$, that is, $L(\epsilon)$ is a bounded linear operator, and $G(\cdot,\epsilon): \mathcal{C}\to X^m$ is a $C^k$($k\geq 3$) function satisfying $G(0,0)=0,\mathrm{D}_{\phi}G(0,0)=0$ for $\phi\in\mathcal{C}$.

Denoting $L_0^*=L(0),D_0=D(0)$, the linearized equation of \eqref{e:2.1} at the origin reads
\begin{align}\label{e:2.2}
\frac{\mathrm{d}}{\mathrm{d}t}u(t)=D_0\Delta u(t)+L_0^*(u_t).
\end{align}

And, we have the following hypotheses (H1)-(H3) (see \cite{FT01, FT, WJH}):
\begin{enumerate}
\item[(H1)] $D_0\Delta$ generates a $\mathcal{C}_0$ semigroup $\{T(t)\}_{t\ge 0}$ on $X^m$ with $|T(t)|\le Me^{\omega t}$ for some $M\ge 1, \omega\in \mathbb{R}$ and all $t\ge0$, and $T(t)$ is a compact operator for $t>0$;
\item[(H2)] the subspaces $\mathcal{B}_{k}\triangleq\left\{\left\langle v(\cdot),\beta_{k}\right\rangle\beta_{k}|v\in \mathcal{C}\right\}\subset\mathcal{C}$ satisfy $L_0^*\left(\mathcal{B}_{k}\right)\subseteq\mathrm{span}\left\{\beta_{k}e_j,1\le j\le m\right\}$, where $\{e_j,1\le j\le m\}$ is the canonical basis of $\mathbb{R}^m$, and
    \begin{equation*}
    \left\langle v(\cdot),\beta_{k} \right\rangle\triangleq \left(\left\langle v_1(\cdot),\beta_{k} \right\rangle,\left\langle v_2(\cdot),\beta_{k} \right\rangle,\cdots, \left\langle v_m(\cdot),\beta_{k} \right\rangle\right)^T, \quad \textrm{for}\;v=(v_1, v_2, \cdots, v_m)^T\in\mathcal{C}, k\in\mathbb{N}_*;
    \end{equation*}
\item[(H3)] $L_0^*$ can be extended to a bounded linear operator from $B\mathcal{C}$ to $X^m$, where
$$B\mathcal{C}=\left\{\psi:[-r,0]\to X^m|\; \psi\; \textrm{is continuous on}\; [-r,0),\exists \lim\limits_{\theta\to 0^-}\psi(\theta)\in X^m \right\},$$
with the sup norm.
\end{enumerate}

Then in $B\mathcal{C}$, Eq.  \eqref{e:2.1} could be rewritten as an abstract ordinary differential equation (ODE),
\begin{align}\label{e:2.3}
\frac{\mathrm{d}}{\mathrm{d}t}u=\mathcal{A}u+X_0F(u,\epsilon),
\end{align}
where $F(\phi,\epsilon)=G(\phi,\epsilon)+(L(\epsilon)-L_0^*)\phi+(D(\epsilon)-D_0)\Delta \phi(0)$, and $\mathcal{A}$ is defined by
\begin{align}\label{e:2.4}
\mathcal{A}:\mathcal{C}_0^1\subset B\mathcal{C}\to B\mathcal{C},\quad
\mathcal{A}\phi=\dot{\phi}+X_0\left[L_0^*(\phi)+D_0\Delta \phi(0)-\dot{\phi}(0)\right],
\end{align}
on $\mathcal{C}_0^1\triangleq\left\{\phi\in \mathcal{C},\dot{\phi}\in \mathcal{C},\phi(0)\in\mathrm{dom}(\Delta)\right\}$, and $X_0\in BV([-r,0],\mathbb{R}^{m\times m})$ satisfies $X_0(\theta)=0$ for $\theta\in[-r,0)$ and $X_0(0)=I$.

Moreover, denote $C\triangleq C\left([-r,0];\mathbb{R}^m\right)$, and for each ${k}\in\mathbb{N}_*$, define $L_{k}:{C}\to\mathbb{R}^m$ by $L_{k}(\psi)\beta_{k}=L_0^*(\psi\beta_{k})$.
Therefore, on $\mathcal{B}_{k}$, linear Eq.  \eqref{e:2.2} is equivalent to the following functional differential equations on $\mathbb{R}^m$,
\begin{align}\label{e:2.5}
\dot{z}(t)=-\mu_{k} D_0z(t)+L_{k}(z_t), \;{k}\in \mathbb{N}_*,
\end{align}
with $z_t=\left\langle u_t(\cdot),\beta_k\right\rangle\in C$. And, the corresponding characteristic equations are
\begin{equation}\label{e:2.6}
\det\Delta(\lambda,\mu_{k})=0,\quad \mathrm{with}\;\Delta(\lambda,\mu_{k})=\lambda I+\mu_kD_0-{L}_k(e^{\lambda\cdot}I),k\in\mathbb{N}_*.
\end{equation}

Define $\eta_k\in BV\left([-r,0];\mathbb{R}^{m\times m}\right)$ satisfying that
\begin{align}\label{e:2.7}
-\mu_kD_0 \psi(0)+L_k(\psi)=\int_{-r}^0 d\eta_k(\theta)\psi(\theta),\qquad\psi\in C, k\in\mathbb{N}_*.
\end{align}
Then, the adjoint bilinear form $(\cdot,\cdot)_k$ on ${C}^*\times{BC}$ with ${C}^*=C([0,r];\mathbb{R}^{m*})$ is defined by
\begin{align}\label{e:2.8}
(\alpha,\beta)_k=\alpha(0)\beta(0)-\int_{-r}^0\int_0^{\theta}\alpha(\zeta-\theta)d\eta_k(\theta)\beta(\zeta)d\zeta, \qquad \alpha\in {C}^*, \beta\in BC, k\in\mathbb{N}_*,
\end{align}
with $BC=\left\{\psi:[-r,0]\to \mathbb{R}^m|\; \psi\; \textrm{is continuous on}\; [-r,0),\exists \lim\limits_{\theta\to 0^-}\psi(\theta)\in \mathbb{R}^m \right\}$.

Also, we have the hypothesis about Turing-Turing bifurcation.
\begin{enumerate}
\item[(H4)]
There exists a neighborhood $(0,0)\in V_0\subset \mathbb{R}^2$ such that for $\epsilon=(\epsilon_1,\epsilon_2)\in V_0$, characteristic equations of the linearized system of \eqref{e:2.1} at the origin have two independent simple real eigenvalues $\gamma_j(\epsilon),j=1,2$ corresponding to different $k_1,k_2\in\mathbb{N}$, satisfying that $\gamma_j(0)=0,\frac{\partial}{\partial \epsilon_j}\gamma_j(0)\neq 0,j=1,2$, and the remaining eigenvalues have non-zero real parts.
\end{enumerate}
Then, $B\mathcal{C}$ is decomposed by $\Lambda=\{0,0\}$ as
\begin{equation*}
B \mathcal{C}=\mathcal{P}\oplus\mathrm{Ker}\,\pi,
\end{equation*}
where $\pi:B \mathcal{C}\to \mathcal{P}$ is the projection defined by
\begin{equation}\label{e:2.9}
\begin{aligned}
\pi\phi=&\sum_{j=1}^2 \Phi_{k_j}\left(\Psi_{k_j},\left\langle \phi(\cdot),\beta_{k_j}\right\rangle\right)_{k_j}\beta_{k_j}, \quad & &\phi\in B\mathcal{C},
\end{aligned}
\end{equation}
satisfying
\begin{equation}\label{e:2.10}
\begin{aligned}
\pi(X_0\alpha)=&\sum_{j=1}^{2}\Phi_{k_j}\Psi_{k_j}(0)\left\langle\alpha,\beta_{k_j} \right\rangle\beta_{k_j},\quad & &\alpha\in X^m,
\end{aligned}
\end{equation}
with
\begin{equation*}
\dot{\Phi}_{k_j}=\Phi_{k_j}B_{k_j}, \quad -\dot{\Psi}_{k_j}=B_{k_j}\Psi_{k_j}, \quad \left(\Psi_{k_j}, \Phi_{k_j}\right)_{k_j}=I, \quad B_{k_j}=0,\quad j=1,2.
\end{equation*}
For simplicity of notations, denote $\phi_j=\Phi_{k_j}, \psi_j=\Psi_{k_j},j=1,2$. Then, by \cite{HJK&LSMY}, we have
\begin{equation}\label{e:phi}
\phi_j(\theta)\equiv\phi_j(0), \theta\in [-r,0],\quad
\psi_j(s)\equiv\psi_j(0), s\in [0,r], \qquad j=1,2,
\end{equation}
satisfying
\begin{equation*}
\Delta(0,\mu_{k_j})\phi_j(0)=0,\quad \psi_j(0)\Delta(0,\mu_{k_j})=0,\quad\left(\psi_j,\phi_j\right)_{k_j}=1,\qquad j=1,2.
\end{equation*}

Thus, $u\in\mathcal{C}_0^1$ is decomposed as
\begin{equation}\label{e:c}
u(t)=\sum_{j=1}^2\phi_j z_{j}(t)\beta_{k_j}+y(t),
\end{equation}
with $z_{j}(t)=\left(\psi_j, \left\langle u(t)(\cdot),\beta_{k_j} \right\rangle \right)_{k_j}\in\mathbb{R},j=1,2$ and $y(t)\in\mathcal{Q}^1\triangleq\mathcal{C}_0^1 \cap \mathrm{Ker}\,\pi$. Since $\pi $ commutes with $A$ in $\mathcal{C}_0^1$, Eq. \eqref{e:2.3} is equivalent to
\begin{equation}\label{e:2.11}
\left\{\begin{aligned}
\dot{z}&=Bz+\Psi(0)
\begin{pmatrix}
 \left\langle F\left(\sum_{j=1}^2\phi_j z_{j} \beta_{k_j}+y,\epsilon\right),\beta_{k_1} \right\rangle \\
 \left\langle F\left(\sum_{j=1}^2\phi_j z_{j} \beta_{k_j}+y,\epsilon\right),\beta_{k_2} \right\rangle
\end{pmatrix},\\
\frac{\mathrm{d}}{\mathrm{d}t}y&=\mathcal{A}_1 y+(I-\pi)X_0 F\left(\sum_{j=1}^2\phi_j z_{j} \beta_{k_j}+y,\epsilon\right),
\end{aligned}\right.
\end{equation}
where $z=(z_1,z_2),B=\mathrm{diag}(0,0), \Psi=(\psi_1,\psi_2)$ and $\mathcal{A}_1$ is defined as $\mathcal{A}_1:\mathcal{Q}^1\subset\mathrm{Ker}\,\pi\to \mathrm{Ker}\,\pi,\,\mathcal{A}_1\phi=\mathcal{A}\phi$ for $\phi\in\mathcal{Q}^1$.

Consider the formal Taylor expansion
\begin{equation}\label{e:2.12}
F(\phi,\epsilon)=\sum_{k\geq2}\frac{1}{j!}F_j(\phi,\epsilon),\quad \phi\in\mathcal{C},\epsilon\in \mathbb{R}^2.
\end{equation}
where $F_j(\cdot,\cdot)$ is the $j$th $\mathrm{Fr}\acute{\mathrm{e}}\mathrm{chet}$ derivation of $F(\cdot,\cdot)$. Then, Eq.  \eqref{e:2.11} is written as,
\begin{equation}\label{e:2.13}
\begin{cases}
\ds\dot{z}=Bz+\sum_{j\geq2}\frac{1}{j!}f_j^1(z,y,\epsilon),\\
\ds\frac{\mathrm{d}}{\mathrm{d}t}y=\mathcal{A}_1y+\sum_{j\geq2 }\frac{1}{j!}f_j^2(z,y,\epsilon),
\end{cases}
\end{equation}
where $f_j^1,f_j^2,\,j\geq2$,
are defined by
\begin{equation}\label{e:2.14}
\left\{\begin{aligned}
f_j^1(z,y)&=\Psi(0)
\begin{pmatrix}
 \left\langle F_j\left(\sum_{j=1}^2\phi_j z_{j}\beta_{k_j}+y,\epsilon\right),\beta_{k_1} \right\rangle\\
 \left\langle F_j\left(\sum_{j=1}^2\phi_j z_{j}\beta_{k_j}+y,\epsilon\right),\beta_{k_2} \right\rangle
\end{pmatrix},\\
f_j^2(z,y)&=(I-\pi)X_0F_j\left(\sum_{j=1}^2\phi_j z_{j}\beta_{k_j}+y,\epsilon\right).
\end{aligned}\right.
\end{equation}

And, by a recursive process through a transformation of variable of the form
\begin{align}\label{e:2.15}
(z,y)=\left(\hat{z},\hat{y}\right)+\frac{1}{j!}\left(U_j^1(\hat{z},\epsilon),U_j^2(\hat{z},\epsilon)\right),
\end{align}
 Eq.  \eqref{e:2.13} is transformed into the normal form
\begin{equation}\label{e:2.16}
\begin{cases}
\ds\dot{z}=Bz+\sum_{j\geq2}\frac{1}{j!}g_j^1(z,y,\epsilon),\\
\ds\frac{\mathrm{d}}{\mathrm{d}t}y=\mathcal{A}_1y+\sum_{j\geq2}\frac{1}{j!}g_j^2(x,y,\epsilon),
\end{cases}
\end{equation}
with
\begin{align}\label{e:2.17}
g_j=\tilde{f}_j-M_jU_j,\quad j\geq2,
\end{align}
where $\tilde{f}_j=\left(\tilde{f}_j^1,\tilde{f}_j^2\right),g_j=\left(g_j^1,g_j^2\right),U_j=\left(U_j^1,U_j^2\right)$, and operator $M_j=\left(M_j^1,M_j^2\right),j\geq2 $ is defined by
\begin{equation}\label{e:2.18}
\begin{aligned}
&M_j^1:V_j^{2+2}\left(\mathbb{R}^2\right)\to V_j^{2+2}\left(\mathbb{R}^2\right),
&\qquad &\left(M_j^1p\right)(z,\epsilon)=D_zp(z,\epsilon)Bz-Bp(z,\epsilon),\\
&M_j^2:V_j^{2+2}\left(\mathcal{Q}^1\right)\subset V_j^{2+2}\left(\mathrm{Ker}\,\pi\right)\to V_j^{2+2}\left(\mathrm{Ker}\,\pi\right),
& &\left(M_j^2h\right)(z,\epsilon)=D_zh(z,\epsilon)Bz-\mathcal{A}_1(h(z,\epsilon)),
\end{aligned}
\end{equation}
where $V_j^{m+p}(Y)$ is the linear space of homogeneous polynomials of degree $j$ in $m+p$ variables $z = (z_1, z_2, \cdots, z_m), \epsilon= (\epsilon_1, \epsilon_2, \cdots, \epsilon_p) $ with coefficients in normed space $Y$, that is,
\begin{equation*}
V_j^{m+p}(Y)=\left\{\sum_{|(q,s)|=j}c_{(q,s)}z^q\epsilon^s:(q,s)\in\mathbb{N}_0^{m+p}, c_{(q,s)}\in Y\right\},
\end{equation*}
and the norm $|\sum_{|(q,s)|=j}c_{(q,s)}z^q\epsilon^s|=\sum_{|(q,s)|=j}|c_{(q,s)}|_Y$. For simplicity, denote $c_{qs}=c_{(q,s)}$.

Denoting $\mathbf{P}$ the projection operator, we also have
\begin{align}\label{e:2.19}
U_j(z,\epsilon)=\left(M_j\right)^{-1}\mathbf{P}_{\mathrm{Im}\left(M_j^1\right)\times \mathrm{Im}\left(M_j^2\right)}\circ\tilde{f}_j(z,0,\epsilon),\; g_j^1(z,0,\epsilon)=\mathbf{P}_{\mathrm{Im}\left(M_j^1\right)^c}\circ\tilde{f}_j^1(z,0,\epsilon),\quad j\ge 2,
\end{align}
where $\tilde{f}_j$ denotes the terms of order $j$ in $(z,y)$ obtained after the computation of normal forms up to order $j-1$.

Then, by writing $D(\epsilon)$ and $L(\epsilon)$ in the following expansions at $\epsilon=0$,
\begin{subequations}\label{e:2.20}
\begin{equation}
\begin{aligned}
D(\epsilon)&=D(0)+\frac{1}{2!}D_1(\epsilon)+\frac{1}{3!}D_2(\epsilon)+\cdots,\\
L(\epsilon)\phi&=L(0)\phi+\frac{1}{2!}L_1(\epsilon)\phi+\frac{1}{3!}L_2(\epsilon)\phi+\cdots, \qquad \phi\in\mathcal{C}, \epsilon\in\mathbb{R}^2,
\end{aligned}
\end{equation}
and $G(\cdot,\epsilon)$ in Taylor expansions in $\phi$ at $\epsilon=0$,
\begin{equation}
G(\phi,0)=\frac{1}{2!}Q(\phi,\phi)+\frac{1}{3!}C(\phi,\phi,\phi)+O\left(||\phi||^4\right),\qquad \phi\in\mathcal{C},
\end{equation}
\end{subequations}
where $Q(\cdot,\cdot)$ and $C(\cdot,\cdot,\cdot)$ are symmetric multi-linear forms, and
\begin{equation}\label{e:2.21}
\begin{aligned}
&D_1(\epsilon)=\sum_{i=1}^2\epsilon_iD_{\epsilon_i},
& &D_2(\epsilon)=\sum_{i=1}^2\sum_{j=1}^2\epsilon_i\epsilon_jD_{\epsilon_i\epsilon_j},\\
&L_1(\epsilon)\phi=\sum_{i=1}^2\epsilon_iL_{\epsilon_i}\phi,
& &L_2(\epsilon)\phi=\sum_{i=1}^2\sum_{j=1}^2\epsilon_i\epsilon_jL_{\epsilon_i\epsilon_j}\phi,
\end{aligned}
\end{equation}
we have the following conclusion about normal form of Turing-Turing bifurcation for PFDEs \eqref{e:2.1}, of which the proof is put in Section 4.
\begin{theorem}\label{thm:2.1}
Assume that (H1)-(H4) hold. Then, the third-order normal form on center manifolds of PFDEs \eqref{e:2.1} at Turing-Turing singularity, reads
\begin{equation*}
\dot{z}=Bz+\frac{1}{2}g_2^1(z,0,\epsilon)+\frac{1}{3!}g_3^1(z,0,\epsilon)+h.o.t.,
\end{equation*}
where $h.o.t.$ stands for higher-order terms, and
\begin{equation*}
\frac{1}{j!}g_j^1(z,0,\epsilon)=\sum_{|q|+|s|=j} \frac{1}{\prod_{i=1}^2 q_i!\prod_{k=1}^2 s_k!}g^1_{qs}z^q\epsilon^s,\quad j=2,3,
\end{equation*}
with
\begin{equation*}
\begin{aligned}
&q=(q_1,q_2)\in\mathbb{N}_0^2,s=(s_1,s_2)\in\mathbb{N}_0^2,z=(z_1,z_2),\epsilon=(\epsilon_1,\epsilon_2),z^q=z_1^{q_1}z_2^{q_2}, \epsilon^s=\epsilon_1^{s_1}\epsilon_2^{s_2},\\
&|q|=q_1+q_2,|s|=s_1+s_2,g_{qs}\triangleq g_{q_1q_2s_1s_2}, g_{qs}=\left(g_{qs}^{1},g_{qs}^{2}\right)^T, g_{qs}^1=\left(g_{qs}^{11},g_{qs}^{12}\right)^T,
\end{aligned}
\end{equation*}
and
\begin{equation}\label{e:2.22}
\begin{aligned}
&g_{2000}^{11}=\psi_1(0)Q(\phi_1,\phi_1)\left\langle\beta_{k_1}^2,\beta_{k_1}\right\rangle;
&\qquad &g_{0200}^{11}=\psi_1(0)Q(\phi_2,\phi_2)\left\langle\beta_{k_2}^2,\beta_{k_1}\right\rangle;\\
&g_{1100}^{11}=\psi_1(0)Q(\phi_1,\phi_2)\left\langle\beta_{k_1}\beta_{k_2},\beta_{k_1}\right\rangle;
& &g_{2000}^{12}=\psi_2(0)Q(\phi_1,\phi_1)\left\langle\beta_{k_1}^2,\beta_{k_2}\right\rangle;\\
&g_{0200}^{12}=\psi_2(0)Q(\phi_2,\phi_2)\left\langle\beta_{k_2}^2,\beta_{k_2}\right\rangle;
& &g_{1100}^{12}=\psi_2(0)Q(\phi_1,\phi_2)\left\langle\beta_{k_1}\beta_{k_2},\beta_{k_2}\right\rangle;\\
&g_{1010}^{11}=\frac{1}{2}\psi_1(0)\left(L_{\epsilon_1}\phi_1-\mu_{k_1}D_{\epsilon_1}\phi_1(0)\right)\left\langle\beta_{k_1},\beta_{k_1}\right\rangle;
& &g_{1001}^{11}=\frac{1}{2}\psi_1(0)\left(L_{\epsilon_2}\phi_1-\mu_{k_1}D_{\epsilon_2}\phi_1(0)\right)\left\langle\beta_{k_1},\beta_{k_1}\right\rangle;\\
&g_{0110}^{11}=\frac{1}{2}\psi_1(0)\left(L_{\epsilon_1}\phi_2-\mu_{k_2}D_{\epsilon_1}\phi_2(0)\right)\left\langle\beta_{k_2},\beta_{k_1}\right\rangle;
& &g_{0101}^{11}=\frac{1}{2}\psi_1(0)\left(L_{\epsilon_2}\phi_2-\mu_{k_2}D_{\epsilon_2}\phi_2(0)\right)\left\langle\beta_{k_2},\beta_{k_1}\right\rangle;\\
&g_{1010}^{12}=\frac{1}{2}\psi_2(0)\left(L_{\epsilon_1}\phi_1-\mu_{k_1}D_{\epsilon_1}\phi_1(0)\right)\left\langle\beta_{k_1},\beta_{k_2}\right\rangle;
& &g_{1001}^{12}=\frac{1}{2}\psi_2(0)\left(L_{\epsilon_2}\phi_1-\mu_{k_1}D_{\epsilon_2}\phi_1(0)\right)\left\langle\beta_{k_1},\beta_{k_2}\right\rangle;\\
&g_{0110}^{12}=\frac{1}{2}\psi_2(0)\left(L_{\epsilon_1}\phi_2-\mu_{k_2}D_{\epsilon_1}\phi_2(0)\right)\left\langle\beta_{k_2},\beta_{k_2}\right\rangle;
& &g_{0101}^{12}=\frac{1}{2}\psi_2(0)\left(L_{\epsilon_2}\phi_2-\mu_{k_2}D_{\epsilon_2}\phi_2(0)\right)\left\langle\beta_{k_2},\beta_{k_2}\right\rangle;
\end{aligned}
\end{equation}
and
\begin{equation}\label{e:2.23}
\begin{split}
g_{3000}^{11}=&\psi_1(0)C(\phi_1,\phi_1,\phi_1)\left\langle\beta_{k_1}^3,\beta_{k_1}\right\rangle+3\psi_1(0)\left\langle Q(\phi_1,h_{2000})\beta_{k_1},\beta_{k_1}\right\rangle;\\
g_{3000}^{12}=&\psi_2(0)C(\phi_1,\phi_1,\phi_1)\left\langle\beta_{k_1}^3,\beta_{k_2}\right\rangle+3\psi_2(0)\left\langle Q(\phi_1,h_{2000})\beta_{k_1},\beta_{k_2}\right\rangle;\\
g_{0300}^{11}=&\psi_1(0)C(\phi_2,\phi_2,\phi_2)\left\langle\beta_{k_2}^3,\beta_{k_1}\right\rangle+3\psi_1(0)\left\langle Q(\phi_2,h_{0200})\beta_{k_2},\beta_{k_1}\right\rangle;\\
g_{0300}^{12}=&\psi_2(0)C(\phi_2,\phi_2,\phi_2)\left\langle\beta_{k_2}^3,\beta_{k_2}\right\rangle+3\psi_2(0)\left\langle Q(\phi_2,h_{0200})\beta_{k_2},\beta_{k_2}\right\rangle;\\
g_{2100}^{11}=&\psi_1(0)C(\phi_1,\phi_1,\phi_2)\left\langle\beta_{k_1}^2\beta_{k_2},\beta_{k_1}\right\rangle+2\psi_1(0)\left\langle Q(\phi_1,h_{1100})\beta_{k_1},\beta_{k_1}\right\rangle
+\psi_1(0)\left\langle Q(\phi_2,h_{2000})\beta_{k_2},\beta_{k_1}\right\rangle;\\
g_{2100}^{12}=&\psi_2(0)C(\phi_1,\phi_1,\phi_2)\left\langle\beta_{k_1}^2\beta_{k_2},\beta_{k_2}\right\rangle+2\psi_2(0)\left\langle Q(\phi_1,h_{1100})\beta_{k_1},\beta_{k_2}\right\rangle
+\psi_2(0)\left\langle Q(\phi_2,h_{2000})\beta_{k_2},\beta_{k_2}\right\rangle;\\
g_{1200}^{11}=&\psi_1(0)C(\phi_1,\phi_2,\phi_2)\left\langle\beta_{k_1}^2\beta_{k_2},\beta_{k_2}\right\rangle+\psi_1(0)\left\langle Q(\phi_1,h_{0200})\beta_{k_1},\beta_{k_1}\right\rangle+2\psi_1(0)\left\langle Q(\phi_2,h_{1100})\beta_{k_2},\beta_{k_1}\right\rangle;\\
g_{1200}^{12}=&\psi_2(0)C(\phi_1,\phi_2,\phi_2)\left\langle\beta_{k_1}\beta_{k_2}^2,\beta_{k_2}\right\rangle+\psi_2(0)\left\langle Q(\phi_1,h_{0200})\beta_{k_1},\beta_{k_2}\right\rangle
+2\psi_2(0)\left\langle Q(\phi_2,h_{1100})\beta_{k_2},\beta_{k_2}\right\rangle;
\end{split}
\end{equation}
where $\phi_j,\psi_j,j=1,2$ are given in \eqref{e:phi}, and $\left\{\beta_{k}\right\}_{{k}\in\mathbb{N}_*},\{\mu_{k}\}_{{k}\in\mathbb{N}_*}$ are eigenfunctions and eigenvalues of $-\Delta$ on $\Omega$ respectively, and $Q(\cdot,\cdot),C(\cdot,\cdot,\cdot)$ are defined by \eqref{e:2.20}, and $L_{\epsilon_j}(\cdot),D_{\epsilon_j},j=1,2$ are defined by \eqref{e:2.21}, and $h_{qs}\left(|q|+|s|=2\right)$ with $h_{qs}\triangleq h_{q_1q_2s_1s_2}$ are determined by Eq. \eqref{e:4.10}.
\end{theorem}
\begin{remark}
We omit quadratic coefficients $g_{qs}^1\;\left(|q|+|s|=2,|s|\ge2\right)$ and cubic coefficients $g_{qs}^1\;\left(|q|+|s|=3,|s|\ge1\right)$ involving in higher-order terms of perturbation parameter $\epsilon$, given that higher-order perturbations have few effects on bifurcation set of normal forms.
\end{remark}

Actually, formulas \eqref{e:2.22} and \eqref{e:2.23} could be further simplified, by choosing $\Omega\subset \mathbb{R}$ and Neumann boundary conditions. Then, three normal forms describing different spatial patterns are deduced from the third-order normal form in Theorem \ref{thm:2.1}. Specifically, let $\Omega=(0,l\pi),l>0$, then eigenvalue problem
\begin{equation*}
\beta''+\mu\beta=0,\quad \beta'(0)=\beta'(l\pi)=0,
\end{equation*}
has eigenvalues $\mu_{k}=\frac{k^2}{l^2}(k\in\mathbb{N}_0)$ and the corresponding normalized eigenfunctions
\begin{equation*}
 \begin{gathered}
 \beta_{ k}=\begin{cases}1,&k=0,\\
 \sqrt{2}\cos\frac{k}{l}x,&k\in \mathbb{N},\end{cases}
 \end{gathered}
\end{equation*}
with $\left\langle \beta_{ k},\beta_{ j} \right\rangle=\frac{1}{l\pi}\int_0^{l\pi}\beta_{ k}(x)\beta_{ j}(x) \mathrm{d}x, {k,j}\in\mathbb{N}_0$.

And, a direct calculations yields,
\begin{align*}
 \left\langle Q(\phi_i,h_{qs})\beta_{k_i},\beta_{k_j} \right\rangle =Q\left(\phi_i,\frac{1}{\sqrt{2}^{\delta(i-j)}}h^{|k_j-k_i|}_{qs}+\frac{1}{\sqrt{2}}h_{qs}^{|k_j+k_i|}\right),\quad i,j=1,2,
\end{align*}
where $h_{qs}^{ k}\triangleq\left\langle h_{qs},\beta_{ k} \right\rangle$, and
\begin{equation*}
\delta(i-j)=
\begin{cases}
1,\quad \mathrm{if}\; i\ne j,\\
0,\quad \mathrm{if}\; i= j.
\end{cases}
\end{equation*}

Then, there are three cases in total, depending on the relationship of $k_1$ and $k_2$. For convenience, suppose $k_2> k_1$. Firstly, we consider the case $k_2=2k_1$, and a direct calculation yields the following lemma.
\begin{lemma}\label{lems:2.3}
For $k_2=2k_1,k_1\in\mathbb{N}$ and $i,j=1,2$, it is not difficult to verity that,
\begin{equation*}
\begin{aligned}
\left\langle\beta_{k_i}, \beta_{k_j}\right\rangle=\begin{cases}1, &i=j,\\ 0, &i\neq j,\end{cases}\qquad
\left\langle\beta_{k_i}^m, \beta_{k_j}^n\right\rangle=\begin{cases}\frac{3}{2},&i=j,m+n=4,m,n\in\mathbb{N}_0,\\ 1, &i\neq j,m=n=2,\\
\frac{\sqrt{2}}{2}, &\frac{m}{k_j}=\frac{n}{k_i}, m+n=3,m,n\in\mathbb{N}_0,\\
0, &\mathrm{others}, 3\le m+n\le 4,m,n\in\mathbb{N}_0.\end{cases}
\end{aligned}
\end{equation*}
\end{lemma}

According to Lemma \ref{lems:2.3}, we have the following conclusion.
\begin{proposition} \label{pro:2.4}
For $k_2=2k_1,k_1\in\mathbb{N}$, on spatial domain $\Omega=(0,l\pi),l>0$ and under Neumann boundary conditions, the third-order normal form of Turing-Turing bifurcation reads
\begin{equation*}
\begin{cases}
\dot{z}_1=\left(g_{1010}^{11}\epsilon_1+g_{1001}^{11}\epsilon_2\right)z_1+g_{1100}^{11}z_1z_2+\frac{1}{6}g_{3000}^{11}z_1^3+\frac{1}{2}g_{1200}^{11}z_1z_2^2+h.o.t.,\\
\dot{z}_2=\left(g_{0110}^{12}\epsilon_1+g_{0101}^{12}\epsilon_2\right)z_2+\frac{1}{2}g_{2000}^{12}z_1^2+\frac{1}{2}g_{2100}^{12}z_1^2z_2+\frac{1}{6}g_{0300}^{12}z_2^3+h.o.t.,
\end{cases}
\end{equation*}
where $h.o.t.$ stands for higher-order terms, and
\begin{equation}\label{2.24}
\begin{aligned}
g_{1100}^{11}=&\frac{\sqrt{2}}{2}\psi_1(0)Q(\phi_1,\phi_2);\qquad\qquad\quad\;\;\;\,
g_{2000}^{12}=\frac{\sqrt{2}}{2}\psi_2(0)Q(\phi_1,\phi_1);\\
g_{1010}^{11}=&\frac{1}{2}\psi_1(0)\left(L_{\epsilon_1}\phi_1-\mu_{k_1}D_{\epsilon_1}\phi_1(0)\right);\quad\;\;\,
g_{1001}^{11}=\frac{1}{2}\psi_1(0)\left(L_{\epsilon_2}\phi_1-\mu_{k_1}D_{\epsilon_2}\phi_1(0)\right);\\
g_{0110}^{12}=&\frac{1}{2}\psi_2(0)\left(L_{\epsilon_1}\phi_2-\mu_{k_2}D_{\epsilon_1}\phi_2(0)\right);\quad\;\;\,
g_{0101}^{12}=\frac{1}{2}\psi_2(0)\left(L_{\epsilon_2}\phi_2-\mu_{k_2}D_{\epsilon_2}\phi_2(0)\right);\\
g_{3000}^{11}=&\frac{3}{2}\psi_1(0)C(\phi_1,\phi_1,\phi_1)+3\psi_1(0)Q\left(\phi_1,h_{2000}^{0}+\frac{1}{\sqrt{2}}h_{2000}^{2k_1}\right);\\
g_{1200}^{11}=&\psi_1(0)C(\phi_1,\phi_2,\phi_2)+\frac{2}{\sqrt{2}}\psi_1(0)Q\left(\phi_2,h^{k_2-k_1}_{1100}+h_{1100}^{k_2+k_1}\right)
+\psi_1(0)Q\left(\phi_1,h^{0}_{0200}\right);\\
g_{2100}^{12}=&\psi_2(0)C(\phi_1,\phi_1,\phi_2)+\frac{2}{\sqrt{2}}\psi_2(0)Q\left(\phi_1,h^{k_2-k_1}_{1100}+h_{1100}^{k_2+k_1}\right)
+\psi_2(0)Q\left(\phi_2,h^{0}_{2000}\right);\\
g_{0300}^{12}=&\frac{3}{2}\psi_2(0)C(\phi_2,\phi_2,\phi_2)+3\psi_2(0)Q\left(\phi_2,h^{0}_{0200}+\frac{1}{\sqrt{2}}h_{0200}^{2k_2}\right);
\end{aligned}
\end{equation}
with
\begin{equation}\label{e:2.25}
\begin{aligned}
h_{2000}^{0}(\theta)=&\Delta(0,0)^{-1}Q(\phi_1,\phi_1);
&h_{0200}^{0}(\theta)=&\Delta(0,0)^{-1}Q(\phi_2,\phi_2);\\
h_{1100}^{k_2+k_1}(\theta)=&\frac{\sqrt{2}}{2}\Delta\left(0,\mu_{k_2+k_1}\right)^{-1}Q(\phi_1,\phi_2);
&h_{1100}^{k_2-k_1}(\theta)=&h_{1100}^{k_1}(\theta)=\frac{\sqrt{2}}{2}\theta\phi_1\psi_1(0)Q(\phi_1,\phi_2)+h^{k_1}_{1100}(0);\\
h_{2000}^{2k_1}(\theta)=&\frac{\sqrt{2}}{2}\theta\phi_2\psi_2(0)Q(\phi_1,\phi_1)+h^{2k_1}_{2000}(0);
&h_{0200}^{2k_2}(\theta)=&\frac{\sqrt{2}}{2}\Delta\left(0,\mu_{2k_2}\right)^{-1}Q(\phi_2,\phi_2);
\end{aligned}
\end{equation}
where
\begin{equation}\label{e:2.26}
\begin{aligned}
\Delta\left(0,\mu_{2k_1}\right)h^{2k_1}_{2000}(0)=&\frac{\sqrt{2}}{2}\left(Q(\phi_1,\phi_1)-\phi_2(0)\psi_2(0)Q(\phi_1,\phi_1)\right);\\
\Delta\left(0,\mu_{k_1}\right)h^{k_1}_{1100}(0)=&\frac{\sqrt{2}}{2}\left(Q(\phi_1,\phi_2)-\phi_1(0)\psi_1(0)Q(\phi_1,\phi_2)\right);
\end{aligned}
\end{equation}
satisfy $\left(\psi_2, h^{2k_1}_{2000}(\theta)\right)_{k_2}=0$ and $\left(\psi_1, h^{k_1}_{1100}(\theta)\right)_{k_1}=0$.
\end{proposition}

Analogously, when $k_2=3k_1$ for $k_1\in\mathbb{N}$, we have the following lemma.
\begin{lemma}\label{lems:2.5}
For $k_2=3k_1,k_1\in\mathbb{N}$ and $i,j=1,2$, we have the following results,
\begin{equation*}
\begin{aligned}
\left\langle\beta_{k_i}, \beta_{k_j}\right\rangle=\begin{cases}1,&i=j,\\ 0, &i\neq j,\end{cases}\qquad
\left\langle\beta_{k_i}^m, \beta_{k_j}^n\right\rangle=\begin{cases}\frac{3}{2},&i=j,m+n=4,m,n\in\mathbb{N}_0,\\ 1, &i\neq j,m=n=2,\\
\frac{1}{2}, &\frac{m}{k_j}=\frac{n}{k_i}=3, m+n=4,m,n\in\mathbb{N}_0,\\
0, &\mathrm{others}, 3\le m+n\le 4,m,n\in\mathbb{N}_0.\end{cases}
\end{aligned}
\end{equation*}
\end{lemma}

Then based on Lemma \ref{lems:2.5}, the following conclusion holds.
\begin{proposition} \label{pro:2.6}
For $k_2=3k_1,k_1\in\mathbb{N}$, on spatial domain $\Omega=(0,l\pi),l>0$ and under Neumann boundary conditions, the third-order normal form of Turing-Turing bifurcation reads
\begin{equation*}
\begin{cases}
\dot{z}_1=\left(g_{1010}^{11}\epsilon_1+g_{1001}^{11}\epsilon_2\right)z_1+\frac{1}{6}g_{3000}^{11}z_1^3+\frac{1}{2}g_{1200}^{11}z_1z_2^2+\frac{1}{2}g_{2100}^{11}z_1^2z_2+h.o.t.,\\
\dot{z}_2=\left(g_{0110}^{12}\epsilon_1+g_{0101}^{12}\epsilon_2\right)z_2+\frac{1}{6}g_{3000}^{12}z_1^3+\frac{1}{2}g_{2100}^{12}z_1^2z_2+\frac{1}{6}g_{0300}^{12}z_2^3+h.o.t.,
\end{cases}
\end{equation*}
where $h.o.t.$ stands for higher-order terms, and
\begin{equation}\label{e:2.27}
\begin{aligned}
g_{1010}^{11}=&\frac{1}{2}\psi_1(0)\left(L_{\epsilon_1}\phi_1-\mu_{k_1}D_{\epsilon_1}\phi_1(0)\right);\quad\;\;\,
g_{1001}^{11}=\frac{1}{2}\psi_1(0)\left(L_{\epsilon_2}\phi_1-\mu_{k_1}D_{\epsilon_2}\phi_1(0)\right);\\
g_{0110}^{12}=&\frac{1}{2}\psi_2(0)\left(L_{\epsilon_1}\phi_2-\mu_{k_2}D_{\epsilon_1}\phi_2(0)\right);\quad\;\;\,
g_{0101}^{12}=\frac{1}{2}\psi_2(0)\left(L_{\epsilon_2}\phi_2-\mu_{k_2}D_{\epsilon_2}\phi_2(0)\right);\\
g_{3000}^{11}=&\frac{3}{2}\psi_1(0)C(\phi_1,\phi_1,\phi_1)+3\psi_1(0)Q\left(\phi_1,h^{0}_{2000}+\frac{1}{\sqrt{2}}h_{2000}^{2k_1}\right);\\
g_{1200}^{11}=&\psi_1(0)C(\phi_1,\phi_2,\phi_2)+\frac{2}{\sqrt{2}}\psi_1(0)Q\left(\phi_2,h^{k_2-k_1}_{1100}+h_{1100}^{k_2+k_1}\right)
+\psi_1(0)Q\left(\phi_1,h^{0}_{0200}\right);\\
g_{2100}^{11}=&\frac{1}{2}\psi_1(0)C(\phi_1,\phi_1,\phi_2)+\frac{2}{\sqrt{2}}\psi_1(0)Q\left(\phi_1,h_{1100}^{2k_1}\right)
+\frac{1}{\sqrt{2}}\psi_1(0)Q\left(\phi_1,h_{2000}^{k_2-k_1}\right);\\
g_{3000}^{12}=&\frac{1}{2}\psi_2(0)C(\phi_1,\phi_1,\phi_1)+\frac{3}{\sqrt{2}}\psi_2(0)Q\left(\phi_1,h^{k_2-k_1}_{2000}\right);\\
g_{2100}^{12}=&\psi_2(0)C(\phi_1,\phi_1,\phi_2)+\frac{2}{\sqrt{2}}\psi_2(0)Q\left(\phi_1,h^{k_2-k_1}_{1100}+h_{1100}^{k_2+k_1}\right)
+\psi_2(0)Q\left(\phi_2,h^{0}_{2000}\right);\\
g_{0300}^{12}=&\frac{3}{2}\psi_2(0)C(\phi_2,\phi_2,\phi_2)+3\psi_2(0)Q\left(\phi_2,h^{0}_{0200}+\frac{1}{\sqrt{2}}h_{0200}^{2k_2}\right);
\end{aligned}
\end{equation}
where
\begin{equation}\label{e:2.28}
\begin{aligned}
h_{2000}^{0}(\theta)=&\Delta(0,0)^{-1}Q(\phi_1,\phi_1);
&\qquad h_{0200}^{0}(\theta)=&\Delta(0,0)^{-1}Q(\phi_2,\phi_2);\\
h_{1100}^{k_2-k_1}(\theta)=&h_{1100}^{2k_1}(\theta)=\frac{\sqrt{2}}{2}\Delta\left(0,\mu_{2k_1}\right)^{-1}Q(\phi_1,\phi_2);
&h_{1100}^{k_2+k_1}(\theta)=&\frac{\sqrt{2}}{2}\Delta\left(0,\mu_{k_2+k_1}\right)^{-1}Q(\phi_1,\phi_2);\\
h_{2000}^{k_2-k_1}(\theta)=&h_{2000}^{2k_1}(\theta)=\frac{\sqrt{2}}{2}\Delta\left(0,\mu_{2k_1}\right)^{-1}Q(\phi_1,\phi_1);
&h_{0200}^{2k_2}(\theta)=&\frac{\sqrt{2}}{2}\Delta\left(0,\mu_{2k_2}\right)^{-1}Q(\phi_2,\phi_2).
\end{aligned}
\end{equation}
\end{proposition}

For $k_2>k_1, k_1,k_2\in \mathbb{N}$ and $k_2\ne 3k_1,k_2\ne 2k_1$, we also have
\begin{lemma}\label{lems:2.7}
For $k_2>k_1, k_1,k_2\in \mathbb{N}$ satisfying that $k_2\ne 3k_1,k_2\ne 2k_1$, a direct calculation yields
\begin{equation*}
\begin{aligned}
\left\langle\beta_{k_i}, \beta_{k_j}\right\rangle=\begin{cases}1,&i=j,\\ 0, &i\neq j,\end{cases}\qquad
\left\langle\beta_{k_i}^m, \beta_{k_j}^n\right\rangle=\begin{cases}\frac{3}{2},&i=j,m+n=4,m,n\in\mathbb{N}_0,\\ 1, &i\neq j,m=n=2,\\0, &\mathrm{others}, 3\le m+n\le4,m,n\in\mathbb{N}_0,\end{cases}
\end{aligned}
\end{equation*}
where $i,j=1,2$.
\end{lemma}

By Lemma \ref{lems:2.7}, we obtain
\begin{proposition} \label{pro:2.8}
For $k_2>k_1,k_1,k_2\in \mathbb{N}$, $k_2\ne 3k_1$ and $k_2\ne 2k_1$, on spatial domain $\Omega=(0,l\pi),l>0$ and under Neumann boundary conditions, the third-order normal form of Turing-Turing bifurcation reads
\begin{equation*}
\begin{cases}
\dot{z}_1=\left(g_{1010}^{11}\epsilon_1+g_{1001}^{11}\epsilon_2\right)z_1+\frac{1}{6}g_{3000}^{11}z_1^3+\frac{1}{2}g_{1200}^{11}z_1z_2^2+h.o.t.,\\
\dot{z}_2=\left(g_{0110}^{12}\epsilon_1+g_{0101}^{12}\epsilon_2\right)z_2+\frac{1}{2}g_{2100}^{12}z_1^2z_2+\frac{1}{6}g_{0300}^{12}z_2^3+h.o.t.,
\end{cases}
\end{equation*}
where $h.o.t.$ stands for higher-order terms, and
\begin{equation}\label{e:2.29}
\begin{aligned}
g_{1010}^{11}=&\frac{1}{2}\psi_1(0)\left(L_{\epsilon_1}\phi_1-\mu_{k_1}D_{\epsilon_1}\phi_1(0)\right);\quad\;\;\,
g_{1001}^{11}=\frac{1}{2}\psi_1(0)\left(L_{\epsilon_2}\phi_1-\mu_{k_1}D_{\epsilon_2}\phi_1(0)\right);\\
g_{0110}^{12}=&\frac{1}{2}\psi_2(0)\left(L_{\epsilon_1}\phi_2-\mu_{k_2}D_{\epsilon_1}\phi_2(0)\right);\quad\;\;\,
g_{0101}^{12}=\frac{1}{2}\psi_2(0)\left(L_{\epsilon_2}\phi_2-\mu_{k_2}D_{\epsilon_2}\phi_2(0)\right);\\
g_{3000}^{11}=&\frac{3}{2}\psi_1(0)C(\phi_1,\phi_1,\phi_1)+3\psi_1(0)Q\left(\phi_1,h^{0}_{2000}+\frac{1}{\sqrt{2}}h_{2000}^{2k_1}\right);\\
g_{1200}^{11}=&\psi_1(0)C(\phi_1,\phi_2,\phi_2)+\frac{2}{\sqrt{2}}\psi_1(0)Q\left(\phi_2,h^{k_2-k_1}_{1100}+h_{1100}^{k_2+k_1}\right)
+\psi_1(0)Q\left(\phi_1,h^{0}_{0200}\right);\\
g_{2100}^{12}=&\psi_2(0)C(\phi_1,\phi_1,\phi_2)+\frac{2}{\sqrt{2}}\psi_2(0)Q\left(\phi_1,h^{k_2-k_1}_{1100}+h_{1100}^{k_2+k_1}\right)
+\psi_2(0)Q\left(\phi_2,h^{0}_{2000}\right);\\
g_{0300}^{12}=&\frac{3}{2}\psi_2(0)C(\phi_2,\phi_2,\phi_2)+3\psi_2(0)Q\left(\phi_2,h^{0}_{0200}+\frac{1}{\sqrt{2}}h_{0200}^{2k_2}\right);
\end{aligned}
\end{equation}
where
\begin{equation}\label{e:2.30}
\begin{aligned}
h_{2000}^{0}(\theta)=&\Delta(0,0)^{-1}Q(\phi_1,\phi_1);
&h_{0200}^{0}(\theta)=&\Delta(0,0)^{-1}Q(\phi_2,\phi_2);\\
h_{1100}^{k_2-k_1}(\theta)=&\frac{\sqrt{2}}{2}\Delta\left(0,\mu_{k_2-k_1}\right)^{-1}Q(\phi_1,\phi_2);
&\qquad h_{1100}^{k_2+k_1}(\theta)=&\frac{\sqrt{2}}{2}\Delta\left(0,\mu_{k_2+k_1}\right)^{-1}Q(\phi_1,\phi_2);\\
h_{2000}^{2k_1}(\theta)=&\frac{\sqrt{2}}{2}\Delta\left(0,\mu_{2k_1}\right)^{-1}Q(\phi_1,\phi_1);
&h_{0200}^{2k_2}(\theta)=&\frac{\sqrt{2}}{2}\Delta\left(0,\mu_{2k_2}\right)^{-1}Q(\phi_2,\phi_2).
\end{aligned}
\end{equation}
\end{proposition}

We emphasize that Propositions \ref{pro:2.4}, \ref{pro:2.6} and \ref{pro:2.8} also apply to computing coefficients of normal forms for partial differential equations.
Moreover, the process of computing coefficients of normal forms for Turing-Turing bifurcation utilizing these concise formulas could be implemented by computer programs.

\section{Application to a diffusive predator-prey model}
In this section, we reveal superposition patterns resulting from Turing-Turing bifurcation via investigating dynamics of a diffusive predator-prey model, with the aid of these three normal forms.

In 1973, May \cite{MR} proposed the following model,
\begin{equation}
\left\{
\begin{aligned}\label{e:b1}
&\frac{\mathrm{d} u}{\mathrm{d} t}=r u\left(1-\frac{u}{K}\right)-\frac{muv}{1+au},\\
&\frac{\mathrm{d} v}{\mathrm{d} t}=sv\left(1-\frac{v}{hu}\right),
\end{aligned}
\right.
\end{equation}
with saturating predator functional response being of Holling type II, where $r$ and $s$ denote their intrinsic growth rates, respectively; $K$ is the carrying capacity of the prey's environment; $h$ is a measure of the food quality of the prey for conversion into predator growth depending on the density of the prey population; $m$ and $a$ describe the effects of capture rate and handling time, respectively.
Then, given that interference among predators exists in nature, Holling type II functional response is modified to Crowley-Martin functional response \cite{CPH&MEK}, and model \eqref{e:b1} is converted to
\begin{equation}
\left\{
\begin{aligned}\label{e:b2}
&\frac{\mathrm{d} u}{\mathrm{d} t}=r u\left(1-\frac{u}{K}\right)-\frac{muv}{(1+au)(1+bu)},\\
&\frac{\mathrm{d} v}{\mathrm{d} t}=sv\left(1-\frac{v}{hu}\right),
\end{aligned}
\right.
\end{equation}
where $b$ describes the magnitude of interference among predators. Especially, if $b=0$, Crowley-Martin functional response reduces to Holling type II functional response, then model \eqref{e:b2} also reduces to model \eqref{e:b1}. Thus, it is more reasonable to consider model \eqref{e:b2}, especially when the predator feeding rate is decreased by higher predator density even when prey density is high.
Furthermore, due to the fact that the distribution of species is generally spatially heterogeneous and therefore the species will migrate toward regions of lower population density to improve the possibility of survival, Shi and Ruan \cite{SHB&RSG} introduced diffusion terms into model \eqref{e:b2} and proposed the following diffusive predator-prey model subject to Neumann boundary conditions,
\begin{equation}
\left\{
\begin{aligned}\label{e:3.0}
&\frac{\partial u}{\partial t}-d_1\Delta u =ru\left(1-\frac{u}{K}\right)-\frac{muv}{(1+au)(1+bv)},& &\quad x\in \Omega,\, t>0,\\
&\frac{\partial v}{\partial t}-d_2\Delta v =sv\left(1-\frac{v}{hu}\right),& &\quad x\in \Omega,\, t>0,\\
&\frac{\partial u}{\partial \nu}=\frac{\partial v}{\partial \nu}=0,& &\quad x\in \partial\Omega,\, t>0,\\
&u(x,0)=u_0(x)\ge 0, \, v(x,0)=v_0(x)\ge 0,& &\quad x\in \Omega.
\end{aligned}\right.
\end{equation}
where $u(x,t)$ and $v(x,t)$ represent the population densities of the prey and predator at location $x \in \Omega $ and time $t\ge 0$, respectively. The domain $\Omega\subset \mathbb{R}^n$ with smooth boundary $\partial \Omega$ is bounded. And $\nu$ is the outward unit normal vector. Diffusion coefficients $d_1$ and $d_2$ are positive, and nonnegative and continuous functions $u_0(x)$ and $v_0(x)$ denote initial functions.

Actually, Hopf bifurcation and Turing instability of system \eqref{e:3.0}, have been discussed in \cite{SHB&RSG}. Also, conditions of the occurrences of Turing-Hopf bifurcation and \emph{Turing-Turing bifurcation}, and spatiotemporal dynamics of system \eqref{e:3.0} near Turing-Hopf singularity, have been investigated in \cite{CX&JWH}. Especially, we anticipated in \cite{CX&JWH} that under proper conditions, diffusive system \eqref{e:3.0} could exhibit complex spatial patterns that four stable spatially inhomogeneous steady states with different spatial wavelengths coexist. In the following, we compute normal forms of Turing-Turing bifurcation for diffusive predator-prey system \eqref{e:3.0} by utilizing Propositions \ref{pro:2.4}, \ref{pro:2.6} and \ref{pro:2.8}, and further investigate spatial patterns near Turing-Turing singularity by analyzing the obtained normal forms.

By applying the following scalings
\begin{equation*}
rt\mapsto t, \, \frac{u}{K}\mapsto u, \, \frac{mKh}{r}\mapsto m, \, aK\mapsto a, \, bKh\mapsto b, \, \frac{s}{r}\mapsto s, \, \frac{d_1}{r_1}\mapsto d_1, \, \frac{d_2}{r_1}\mapsto d_2,
\end{equation*}
diffusive system \eqref{e:3.0} is transformed into
\begin{equation}
\left\{
\begin{aligned}\label{e:3.1}
&\frac{\partial u}{\partial t}-d_1\Delta u =u\left(1-u\right)-\frac{muv}{(1+au)(1+bv)},& &\quad x\in \Omega,\, t>0,\\
&\frac{\partial v}{\partial t}-d_2\Delta v =sv\left(1-\frac{v}{u}\right),& &\quad x\in \Omega,\, t>0,\\
&\frac{\partial u}{\partial \nu}=\frac{\partial v}{\partial \nu}=0,& &\quad x\in \partial\Omega,\, t>0,\\
&u(x,0)=u_0(x)\ge 0, \, v(x,0)=v_0(x)\ge 0,& &\quad x\in \Omega.
\end{aligned}\right.
\end{equation}
Obviously, $E_1=(1,0)$ is a boundary equilibrium of system \eqref{e:3.1}.
Moreover, according to \cite{CX&JWH, SHB&RSG}, if $a+b \ge ab$, system \eqref{e:3.1} admits a unique interior equilibrium $E_*=(u_*,v_*)$, where $u_*>0,v_*>0$ satisfy
\begin{equation*}
1-u-\frac{mv}{(1+au)(1+bv)}=0,\quad \frac{v}{u}=1.
\end{equation*}

Let $\Omega=(0,\pi)$, then the linearized system of \eqref{e:3.1} at $E_*=(u_*,v_*)$ reads
\begin{equation}
\left\{
\begin{aligned}\label{e:3.2}
\frac{\partial u}{\partial t}-d_1\Delta u &=s_0 (u-u_*)+\sigma (v-v_*),\\
\frac{\partial v}{\partial t}-d_2\Delta v &=s(u-u_*)-s(v-v_*),
\end{aligned}\right.
\end{equation}
where
\begin{equation}\label{e:3.3}
s_0=u_*\left(\frac{amu_*}{(1+au_*)^2(1+bu_*)}-1\right), \quad \sigma=-\frac{mu_*}{(1+au_*)(1+bu_*)^2}<0.
\end{equation}

Thus, the characteristic equations of system \eqref{e:3.2} are
\begin{equation}\label{e:3.4}
P_k(\lambda)\triangleq\lambda^2-\Theta(k)\lambda+\Delta(k)=0,\quad k\in \mathbb{N}_0,
\end{equation}
where
\begin{align*}
\Theta(k)\triangleq s_0-s-(d_1+d_2)k^2, \quad
\Delta(k)\triangleq d_1d_2k^4+(sd_1-s_0d_2)k^2-s(s_0+\sigma).
\end{align*}

Then, we have the following conclusions about Turing-Turing bifurcation, due to \cite{CX&JWH}.
\begin{lemma}\label{lems:3.1}(See \cite{CX&JWH})
Assume that $a+b\ge ab$ and $s_0>0$. If $s_0+\sigma < 0$, system \eqref{e:3.1} undergoes codimension-2  $(i,j)$-mode Turing-Turing bifurcation at $E_*$ when $\left(d_1,s\right)=\left(d_{i,j}^*, s_{i,j}^*\right)$, where
\begin{equation*}
\begin{aligned}
d_{i,j}^*&=\frac{\left(i^2+j^2\right)(s_0+\sigma)+\sqrt{\left(i^2+j^2\right)^2(s_0+\sigma)^2-4i^2j^2(s_0+\sigma)s_0}}{2i^2j^2}>0,\\ s_{i,j}^*&=\frac{s_0d_2i^2-d_{i,j}^*d_2i^4}{d_{i,j}^*i^2-(s_0+\sigma)}>0,\quad i,j\in \mathbb{N}.
\end{aligned}
\end{equation*}
\end{lemma}
\begin{lemma}\label{lems:3.2}(See \cite{CX&JWH})
Assume that conditions in Lemma \ref{lems:3.1} hold. For $k> k_0^*$, $(k-1,k)$-mode Turing-Turing bifurcation points $\left(d_{k-1,k}^*, s_{k-1,k}^*\right)$ are locating on boundary of parameter region of the stability for coexistence equilibrium, satisfying that characteristic Eq. \eqref{e:3.4} have two independent zero roots with the remaining roots having negative real parts, where
\begin{equation*}
k^*_0=\max\left\{ \xi:d_{\xi}^*=\max_{k\in \mathbb{N}}d_k^*\right\},\qquad\mathrm{with}\; d_k^*=\frac{s_0}{k^2}\left(1+\frac{\sigma}{d_2k^2+s_0}\right),~k\in\mathbb{N}.
\end{equation*}
\end{lemma}

\subsection{Normal form of Turing-Turing bifurcation}
Then, we compute normal form of system \eqref{e:3.1} at Turing-Turing singularity.
Introducing perturbation parameters $\epsilon_1,\epsilon_2$ by letting $d_1=d_{k_0-1,k_0}^*+\epsilon_1,s=s_{k_0-1,k_0}^*+\epsilon_2$ for $k_0>k_0^*$, which satisfy that system \eqref{e:3.1} undergoes $\left(k_0-1,k_0\right)$-mode Turing-Turing bifurcation at the coexistence equilibrium $E_*$ when $\epsilon_1=0,\epsilon_2=0$, and transforming $E_*$ to the origin, system \eqref{e:3.1} becomes
\begin{equation}
\begin{cases}\label{e:3.5}
\ds \frac{\partial u}{\partial t}-\left(d_{k_0-1,k_0}^*+\epsilon_1\right)\Delta u =(u+u_*)\left(1-(u+u_*)\right)-\frac{m(u+u_*)(v+v_*)}{(1+a(u+u_*))(1+b(v+v_*))},\\
\ds \frac{\partial v}{\partial t}-d_2\Delta v =\left(s_{k_0-1,k_0}^*+\epsilon_2\right)(v+v_*)\left(1-\frac{v+v_*}{u+u_*}\right).
\end{cases}
\end{equation}
Defining $U(t) = (u(t),v(t))^T$, system \eqref{e:3.5} is written as an abstract differential equation
\begin{equation*}
\frac{\mathrm{d}U(t)}{\mathrm{d}t}=D(\epsilon)\Delta U(t)+L(\epsilon)U(t)+F(U(t),\epsilon),
\end{equation*}
in the phase space $\mathscr{X}$,
\begin{equation*}
\mathscr{X}\triangleq \left\{ (u,v)^T\in \left(W^{2,2}(0,\pi)\right)^2|\frac{\partial u}{\partial x}|_{x=0,\pi}=\frac{\partial v}{\partial x}|_{x=0,\pi}=0\right\},
\end{equation*}
where
\begin{equation*}
\begin{split}
D(\epsilon)&=\begin{pmatrix}d_{k_0-1,k_0}^*+\epsilon_1&0\\0&d_2\end{pmatrix},\qquad\qquad
L(\epsilon)=\begin{pmatrix}s_0&\sigma\\s_{k_0-1,k_0}^*+\epsilon_2&-\left(s_{k_0-1,k_0}^*+\epsilon_2\right)\end{pmatrix},\\
F(\varphi,\epsilon)&=\begin{pmatrix}(\varphi_1+u_*)\left(1-(\varphi_1+u_*)\right)-\frac{m(\varphi_1+u_*)(\varphi_2+v_*)}{(1+a(\varphi_1+u_*))(1+b(\varphi_2+v_*))}
-s_0\varphi_1-\sigma \varphi_2\\
\left(s_{k_0-1,k_0}^*+\epsilon_2\right)(\varphi_2+v_*)\left(1-\frac{(\varphi_2+v_*)}{(\varphi_1+u_*)}\right)-\left(s_{k_0-1,k_0}^*+\epsilon_2\right)\varphi_1+\left(s_{k_0-1,k_0}^*+\epsilon_2\right)\varphi_2\end{pmatrix},
\end{split}
\end{equation*}
with $\varphi=(\varphi_1,\varphi_2)^T \in \mathscr{X}$.

According to expansions \eqref{e:2.20}, we have
\begin{equation*}
\begin{split}
D(0)&=\begin{pmatrix}d_{k_0-1,k_0}^*&0\\ 0&d_2\end{pmatrix},
D_1(\epsilon)=\begin{pmatrix}\epsilon_1&0\\ 0&0\end{pmatrix},
L(0)=\begin{pmatrix}s_0&\sigma\\ s_{k_0-1,k_0}^*&-s_{k_0-1,k_0}^*\end{pmatrix},
L_1(\epsilon)=\begin{pmatrix}0&0\\ \epsilon_2&-\epsilon_2\end{pmatrix},\\
Q(\varphi,\varphi)&=\begin{pmatrix}\left(-2+\frac{2mv_*a}{(au_*+1)^2(bv_*+1)}-\frac{2mu_*v_*a^2}{(au_*+1)^3(bv_*+1)}\right)\varphi_1^2-\frac{2m}{(au_*+1)^2(bv_*+1)^2}\varphi_1\varphi_2
+\frac{2mu_*b}{(au_*+1)(bv_*+1)^3}\varphi_2^2\\ \frac{2s_{k_0-1,k_0}^*v_*}{u_*^2}\varphi_1^2-\frac{4s_{k_0-1,k_0}^*}{u_*}\varphi_1\varphi_2-\frac{2s_{k_0-1,k_0}^*v_*^2}{u_*^3}\varphi_2^2\end{pmatrix},\\
C(\varphi,\varphi, \varphi)&=\begin{pmatrix}-\frac{6mv_*a^2}{(au_*+1)^4(bv_*+1)}\varphi_1^3-\frac{6ma}{(au_*+1)^3(bv_*+1)^2}\varphi_1^2\varphi_2
+\frac{6mb}{(au_*+1)^2(bv_*+1)^3}\varphi_1\varphi_2^2-\frac{6mu_*b^2}{(au_*+1)(bv_*+1)^4}\varphi_2^3\\ \frac{6s_{k_0-1,k_0}^*v_*^2}{u_*^4}\varphi_1^3-\frac{12s_{k_0-1,k_0}^*v_*}{u_*^3}\varphi_1^2\varphi_2-\frac{6s_{k_0-1,k_0}^*}{u_*^2}\varphi_1\varphi_2^2\end{pmatrix},
\end{split}
\end{equation*}
with $\varphi=(\varphi_1,\varphi_2)^T \in \mathscr{X}$.

And, the corresponding characteristic matrices are
\begin{equation*}
\Delta(\lambda,\mu_k)=\lambda I-\tilde{\Delta}_k, \quad \mathrm{with}\;\tilde{\Delta}_k\triangleq\begin{pmatrix}-d_{k_0-1,k_0}^*k^2+s_0&\sigma\\ s_{k_0-1,k_0}^*&-d_2k^2-s_{k_0-1,k_0}^*\end{pmatrix},\quad k\in \mathbb{N}_0.
\end{equation*}

Obviously, $\tilde{\Delta}_{k}$ has two independent zero eigenvalues corresponding to $k_0-1$ and $k_0$ respectively, such that the remaining eigenvalues of $\tilde{\Delta}_{k}$ have negative real parts, according to Lemma \ref{lems:3.2}. Then, a direct calculation yields
\begin{equation}\label{e:3.6-1}
\begin{aligned}
&\phi_1=\begin{pmatrix}1 \\ \frac{s_{k_0-1,k_0}^*}{d_2{(k_0-1)}^2+s_{k_0-1,k_0}^*}\end{pmatrix},
& &\psi_1=\begin{pmatrix}\frac{1}{N_1} \\ \frac{d_{k_0-1,k_0}^*{(k_0-1)}^2-s_0}{s_{k_0-1,k_0}^*N_1}\end{pmatrix}^T,
& &N_1=1+\frac{d_{k_0-1,k_0}^*{(k_0-1)}^2-s_0}{d_2{(k_0-1)}^2+s_{k_0-1,k_0}^*},\\
&\phi_2=\begin{pmatrix}1 \\ \frac{s_{k_0-1,k_0}^*}{d_2k_0^2+s_{k_0-1,k_0}^*}\end{pmatrix},
& &\psi_2=\begin{pmatrix}\frac{1}{N_2} \\ \frac{d_{k_0-1,k_0}^*k_0^2-s_0}{s_{k_0-1,k_0}^*N_2}\end{pmatrix}^T,
& &N_2=1+\frac{d_{k_0-1,k_0}^*k_0^2-s_0}{d_2k_0^2+s_{k_0-1,k_0}^*}.
\end{aligned}
\end{equation}
And, $\Phi=\left(\phi_1,\phi_2\right)$ and $\Psi=\left(\psi_1,\psi_2\right)^T$ satisfy that $\psi_j\phi_j=1$ for $j=1,2$.

Then by Theorem \ref{thm:2.1}, the third-order normal form on center manifolds for system \eqref{e:3.1} at Turing-Turing singularity, has the form
\begin{equation}\label{e:3.6}
\dot{z}=Bz+\frac{1}{2}g_2^1(z,0,\epsilon)+\frac{1}{3!}g_3^1(z,0,\epsilon)+h.o.t.,
\end{equation}
with
\begin{equation*}
\begin{aligned}
\frac{1}{j!}g_j^1(z,0,\epsilon)=\sum_{|q|+|s|=j} \frac{1}{\prod_{i=1}^2 q_i!\prod_{k=1}^2 s_k!}g^1_{qs}z^q\epsilon^s, \quad j=2,3,
\end{aligned}
\end{equation*}
where $h.o.t.$ stands for higher-order terms, and polynomial coefficients $\{g_{qs}^1:|q|+|s|=2,3\}$ could be computed by utilizing these formulas given in Proposition \ref{pro:2.4}, Proposition \ref{pro:2.6} and Proposition \ref{pro:2.8}, depending on the relationship of $(k_0-1)$ and $k_0$.

\subsection{Spatial patterns arising from Turing-Turing bifurcation}
We consider two sets of different parameters for system \eqref{e:3.1}, which involve in different normal forms and different spatial patterns.

\begin{figure}[!hbt]
\centering
\subfigure[Parameter region of the stability for coexistence equilibrium $E_*$ and bifurcation set in $d_1\text{-}s$ plane. $TT_{2}$ represent $(2,3)-$mode Turing-Turing bifurcation point.]{
\label{fig:1:a}
\begin{minipage}[b]{0.485\textwidth}
  \centering\includegraphics[width=\textwidth]{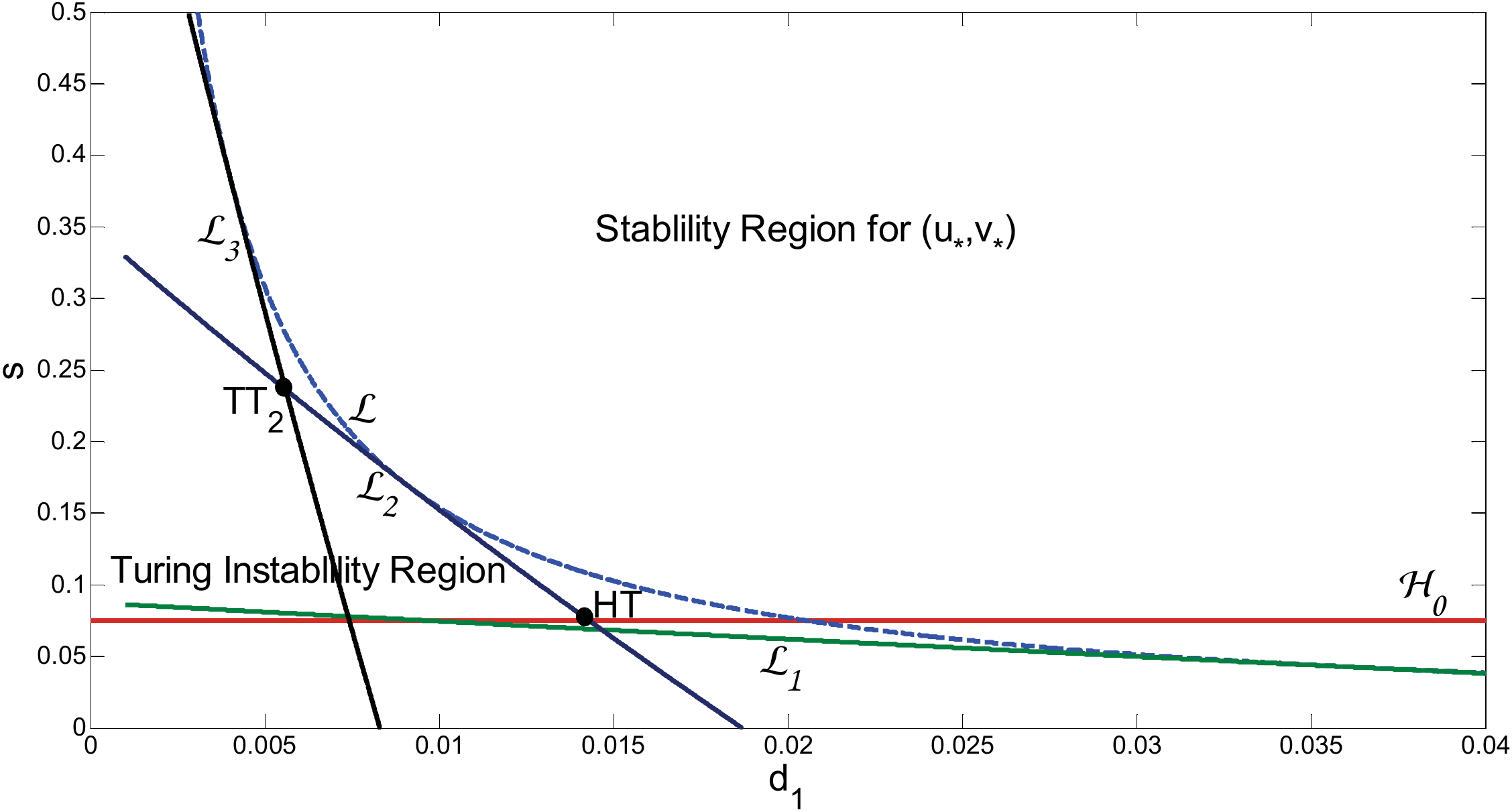}
\end{minipage}
}
\subfigure[Bifurcation set with local bifurcation curves $\mathcal{T}_1,\mathcal{T}_2$ at Turing-Turing point $\left(d_{2,3}^*,s_{2,3}^*\right)$.]{
\label{fig:1:b}
\begin{minipage}[b]{0.475\textwidth}
\centering \includegraphics[width=\textwidth]{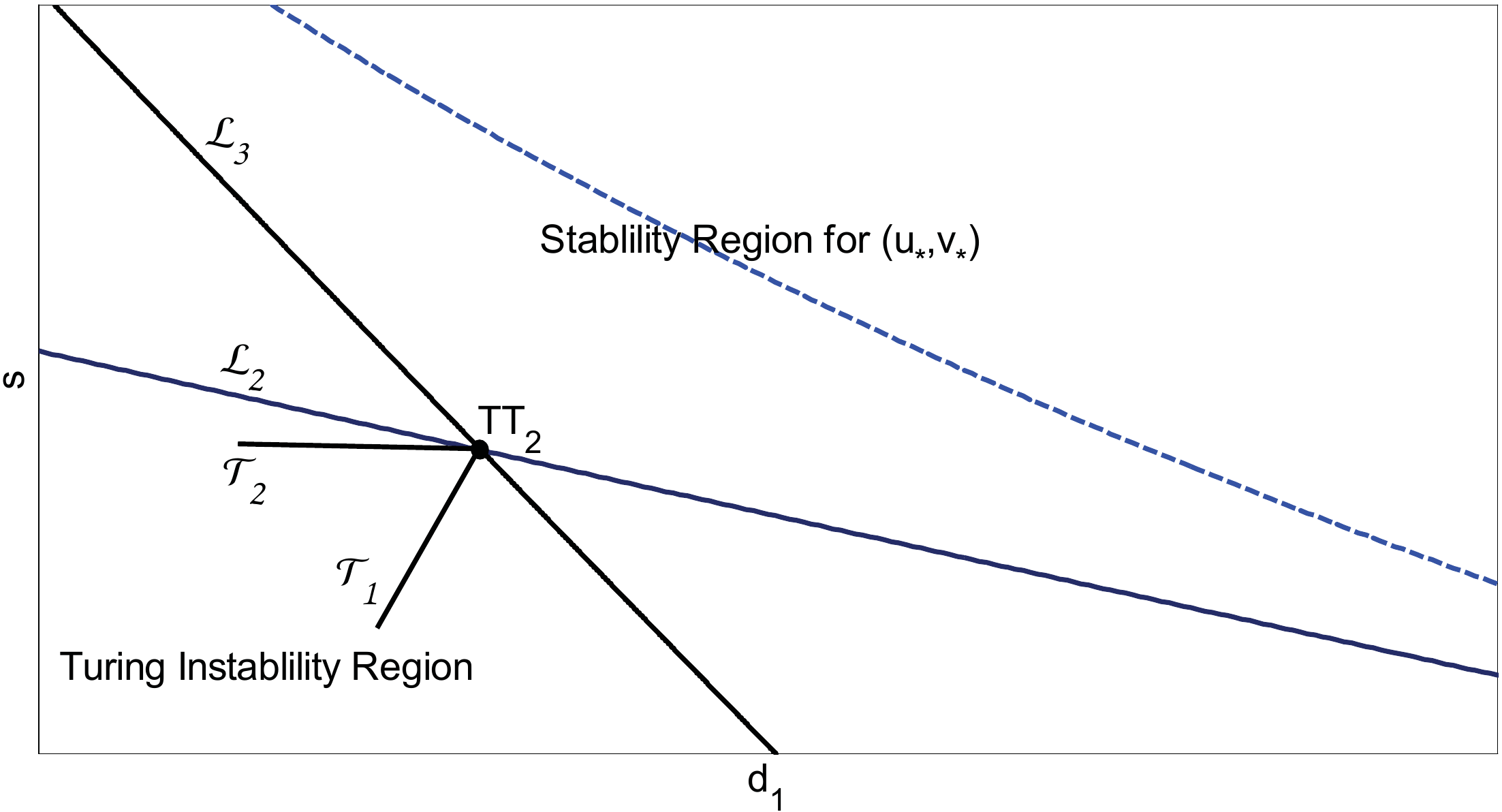}
\end{minipage}
}
\caption{\rm Turing bifurcation curves and bifurcation set (a), and bifurcation set with local bifurcation curves (b) in $d_1\text{-}s$ plane.}
  \label{fig:1}
\end{figure}

For the first set of system parameters, let $m=6,a=3,b=0.5,d_2=0.7$, then $(u_*,v_*)=(0.245,0.245)$. We also derive $s_0=0.0748, \sigma=-0.673$ satisfying $s_0+\sigma<0$. Then, a direct calculation yields $k_0^*=2$, and
\begin{equation*}
\begin{split}
\mathcal{L}_k:s=s_k(d_1)=\frac{s_0d_2k^2-d_1d_2k^4}{d_1k^2-(s_0+\sigma)}=\frac{-0.7d_1 k^4+0.0524k^2}{d_1 k^2+0.598},\quad 0< d_1<\frac{s_0}{k^2},k\in\mathbb{N}.
\end{split}
\end{equation*}
Therefore, Turing bifurcation curves in $d_1\text{-}s$ plane are shown in Fig. \ref{fig:1:a}, see \cite{CX&JWH}.

Choose ${k_1}=2,{k_2}=3$ in Theorem \ref{thm:2.1}, then $\left(d_{2,3}^*,s_{2,3}^*\right)=(0.0056,0.2364)$ is the \emph{$(2,3)-$mode Turing-Turing bifurcation point}.

Applying \emph{Proposition \ref{pro:2.8}}, for these given parameters, normal form truncated to order 3, reads
\begin{equation}\label{e:3.7}
\left\{
\begin{aligned}
\dot{z}_1&=-(4.0702\epsilon_1+0.20782\epsilon_2)z_1-1.6069z_1^3-5.9052z_1z_2^2,\\
\dot{z}_2&=-(9.0336\epsilon_1+0.099536\epsilon_2)z_2-4.1872z_1^2z_2-3.2439z_2^3.
\end{aligned}\right.
\end{equation}
Then, \eqref{e:3.7} has the following equilibria,
\begin{equation*}
\begin{aligned}
&A_0=(0,0),\;\text{for all}\; \epsilon_1,\epsilon_2,\\
&A_1^{\pm}=(\pm\sqrt{-2.5330\epsilon_1-0.1293\epsilon_2},0), \;\text{for}\;\epsilon_2<-19.5852\epsilon_1,\\
&A_2^{\pm}=(0,\pm\sqrt{-2.7848\epsilon_1-0.03068\epsilon_2}), \;\text{for}\; \epsilon_2<-90.7571\epsilon_1,\\ &A_{3}^{\pm\pm}=(\pm h_1,\pm h_2),\;\text{with}\;h_1=\sqrt{-2.0571\epsilon_1+0.004426\epsilon_2}, h_2=\sqrt{-0.1295\epsilon_1-0.03640\epsilon_2}\\
&\qquad\qquad\qquad\qquad\quad \text{for}\; \epsilon_2<-3.5568\epsilon_1,\epsilon_1<0,\;\text{and}\;\epsilon_2>464.7488\epsilon_1,\epsilon_1<0.
\end{aligned}
\end{equation*}
Based on \emph{center manifold theory} \cite{CJ,WJH}, correspondences between equilibria of normal form \eqref{e:3.7} and steady states of diffusive system \eqref{e:3.1} are shown in Table 1. Then, we could reveal long-time dynamical behaviors of system \eqref{e:3.1}, by analyzing long-time dynamics of normal form \eqref{e:3.7}.
\begin{table}[ht]\label{tab:1}
\centering\caption{Correspondences between equilibria of normal form \eqref{e:3.7} and steady states of diffusive system \eqref{e:3.1}.}
{\begin{tabular}{ll}
\hline\noalign{\smallskip}
Equilibria of normal form \eqref{e:3.7} & Steady states of diffusive system \eqref{e:3.1}\\
\hline\\
\qquad\qquad\qquad $A_0$ & Constant steady state $E_*$\\
\qquad\qquad\qquad $A_1^\pm$ & Non-constant steady states with the shape of $\phi_1\cos2x-$like\\
\qquad\qquad\qquad $A_2^\pm$ & Non-constant steady states with the shape of $\phi_2\cos3x-$like\\
\qquad\qquad\qquad $A_3^{\pm\pm}$ & Non-constant steady states with the shape of $h_1\phi_1\cos2x+h_2\phi_2\cos3x-$like\\
\noalign{\smallskip}\hline
\end{tabular}}
\end{table}

\begin{figure}[!hbt]
\centering
\subfigure[Local bifurcation set for normal form \eqref{e:3.7}.]{
\label{fig:2:a}
\begin{minipage}[b]{0.42\textwidth}
  \centering \includegraphics[width=\textwidth]{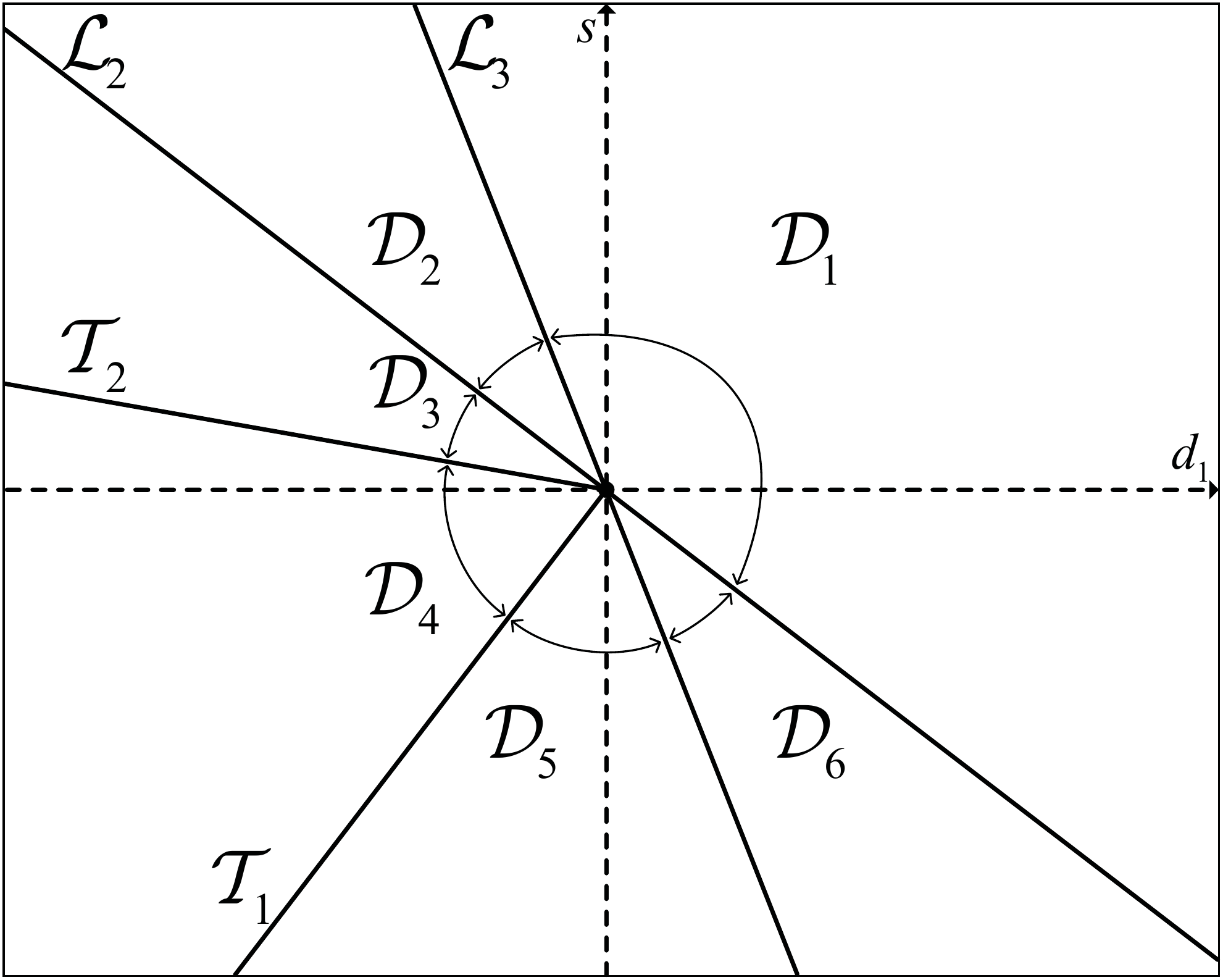}
\end{minipage}
}
\subfigure[Phase portraits for normal form \eqref{e:3.7}.]{
\label{fig:2:b}
\begin{minipage}[b]{0.535\textwidth}
\centering \includegraphics[width=\textwidth]{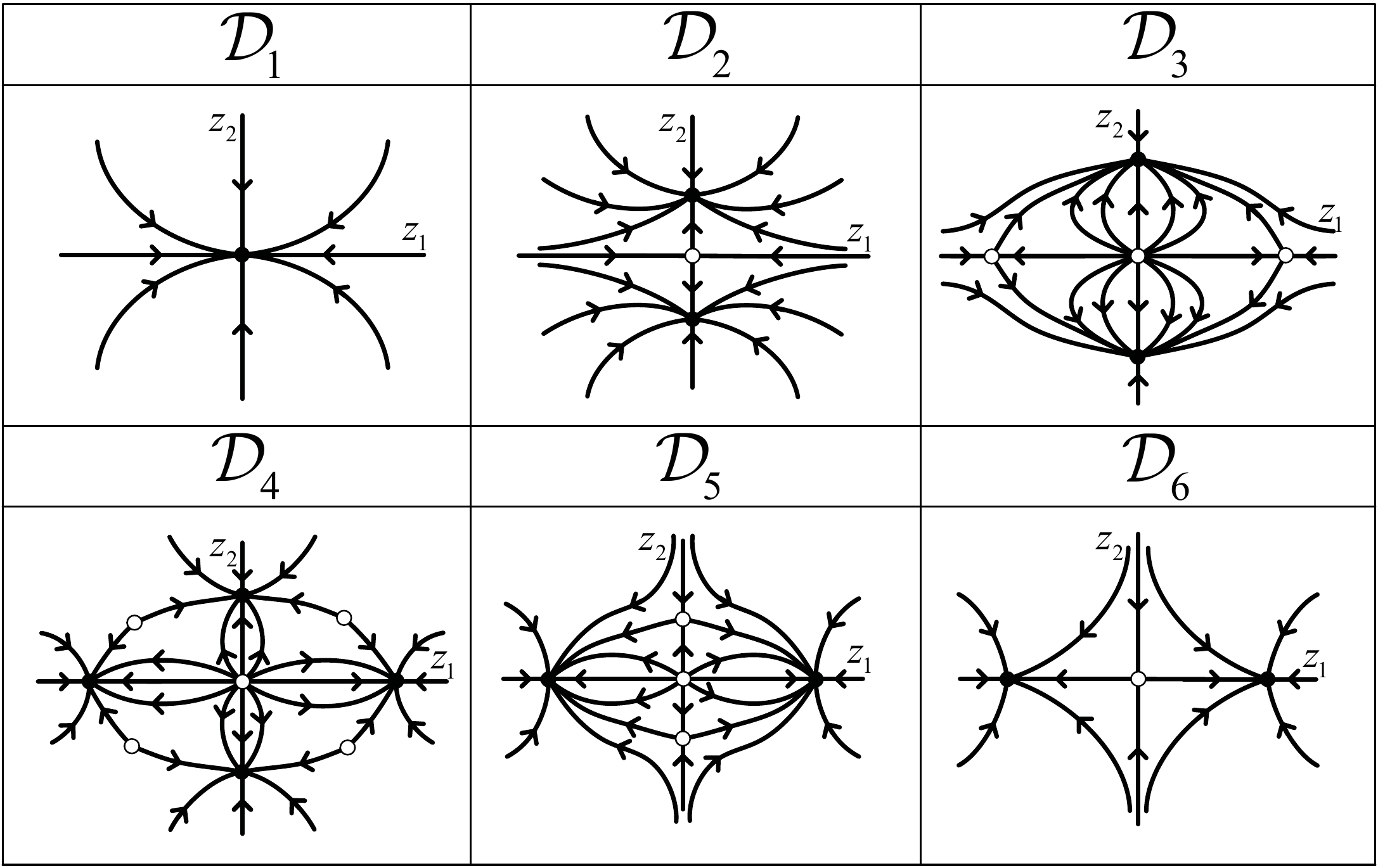}
\end{minipage}
}
\caption{\rm Local bifurcation set (a) at $(2,3)-$mode Turing-Turing point $\left(d_{2,3}^*,s_{2,3}^*\right)$ and the corresponding phase portraits (b).}
\label{fig:2}
\end{figure}

According to Table 7.5.2 in \cite{GJ&HP}, there are twelve unfoldings (see Table 2)
\begin{table}[ht]\label{tab:2}
\centering\caption{Twelve unfoldings of normal form \eqref{e:nm}, according to \cite{GJ&HP}.}
{\begin{tabular}{l l l l l l l l l l l l l}
\hline\noalign{\smallskip}
Case& Ia& Ib& II& III& IVa& IVb& V& VIa& VIb& VIIa& VIIb& VIII\\
\hline\\
$d_0$& $+1$& $+1$& $+1$& $+1$& $+1$& $+1$& $-1$& $-1$& $-1$& $-1$& $-1$& $-1$\\
$b_0$& $+$& $+$& $+$& $-$& $-$& $-$& $+$& $+$& $+$& $-$& $-$& $-$\\
$c_0$& $+$& $+$& $-$& $+$& $-$& $-$& $+$& $-$& $-$& $+$& $+$& $-$\\
$d_0-b_0c_0$& $+$& $-$& $+$& $+$& $+$& $-$& $-$& $+$& $-$& $+$& $-$& $-$\\
\noalign{\smallskip}\hline
\end{tabular}}
\end{table}
for the following normal form
\begin{equation}\label{e:nm}
\begin{cases}
\dot{z}_1=z_1\left(\alpha(\epsilon)+z_1^2+b_0z_2^2\right),\\
\dot{z}_2=z_2\left(\beta(\epsilon)+c_0z_1^2+d_0z_2^2\right).
\end{cases}
\end{equation}
A direct calculation yields that the unfolding for system \eqref{e:3.7} is Case Ib. Therefore, define the following critical bifurcation curves in $d_1\text{-}s$ plane,
\begin{equation*}
\begin{aligned}
&\mathcal{L}_2: s=s^*_{2,3}-19.5852\left(d_1-d_{2,3}^*\right), & &\mathcal{L}_3: s=s^*_{2,3}-90.7571\left(d_1-d_{2,3}^*\right),\\
&\mathcal{T}_1: s=s^*_{2,3}+464.7488\left(d_1-d_{2,3}^*\right), d_1\le d_{2,3}^*, &\quad &\mathcal{T}_2: s=s^*_{2,3}-3.5568\left(d_1-d_{2,3}^*\right), d_1\le d_{2,3}^*.
\end{aligned}
\end{equation*}
And, local bifurcation set and the corresponding phase portraits are shown in Fig. \ref{fig:2}.

Actually, local bifurcation set could be embedded into global bifurcation set, see Fig. \ref{fig:1:b}.

\begin{figure}[!hbt]
\centering
\subfigure[The initial values are $u(0,x)=0.245-0.02\cos 2x, v(0,x)=0.245+0.05\cos 2x$.]{
\label{fig:61:a} 
\begin{minipage}[b]{0.48\textwidth}
  \centering\includegraphics[width=0.49\textwidth]{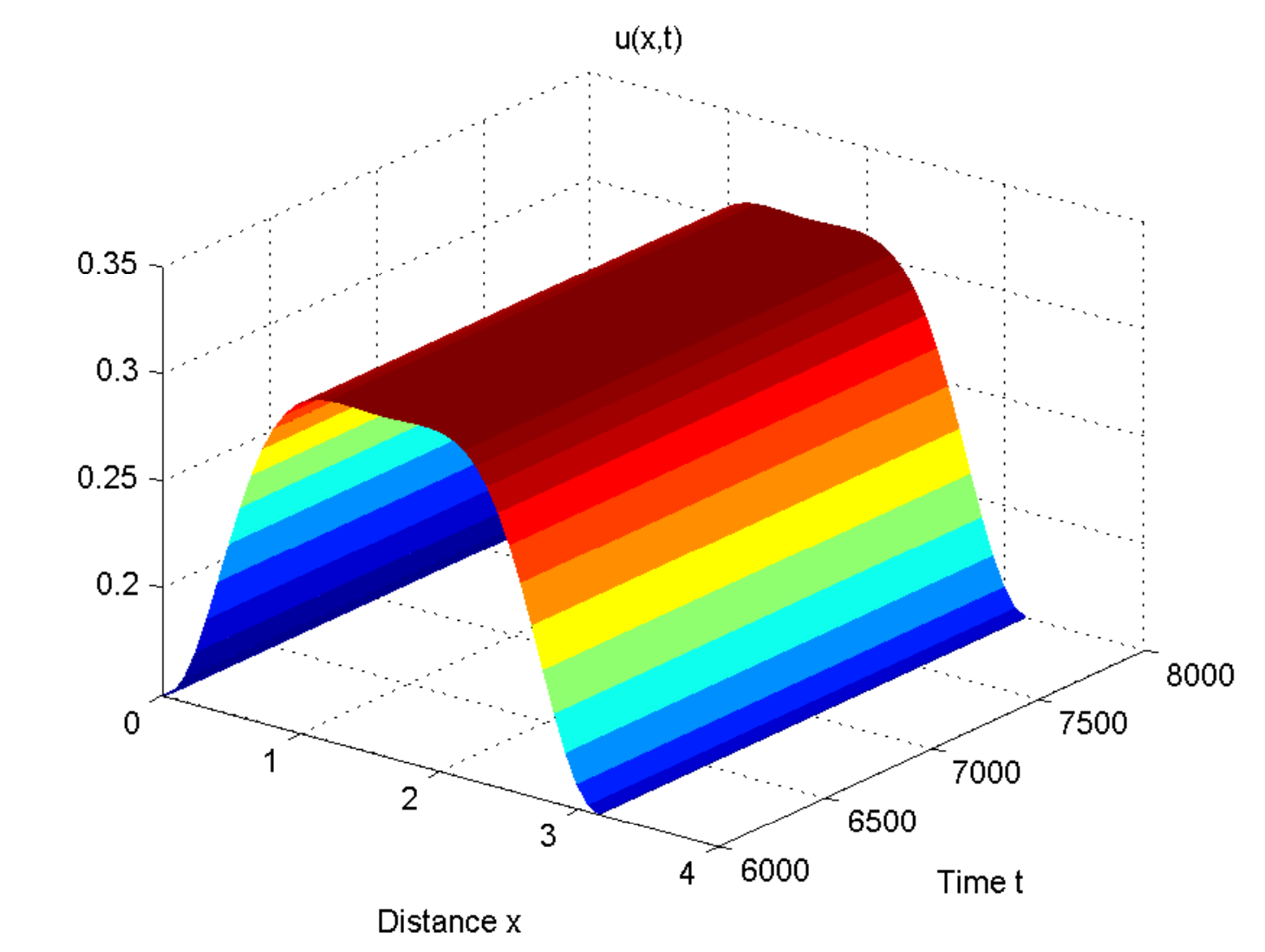}
  \includegraphics[width=0.49\textwidth]{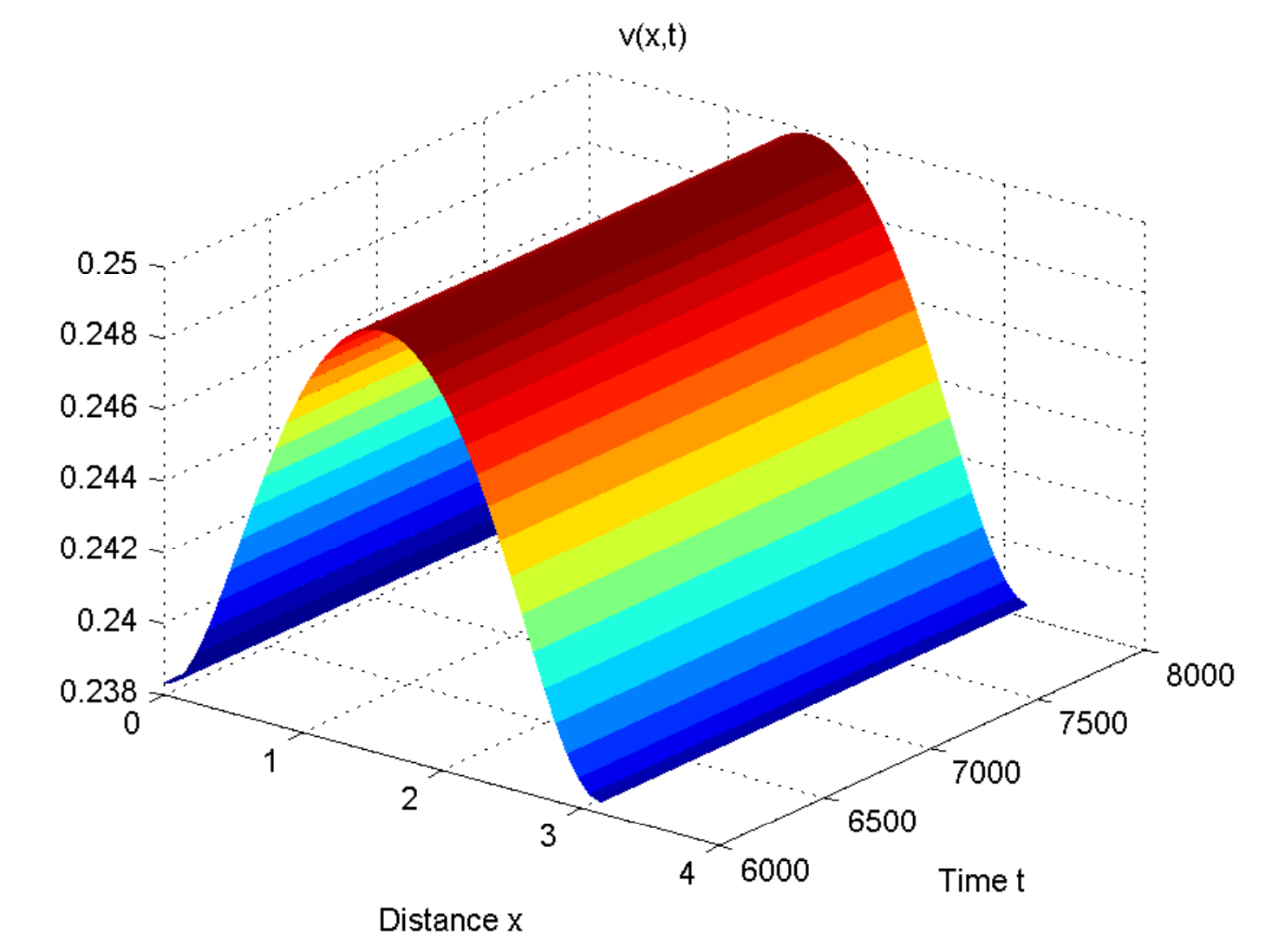}
\end{minipage}
}
\subfigure[The initial values are $u(0,x)=0.245+0.02\cos 2x, v(0,x)=0.245-0.05\cos 2x$.]{
\label{fig:3:b} 
\begin{minipage}[b]{0.48\textwidth}
  \centering\includegraphics[width=0.49\textwidth]{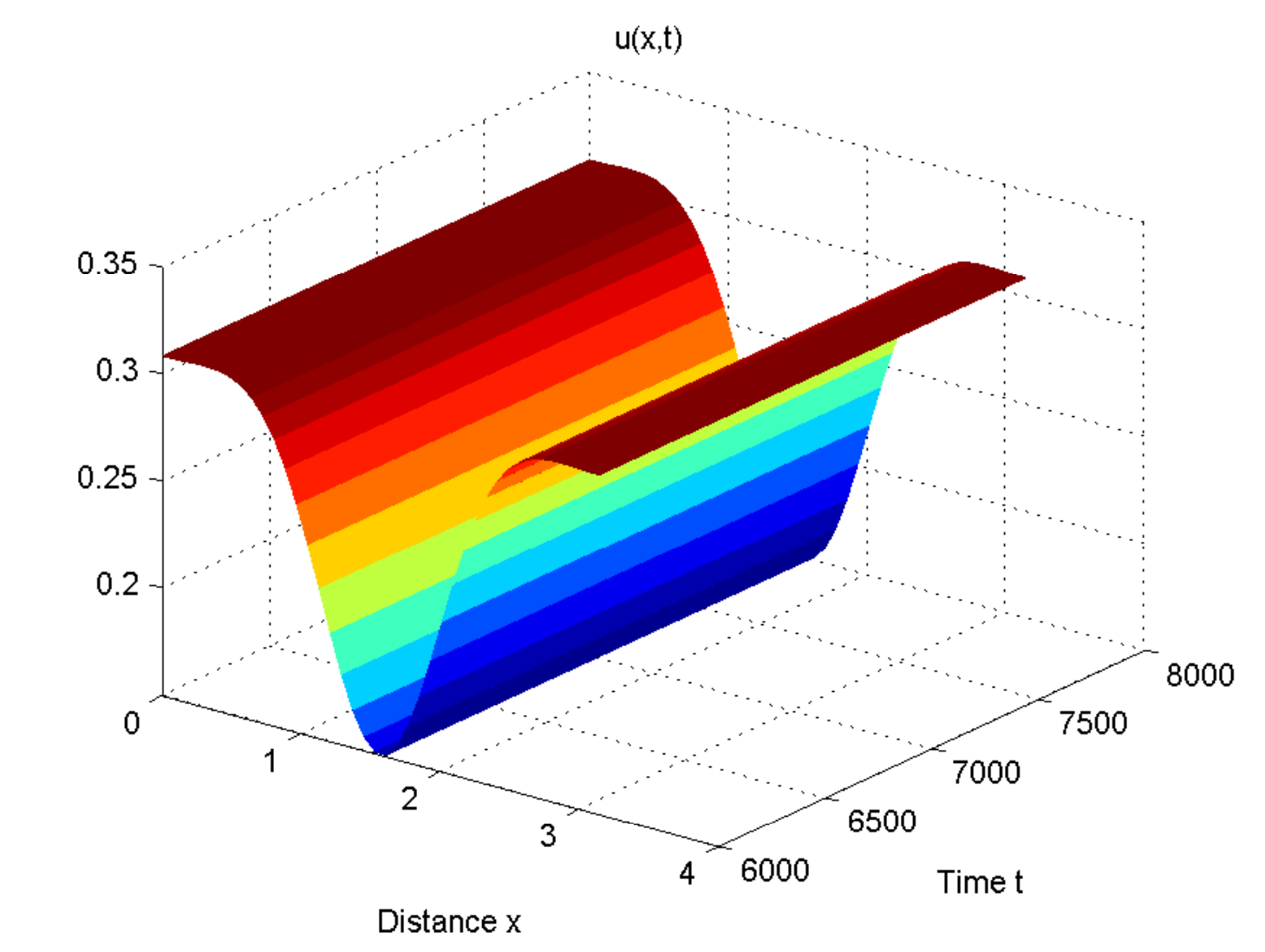}
  \includegraphics[width=0.49\textwidth]{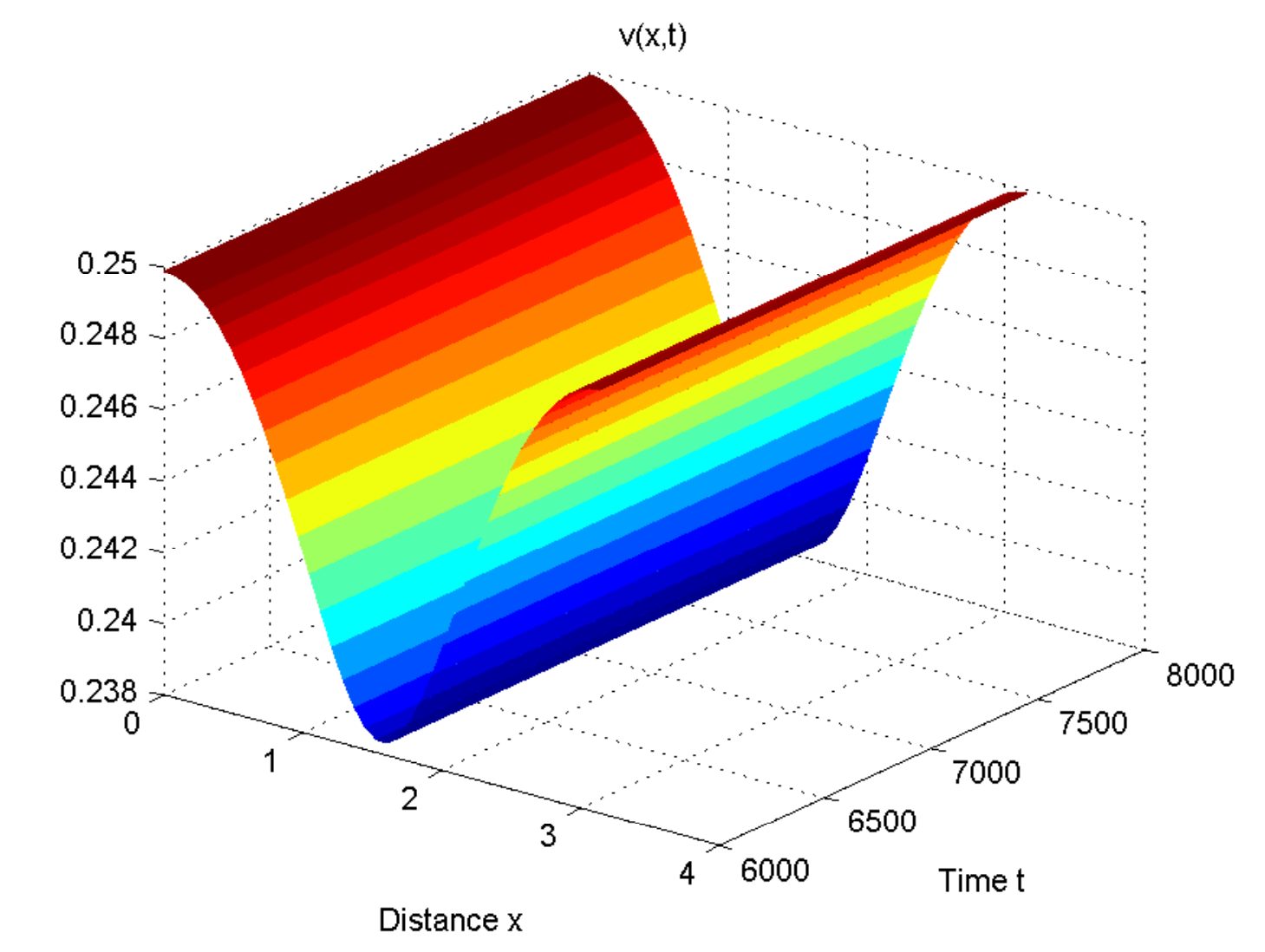}
\end{minipage}
}
\subfigure[The initial values are $u(0,x)=0.245-0.02\cos 3x, v(0,x)=0.245+0.05\cos 3x$.]{
\label{fig:3:c} 
\begin{minipage}[b]{0.48\textwidth}
  \centering\includegraphics[width=0.49\textwidth]{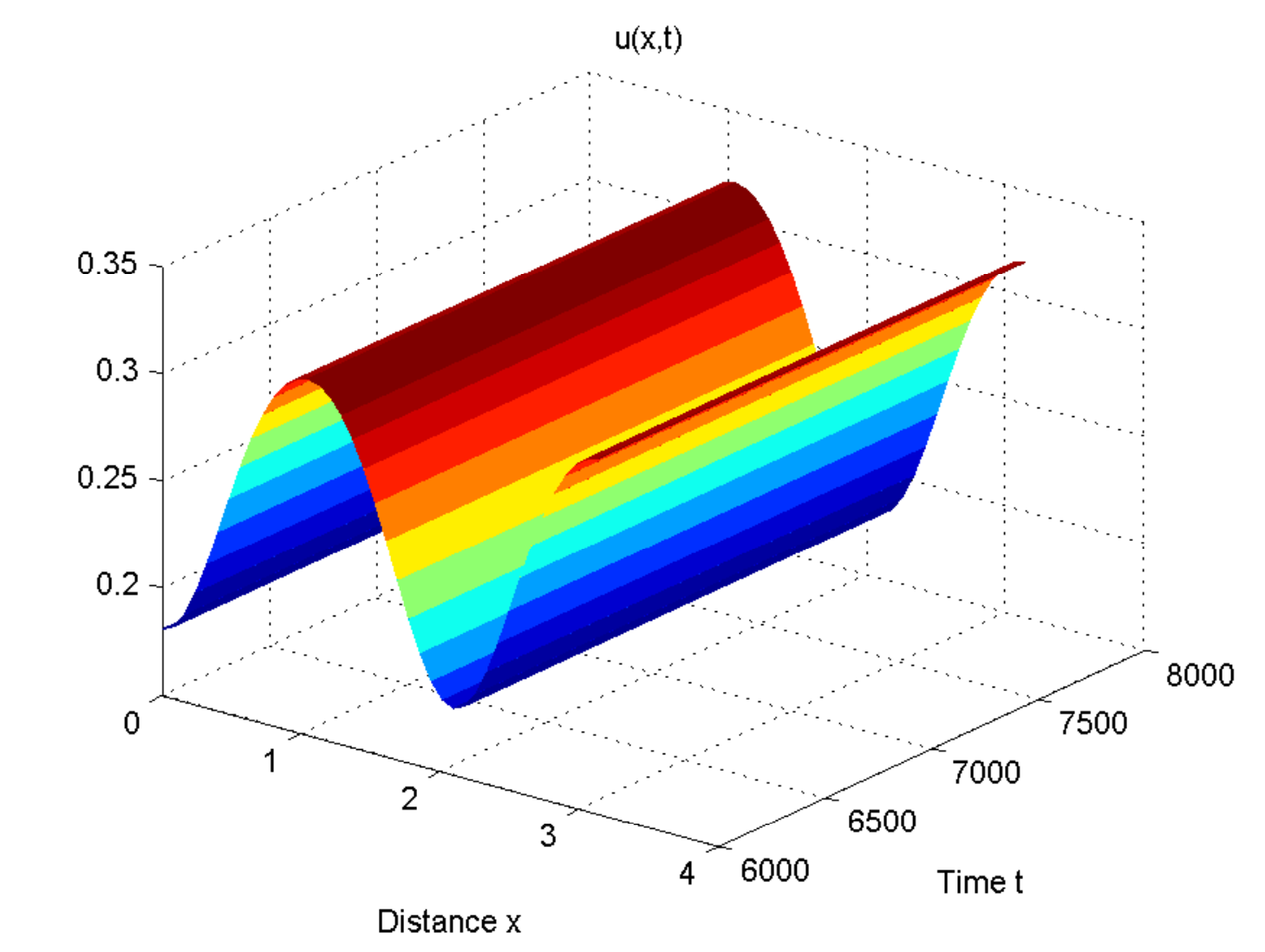}
  \includegraphics[width=0.49\textwidth]{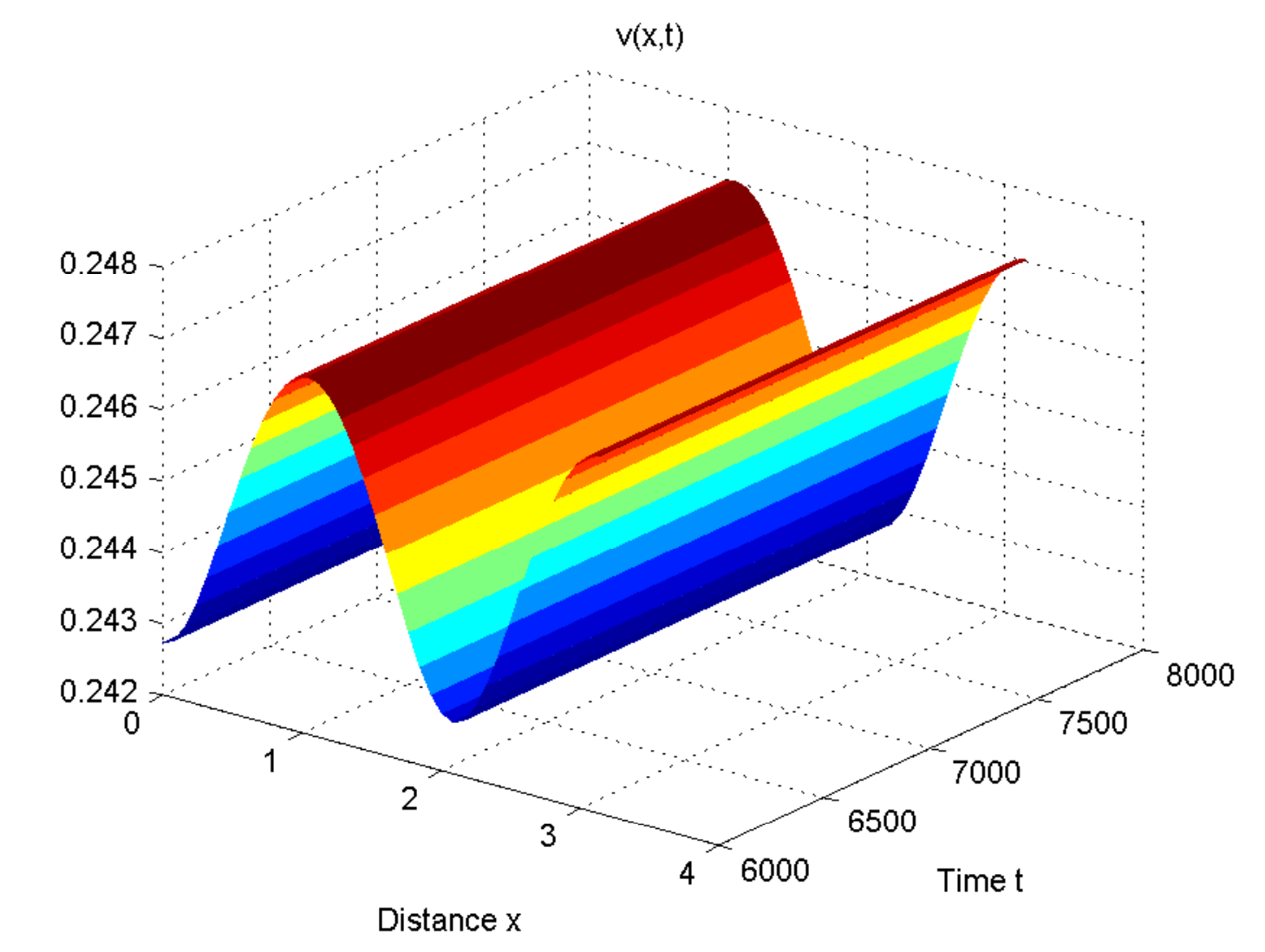}
\end{minipage}
}
\subfigure[The initial values are $u(0,x)=0.245+0.02\cos 3x, v(0,x)=0.245-0.05\cos 3x$.]{
\label{fig:3:d} 
\begin{minipage}[b]{0.48\textwidth}
  \centering\includegraphics[width=0.49\textwidth]{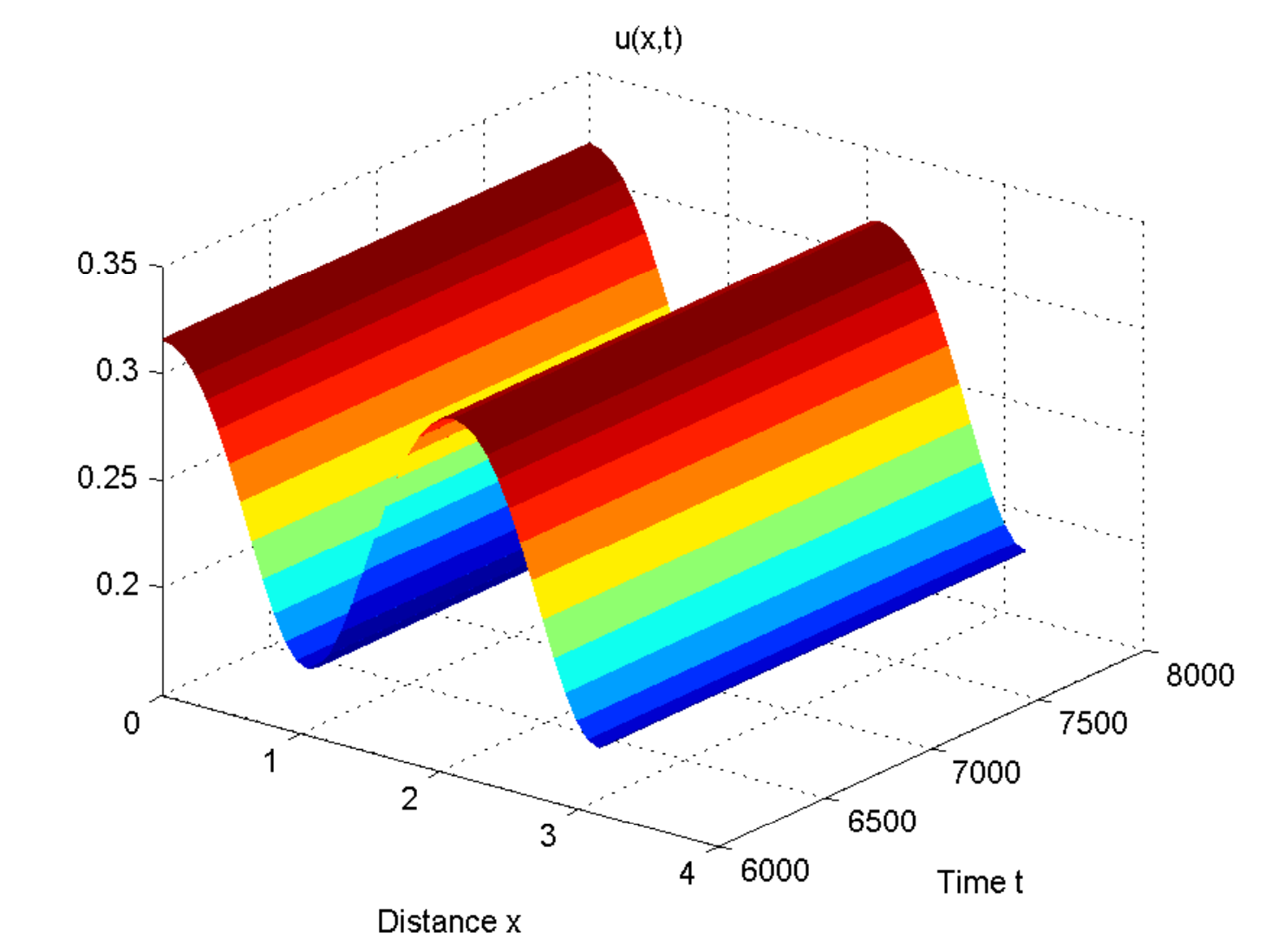}
  \includegraphics[width=0.49\textwidth]{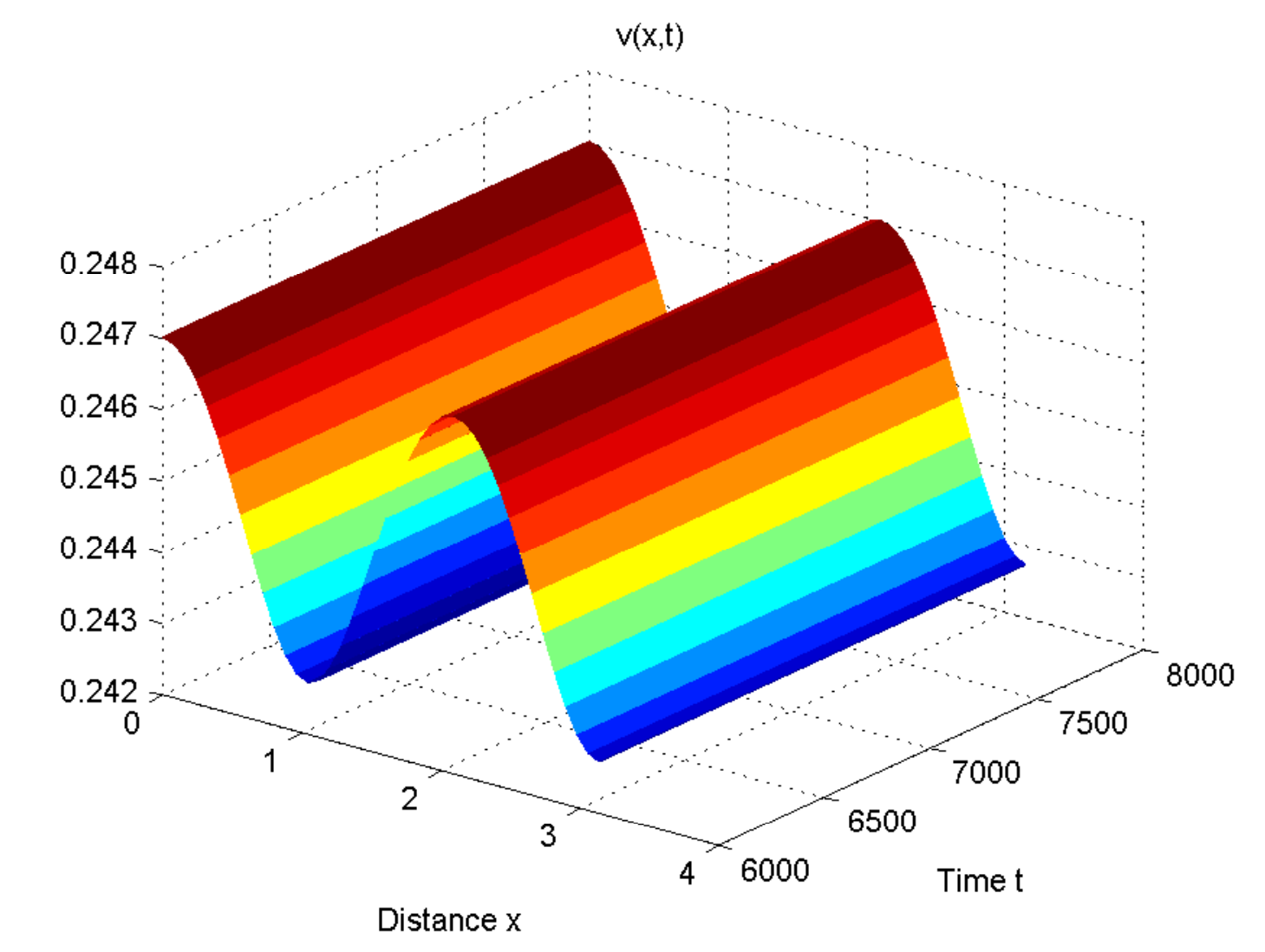}
\end{minipage}
}
\caption{For $\left(d_1,s\right)=(0.0051,0.2064)\in \mathcal{D}_4$, a pair of stable spatially inhomogeneous steady states with the shape of $\phi_1\cos2x-$like and a pair of stable spatially inhomogeneous steady states with the shape of $\phi_2\cos3x-$like, coexist.}
\label{fig:3} 
\end{figure}

Near $(2,3)-$mode Turing-Turing bifurcation point $\left(d_{2,3}^*,s_{2,3}^*\right)$, the $d_1\text{-}s$ parameter plane is divided into six regions. And long-time dynamics of normal form \eqref{e:3.7} can be described by the corresponding phase portraits respectively, when $(d_1,s)$ are chosen in these regions.
Then based on \emph{center manifold theory} \cite{CJ,WJH} and Table 1, we could reveal long-time dynamics of system \eqref{e:3.1}, when parameters are chosen in these six regions, respectively. And, dynamics of diffusive system \eqref{e:3.1} are concluded as follows,
\begin{proposition}\label{pro:4.11}
For fixed $m=6,a=3,b=0.5,d_2=0.7$, diffusive predator-prey system \eqref{e:3.1} exhibits complex spatial patterns, when parameters $(d_1,s)$ are chosen near $(2,3)-$mode Turing-Turing bifurcation point $\left(d_{2,3}^*,s_{2,3}^*\right)=(0.0056,0.2364)$. Here are the results:
\begin{enumerate}
\item For $\left(d_1,s\right)\in \mathcal{D}_1$, the coexistence equilibrium $E_*$ of system \eqref{e:3.1} is asymptotically stable. Otherwise, the coexistence equilibrium $E_*$ is unstable while $(d_1,s)\notin{\mathcal{D}}_1$.
\item When $\left(d_1,s\right)$ crosses from $\mathcal{D}_1$ into $\mathcal{D}_2$, a pair of stable spatially inhomogeneous steady states with the shape of $\phi_2\cos3x-$like bifurcates from the coexistence equilibrium $E_*$ through Turing bifurcation. Therefore, system \eqref{e:3.1} supports \emph{bi-stability}. And, $\phi_2$ is given in \eqref{e:3.6-1}.
\item When $\left(d_1,s\right)$ crosses from $\mathcal{D}_2$ into $\mathcal{D}_3$, the pair of stable spatially inhomogeneous steady states with the shape of $\phi_1\cos3x-$like persists, and a pair of unstable spatially inhomogeneous steady states with the shape of $\phi_2\cos2x-$like bifurcates from the coexistence equilibrium $E_*$ through Turing bifurcation, where $\phi_1$ is given in \eqref{e:3.6-1}. Thus, system \eqref{e:3.1} admits transient spatial patterns and \emph{bi-stability}.
\item When $\left(d_1,s\right)$ crosses from $\mathcal{D}_3$ into $\mathcal{D}_4$, the pair of spatially inhomogeneous steady states with the shape of $\phi_2\cos2x-$like becomes stable and four \emph{saddle-type} superposition steady states with the shape of $\phi_1h_1\cos 2x+\phi_2h_2\cos3x-$like emerge, through Turing bifurcation, where $(h_1,h_2)$ is the coexistence equilibrium of normal form \eqref{e:3.7}.
    Hence, system \eqref{e:3.1} supports the \emph{coexistence of four stable spatially inhomogeneous steady states} (see Fig. \ref{fig:3}), considering that the pair of stable spatially inhomogeneous steady states with the shape of $\phi_1\cos3x-$like still remains.
\item When $\left(d_1,s\right)$ crosses from $\mathcal{D}_4$ into $\mathcal{D}_5$, the pair of stable spatially inhomogeneous steady states with the shape of $\phi_1\cos3x-$like loses its stability by colliding with the four \emph{saddle-type} superposition steady states with the shape of $\phi_1h_1\cos 2x+\phi_2h_2\cos3x-$like through Turing bifurcation. Moreover, the pair of stable spatially inhomogeneous steady states with the shape of $\phi_1\cos2x-$like persists. Then, system \eqref{e:3.1} exhibits transient patterns and \emph{bi-stability}.
\item When $\left(d_1,s\right)$ crosses from $\mathcal{D}_5$ into $\mathcal{D}_6$, the pair of unstable spatially inhomogeneous steady states with the shape of $\phi_1\cos3x-$like collides with the coexistence equilibrium $E_*$ and disappears, through Turing bifurcation. Thus, system \eqref{e:3.1} supports \emph{bistable spatial patterns}, given that the pair of stable spatially inhomogeneous steady states with the shape of $\phi_1\cos2x-$like persists.
\item When $\left(d_1,s\right)$ crosses from $\mathcal{D}_6$ back into $\mathcal{D}_1$, the coexistence equilibrium $E_*$ regains its stability by colliding with the pair of stable spatially inhomogeneous steady states with the shape of $\phi_1\cos2x-$like through Turing bifurcation.
\end{enumerate}
\end{proposition}

As we anticipated in \cite{CX&JWH} that under proper conditions, diffusive predator-prey system \eqref{e:3.1} with Crowley-Martin functional response will exhibit complex spatial patterns that four stable spatially inhomogeneous steady states coexist, we theoretically demonstrate this conjecture.

\begin{figure}[!hbt]
\centering
\subfigure[Bifurcation set in $d_1\text{-}s$ plane. $TT_{1,2}$ is $(1,2)-$mode Turing-Turing bifurcation point.]{
\label{fig:4:a}
\begin{minipage}[b]{0.483\textwidth}
  \centering\includegraphics[width=\textwidth]{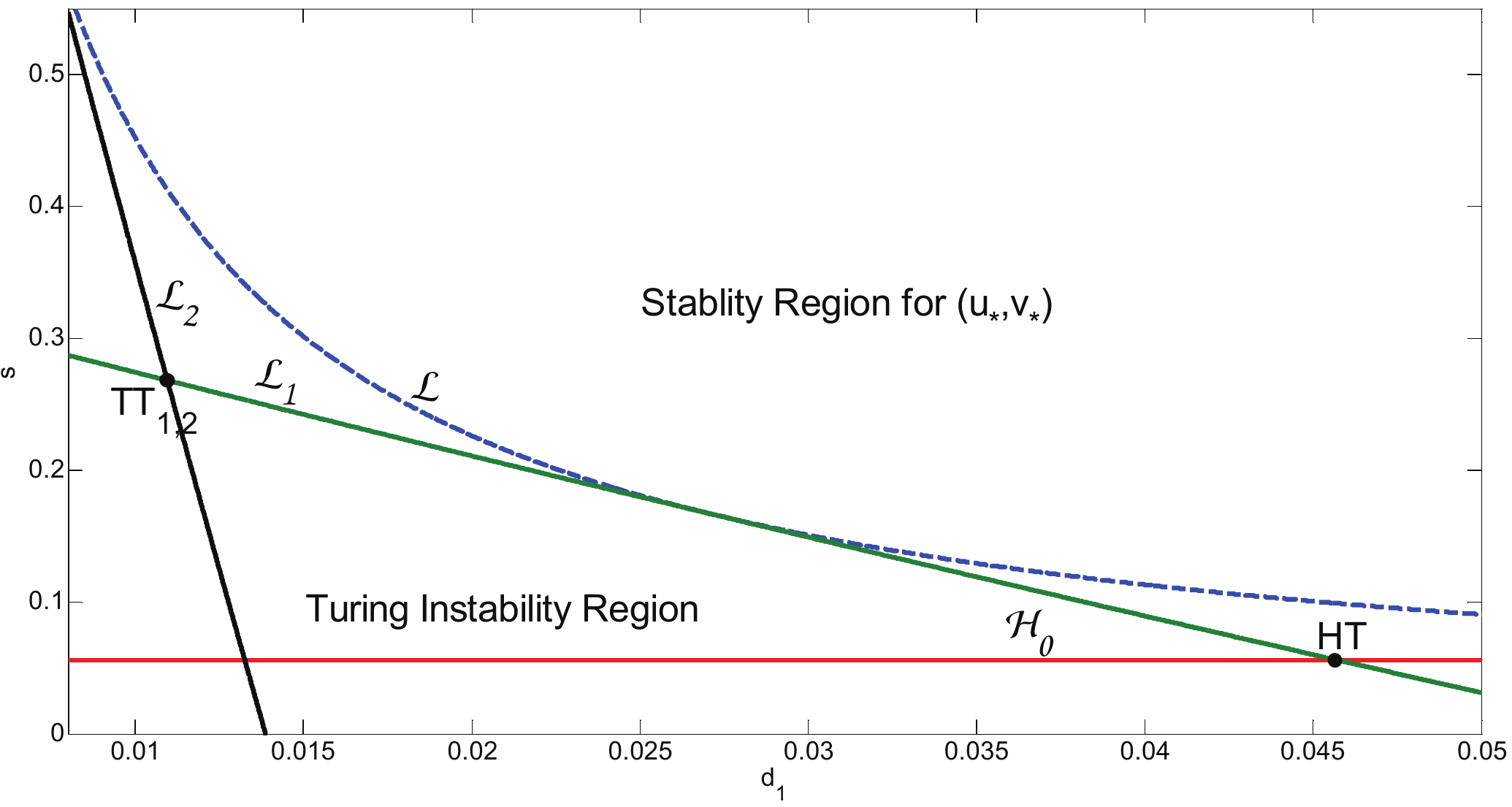}
\end{minipage}
}
\subfigure[Bifurcation set with local bifurcation curves $\mathcal{T}_1$ at Turing-Turing point $\left(d_{1,2}^*,s_{1,2}^*\right)$.]{
\label{fig:4:b}
\begin{minipage}[b]{0.48\textwidth}
\centering \includegraphics[width=\textwidth]{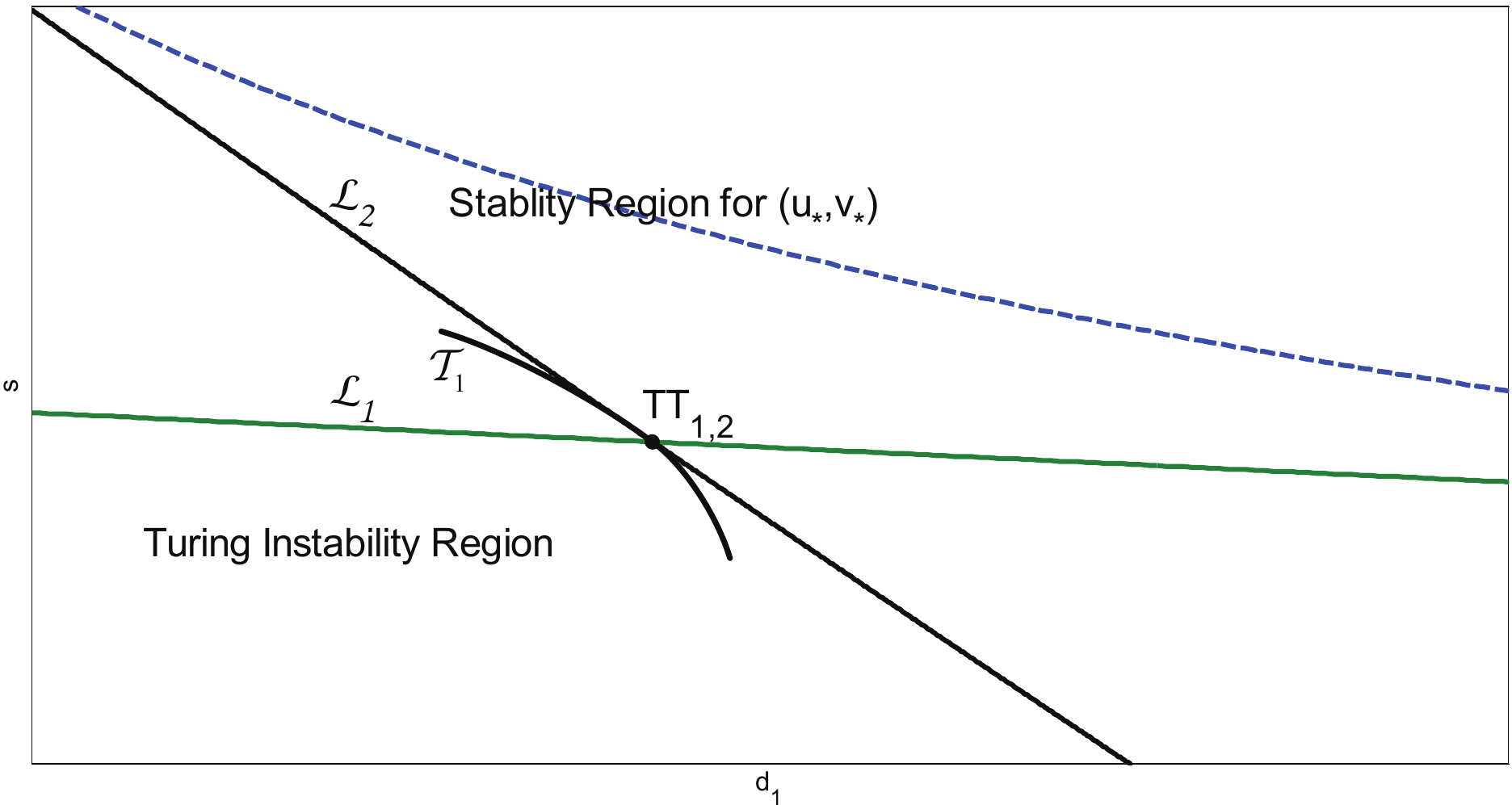}
\end{minipage}
}
\caption{\rm Bifurcation set and Turing bifurcation curves (a), and bifurcation set with local bifurcation curves (b) in $d_1\text{-}s$ plane.}
  \label{fig:4}
\end{figure}

Next, we show that system \eqref{e:3.1} could generate other more complex spatial patterns for another set of parameters. Analogously, let $m=5,a=3,b=0.1,d_2=4$, then $(u_*,v_*)=(0.2716,0.2716)$. Also, $s_0=0.0555, \sigma=-0.7092$ such that $s_0+\sigma<0$. Furthermore, we derive $k_0^*=1$, and
\begin{equation*}
\begin{split}
\mathcal{L}_k:s=s_k(d_1)=\frac{s_0d_2k^2-d_1d_2k^4}{d_1k^2-(s_0+\sigma)}=\frac{-4d_1 k^4+0.2219k^2}{d_1 k^2+0.6537},\quad 0< d_1<\frac{s_0}{k^2},k\in\mathbb{N}.
\end{split}
\end{equation*}
Therefore, the corresponding Turing bifurcation curves in $d_1\text{-}s$ plane, are shown in Fig. \ref{fig:4:a}.

Then, $\left(d_{1,2}^*,s_{1,2}^*\right)=(0.01095,0.2679)$ is the \emph{$(1,2)-$mode Turing-Turing bifurcation point}, when taking ${k_1}=1,{k_2}=2$ in Theorem \ref{thm:2.1}.

\begin{figure}[!hbt]
\centering
\subfigure[Local bifurcation set in $d_1\text{-}s$ plane.]{
\label{fig:5:a}
\begin{minipage}[b]{0.42\textwidth}
  \centering \includegraphics[width=\textwidth]{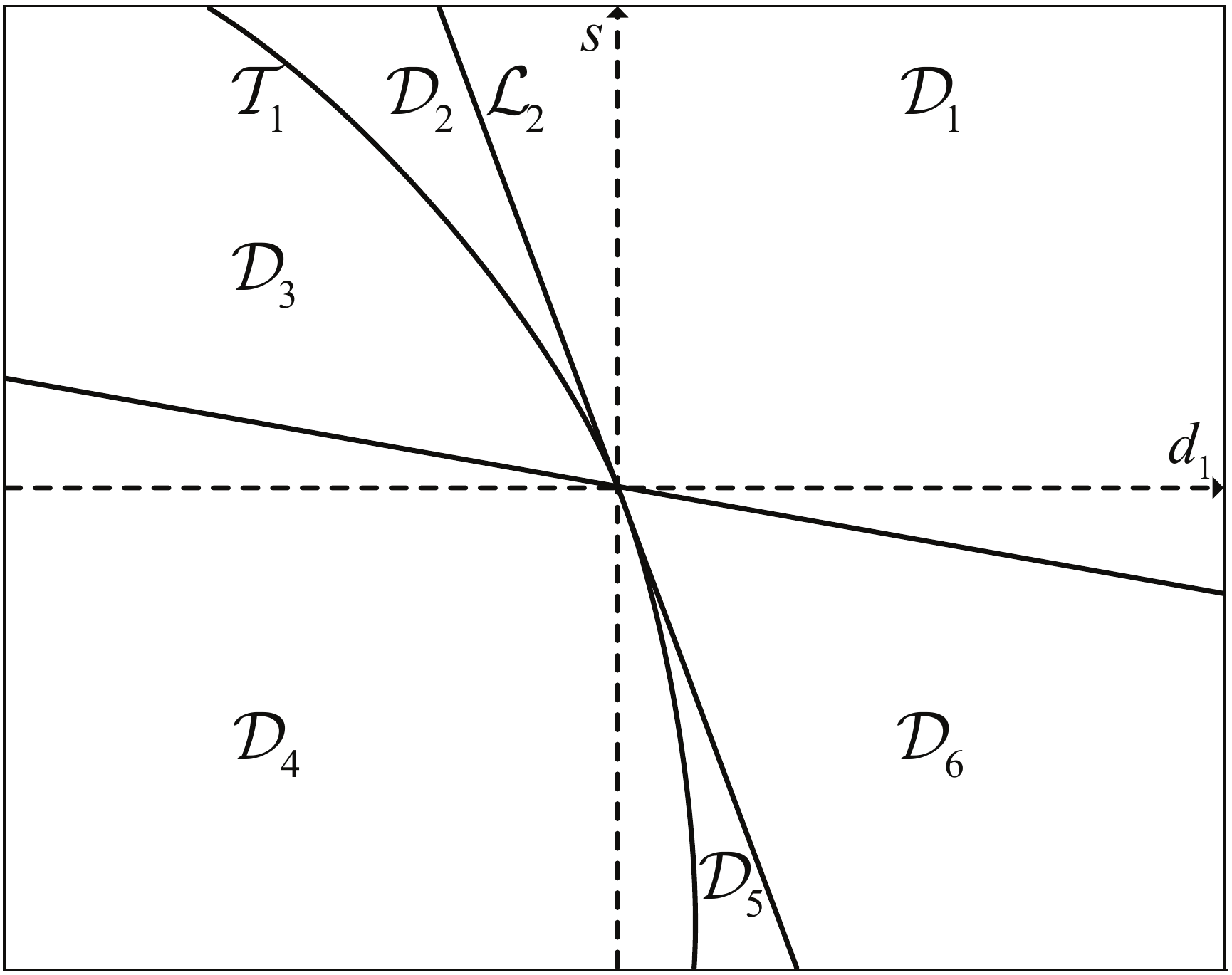}
\end{minipage}
}
\subfigure[Phase portraits.]{
\label{fig:5:b}
\begin{minipage}[b]{0.53\textwidth}
\centering \includegraphics[width=\textwidth]{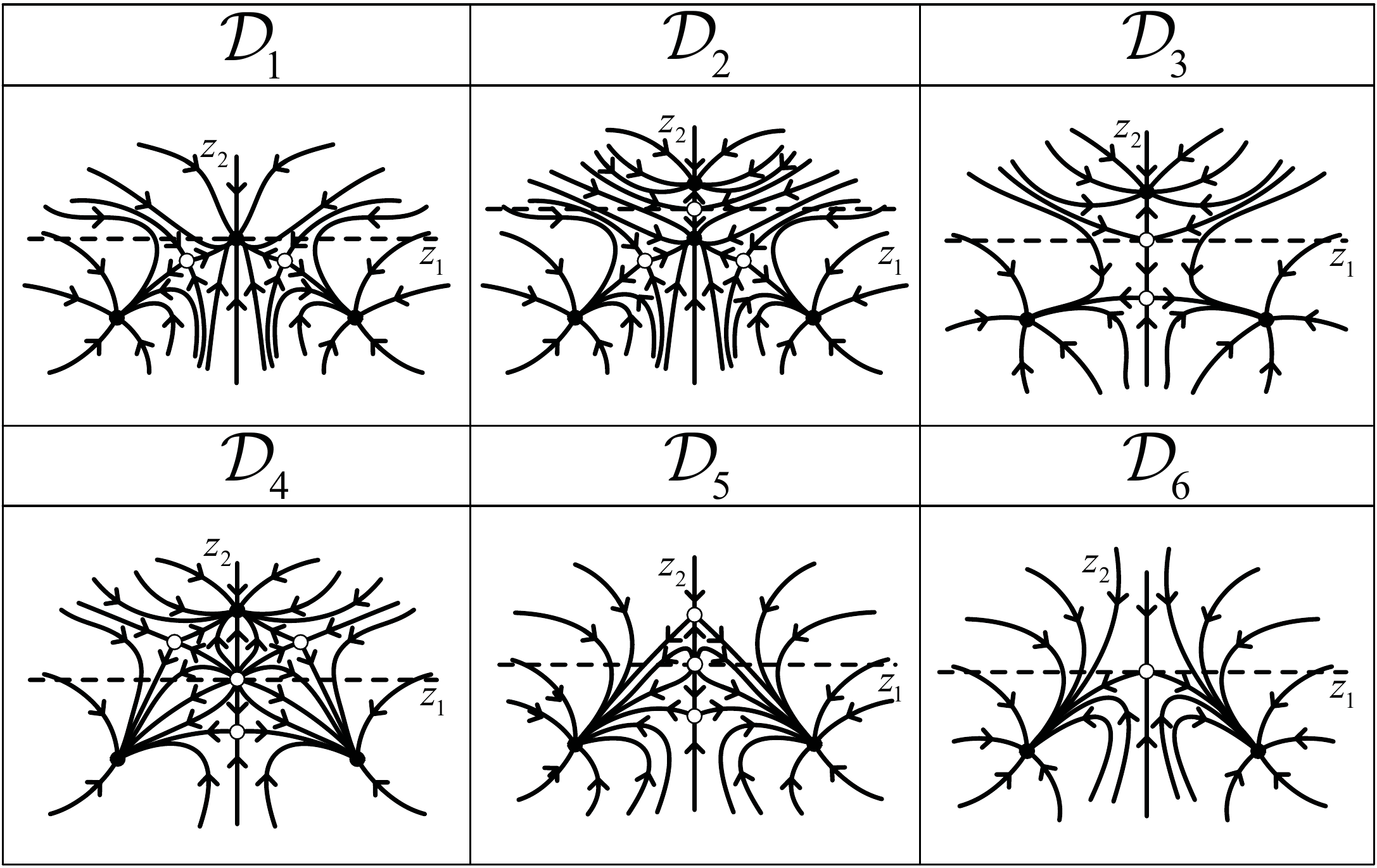}
\end{minipage}
}
\caption{\rm Local bifurcation set (a) at $(1,2)-$mode Turing-Turing point $\left(d_{1,2}^*,s_{1,2}^*\right)$, and the corresponding phase portraits (b) for normal form \eqref{e:3.8}.}
\label{fig:5}
\end{figure}

Applying \emph{Proposition \ref{pro:2.4}}, for this set of given parameters, normal form truncated to order 3 reads
\begin{equation}\label{e:3.8}
\left\{
\begin{aligned}
\dot{z}_1&=-(1.0105\epsilon_1+0.1574\epsilon_2)z_1-0.3461z_1z_2-2.8448z_1^3-2.4251z_1z_2^2,\\
\dot{z}_2&=-(4.0028\epsilon_1+0.04291\epsilon_2)z_2-0.2750z_1^2+1.2199z_1^2z_2-2.5756z_2^3.
\end{aligned}\right.
\end{equation}
By analyzing normal form \eqref{e:3.8}, we have the following critical bifurcation curves in $d_1\text{-}s$ plane,
\begin{equation*}
\begin{aligned}
\mathcal{L}_2:&s=s_{1,2}^*-93.2920\left(d_1-d_{1,2}^*\right),\\
\mathcal{T}_1: &\left(-0.07136+0.07136\sqrt{1+80.9673\vartheta}\right)^3-0.02379\left(-0.07136+0.07136\sqrt{1+80.9673\vartheta}\right)^2\\
&+\left(1.2270\left(d_1-d_{1,2}^*\right)+0.03053\left(s-s_{1,2}^*\right)-0.009255\right)\left(-0.07136+0.07136\sqrt{1+80.9673\vartheta}\right)\\
&+0.02674\vartheta=0,\quad\textrm{with}\; \vartheta\triangleq-\left(1.0105\left(d_1-d_{1,2}^*\right)+0.1574\left(s-s_{1,2}^*\right)\right),
\end{aligned}
\end{equation*}
as well as local bifurcation set and the corresponding phase portraits shown in Fig. \ref{fig:5}.

Analogously, local bifurcation set could also be embedded into global bifurcation set, see Fig. \ref{fig:4:b}.

\begin{figure}[!hbt]
\subfigure[The initial values are $u(0,x)=0.2716-0.1\cos x, v(0,x)=0.2716-0.1\cos(x)$.]{
\label{fig:6:a} 
\begin{minipage}[b]{0.48\textwidth}
  \centering\includegraphics[width=0.5\textwidth]{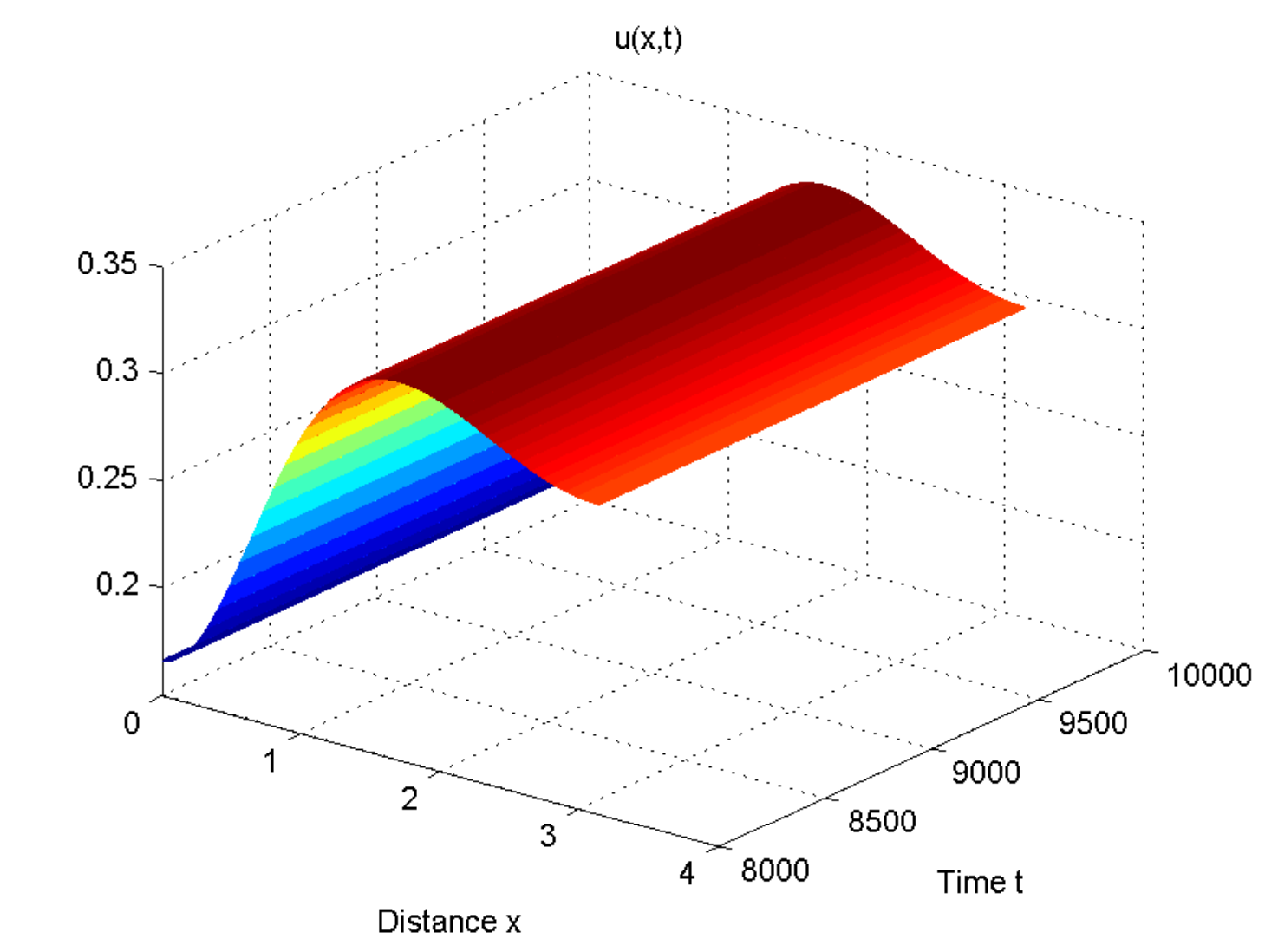}\includegraphics[width=0.5\textwidth]{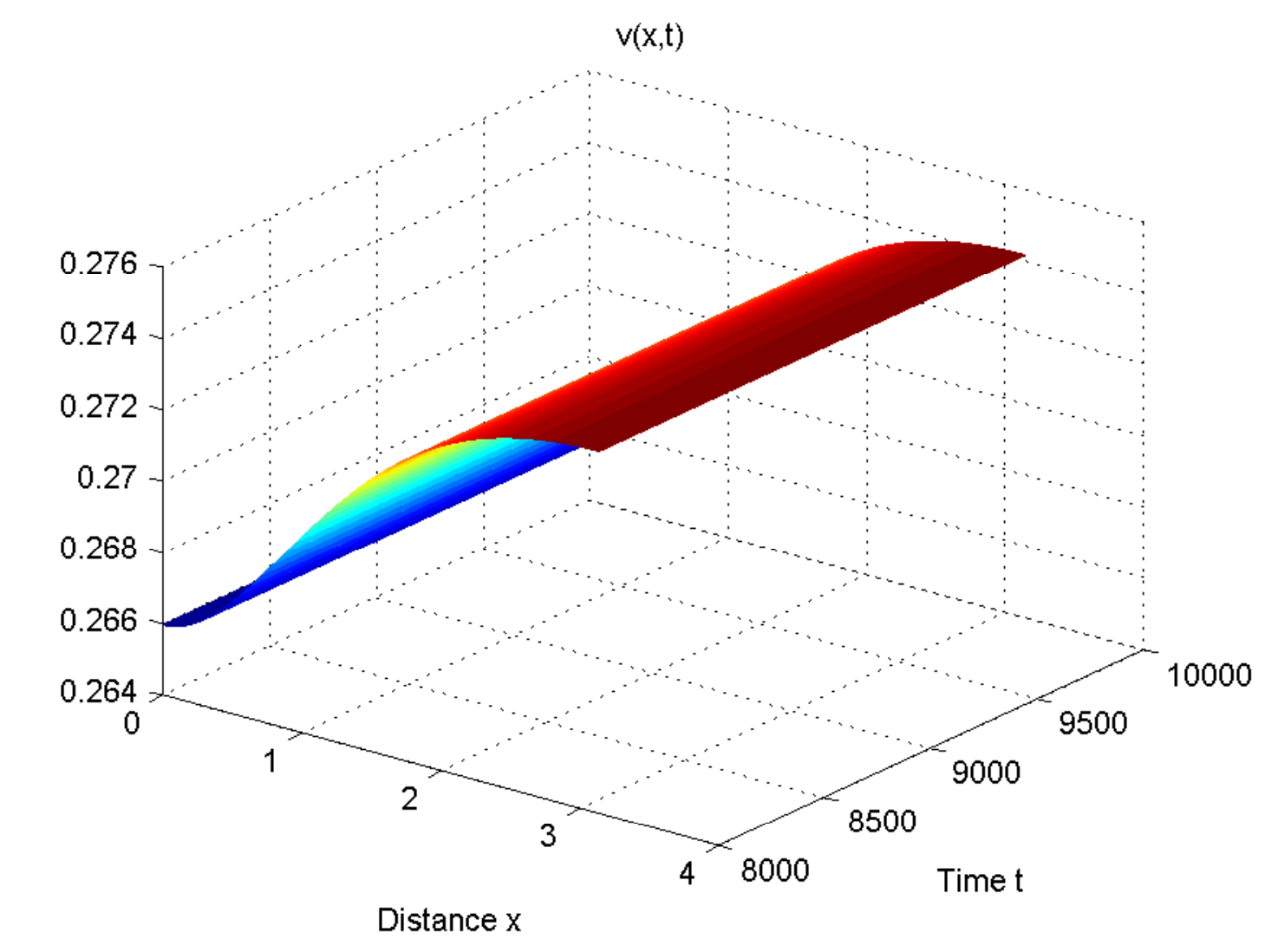}
\end{minipage}
}
\subfigure[The initial values are $u(0,x)=0.2716+0.1\cos x, v(0,x)=0.2716+0.1\cos x$.]{
\label{fig:6:b} 
\begin{minipage}[b]{0.48\textwidth}
  \centering\includegraphics[width=0.5\textwidth]{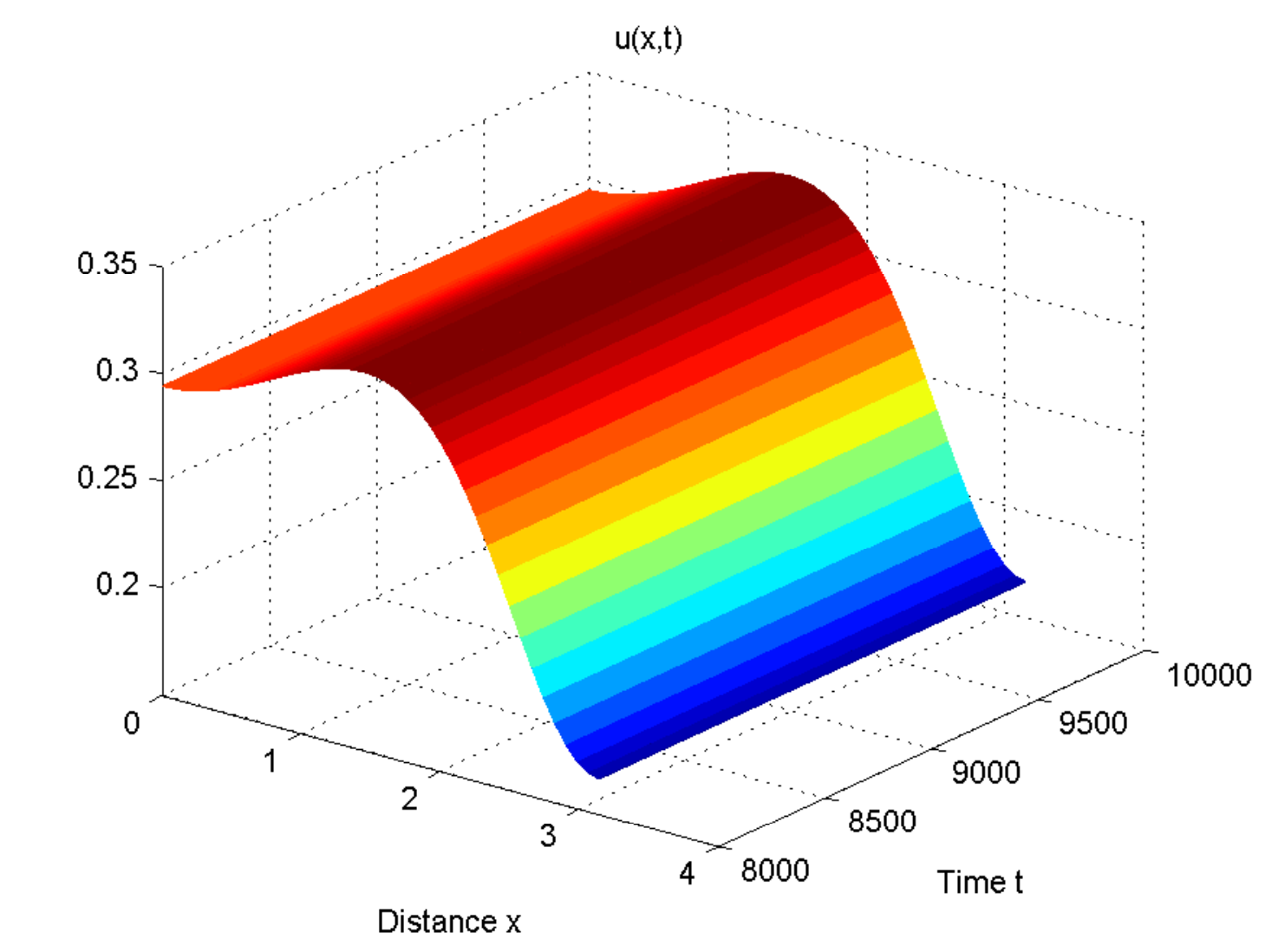}\includegraphics[width=0.5\textwidth]{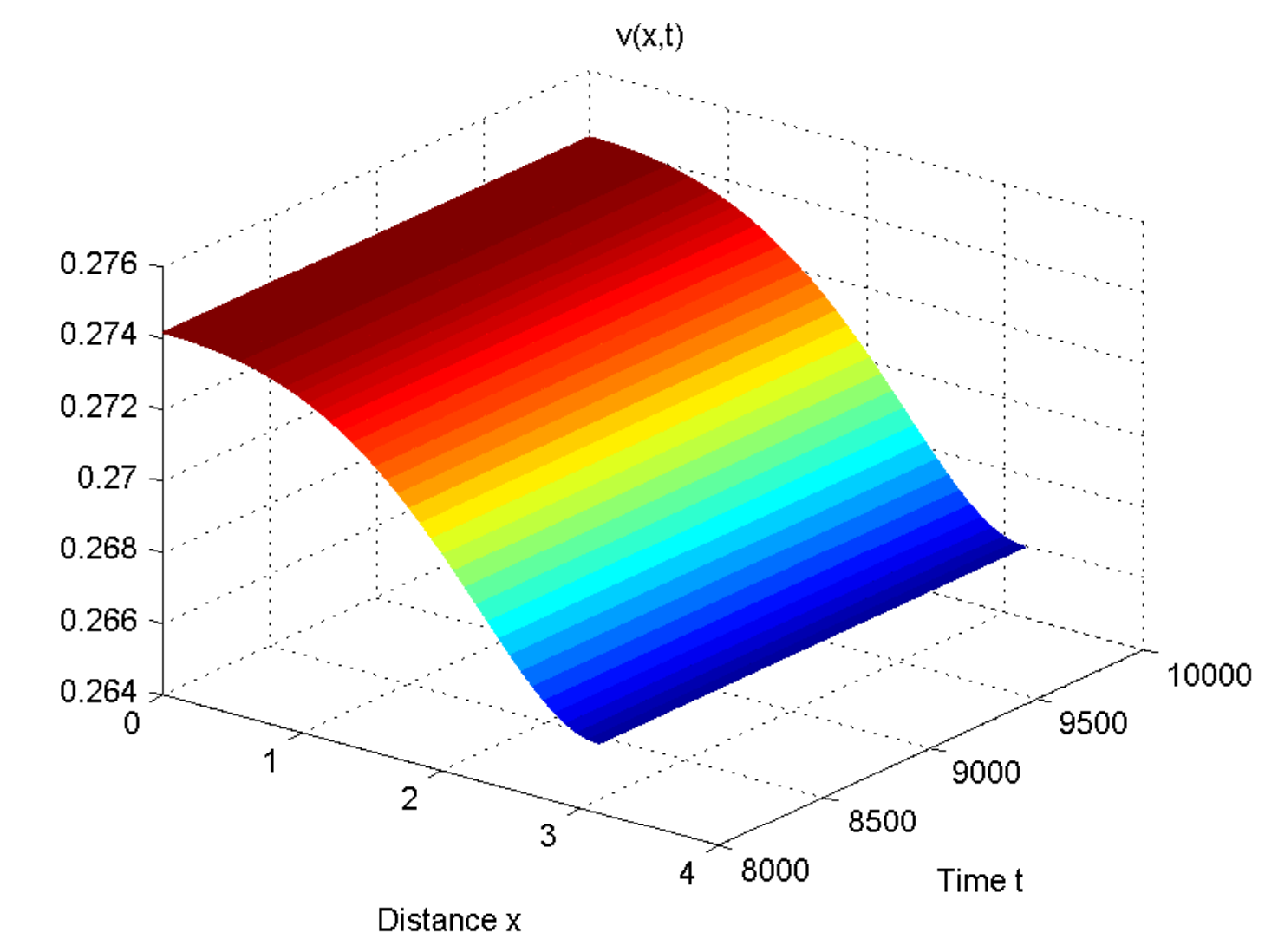}
\end{minipage}
}
\begin{center}
\subfigure[The initial values are $u(0,x)=0.2716-0.02\cos 2x, v(0,x)=0.2716-0.05\cos 2x$.]{
\label{fig:6:c} 
\begin{minipage}[b]{0.48\textwidth}
  \centering\includegraphics[width=0.5\textwidth]{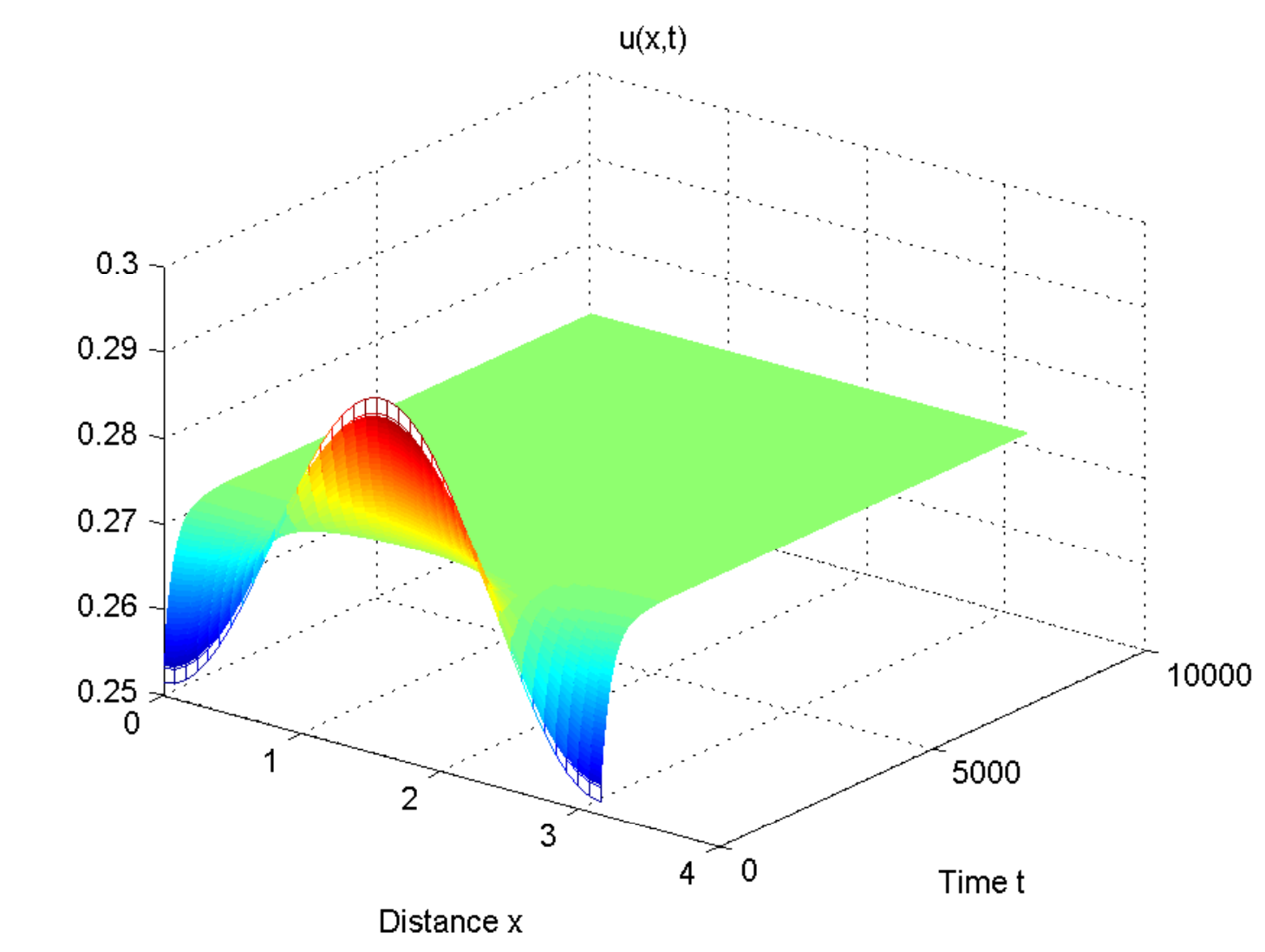}\includegraphics[width=0.5\textwidth]{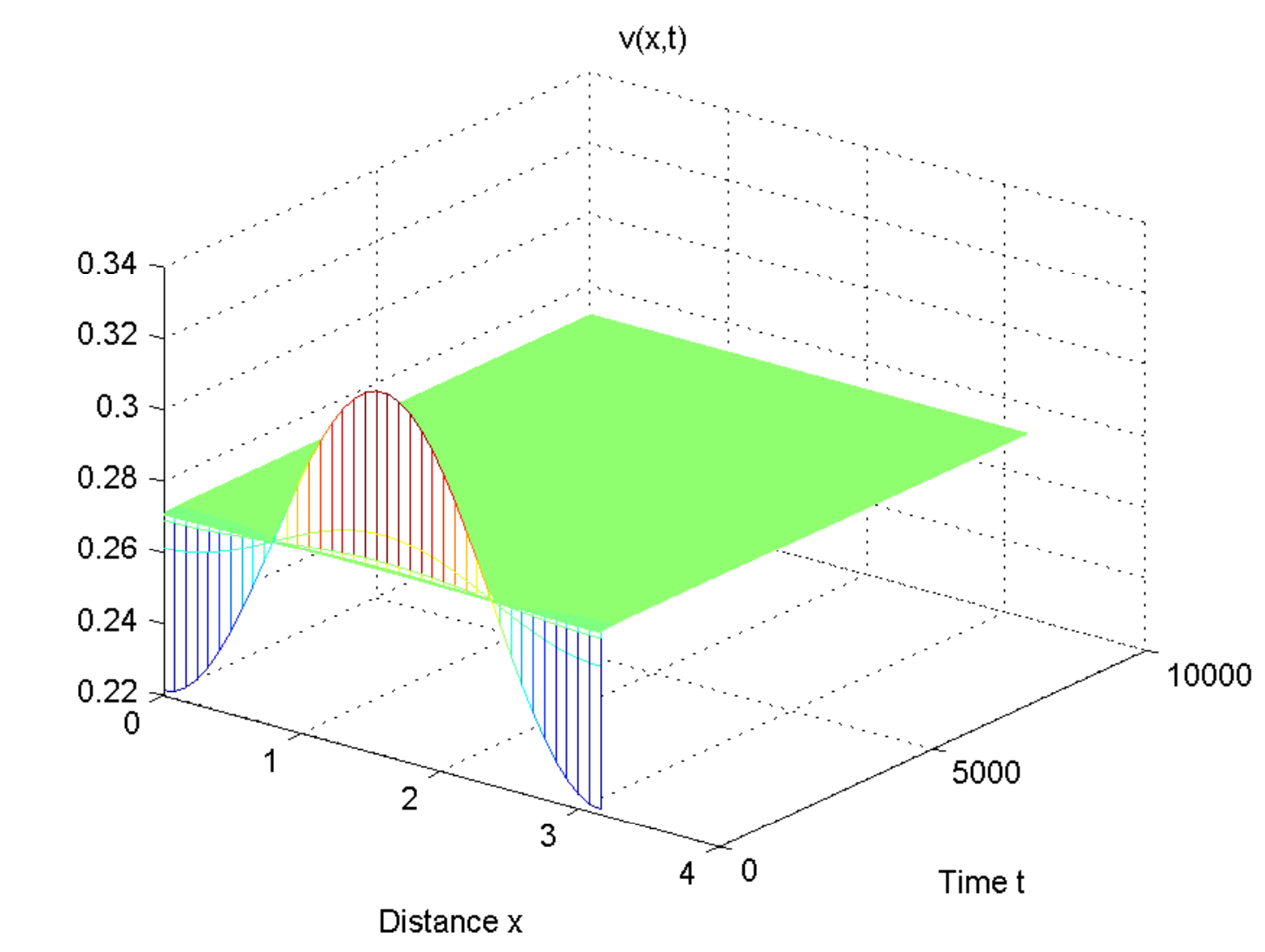}
\end{minipage}
}
\end{center}
\caption{\rm For $\left(d_1,s\right)=(0.01195,0.2679)\in \mathcal{D}_1$, a pair of stable superposition steady states with the shape of $\phi_1h_1\cos x+\phi_2h_2\cos2x-$like and the stable coexistence equilibrium $E_*$, coexist.}
\label{fig:6}
\end{figure}

\begin{figure}[!hbt]
\centering
\subfigure[The initial values are $u(0,x)=0.2716-0.1\cos x, v(0,x)=0.2716-0.1\cos x$.]{
\label{fig:7:a} 
\begin{minipage}[b]{0.48\textwidth}
  \centering\includegraphics[width=0.49\textwidth]{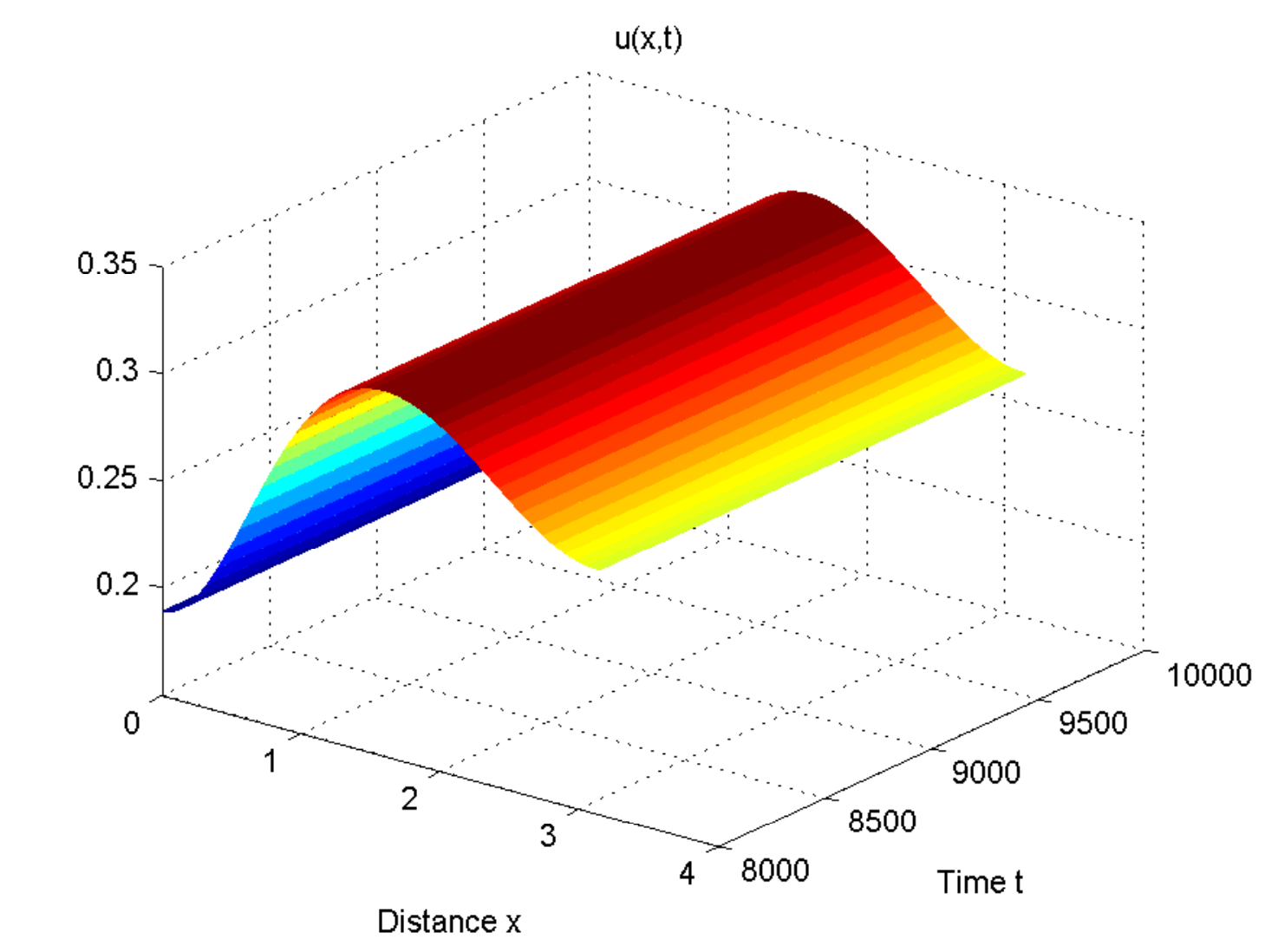}\includegraphics[width=0.49\textwidth]{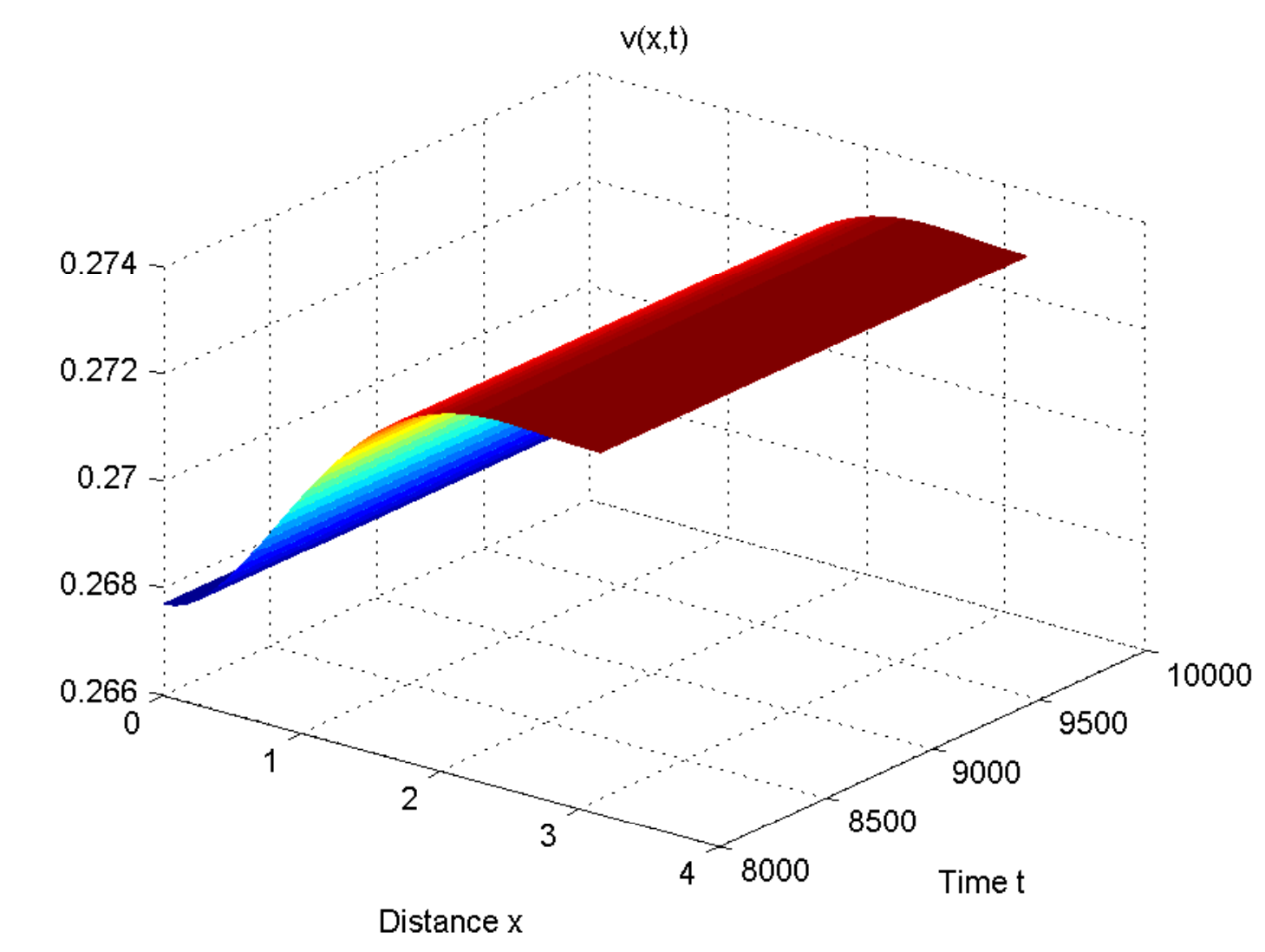}
\end{minipage}
}
\subfigure[The initial values are $u(0,x)=0.2716+0.1\cos x, v(0,x)=0.2716+0.1\cos x$.]{
\label{fig:7:b} 
\begin{minipage}[b]{0.48\textwidth}
  \centering\includegraphics[width=0.49\textwidth]{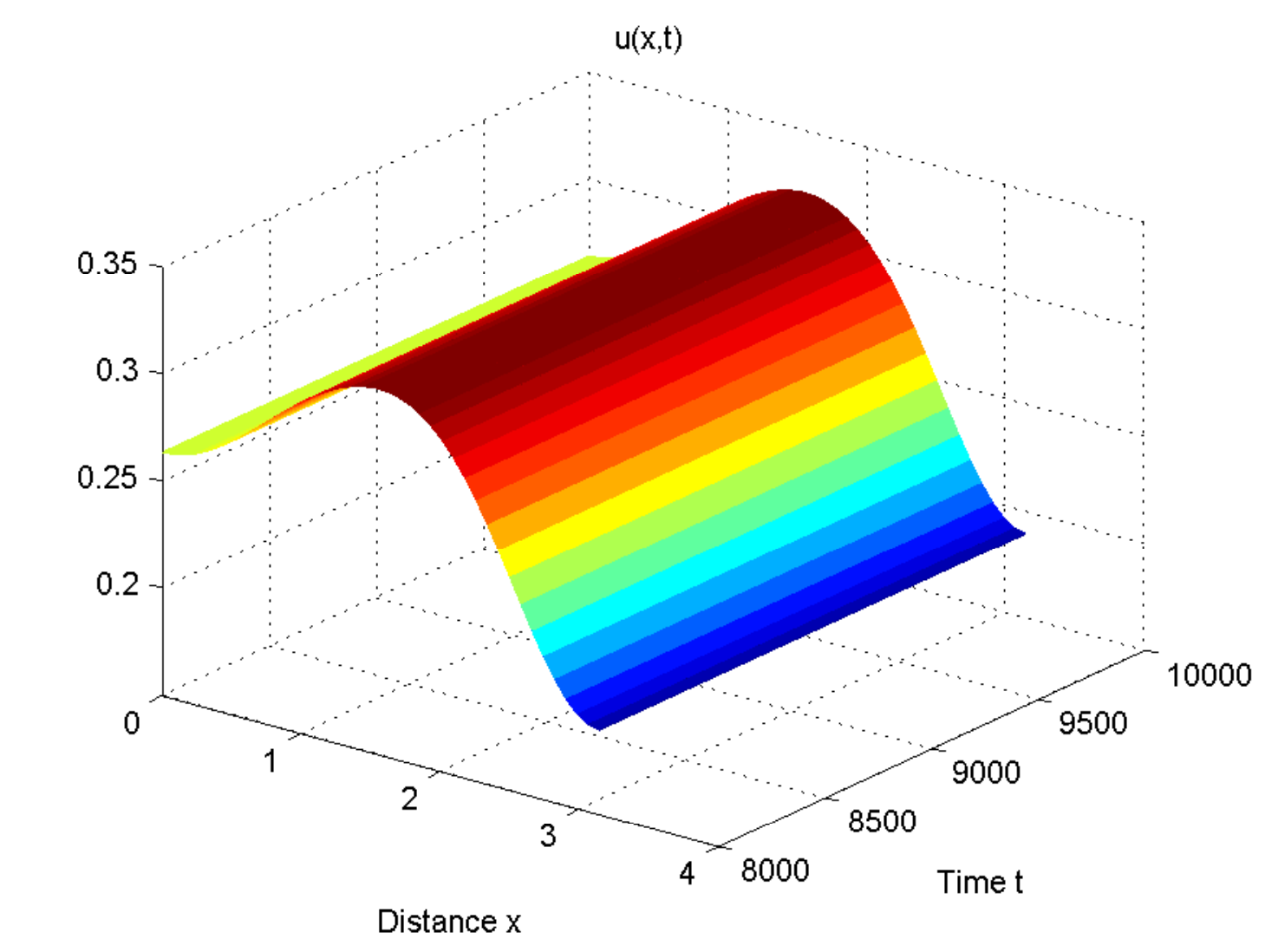}\includegraphics[width=0.49\textwidth]{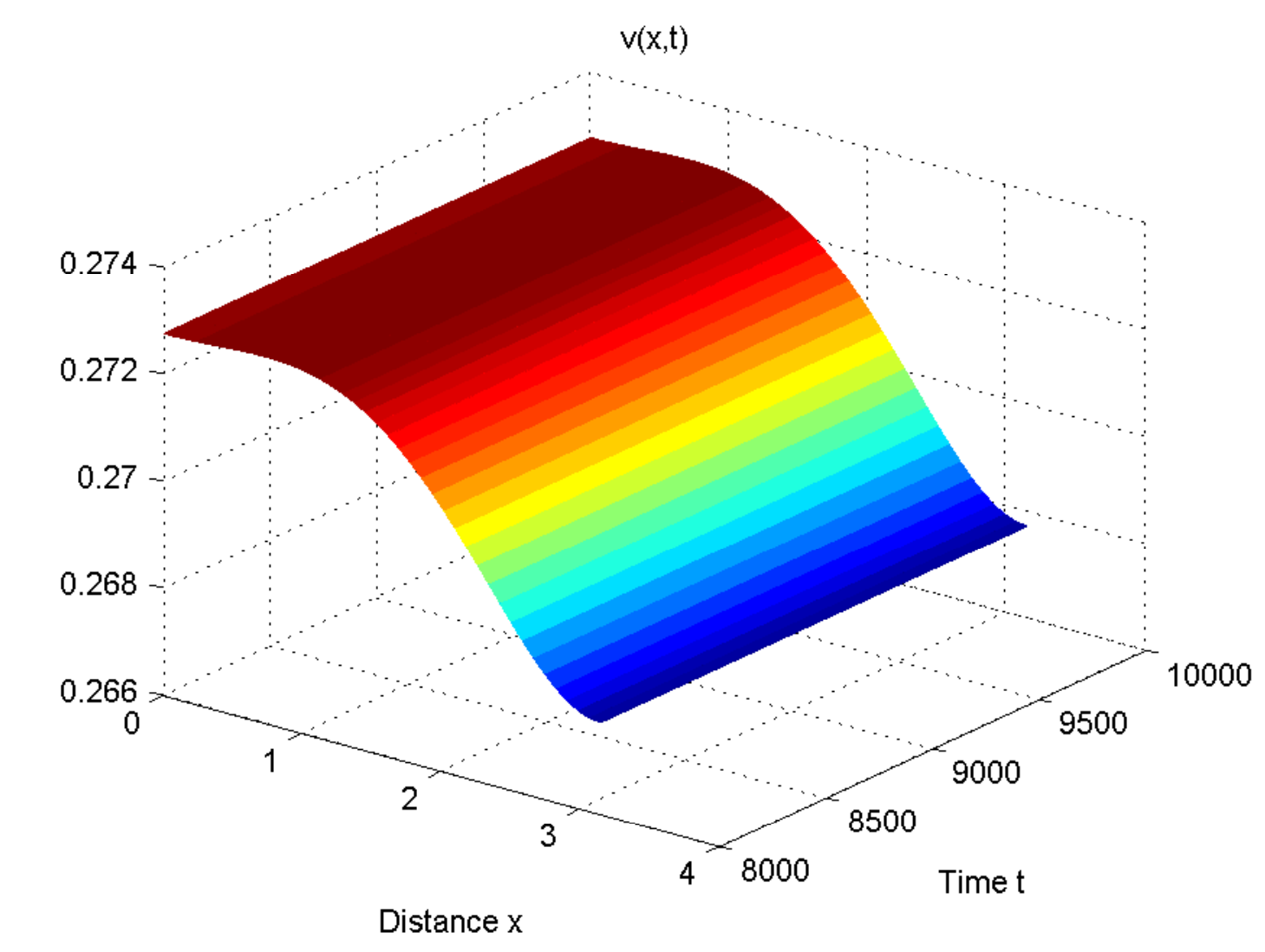}
\end{minipage}
}
\subfigure[The initial values are $u(0,x)=0.2716+0.02\cos 2x, v(0,x)=0.2716+0.05\cos 2x$.]{
\label{fig:7:c} 
\begin{minipage}[b]{0.48\textwidth}
  \centering\includegraphics[width=0.49\textwidth]{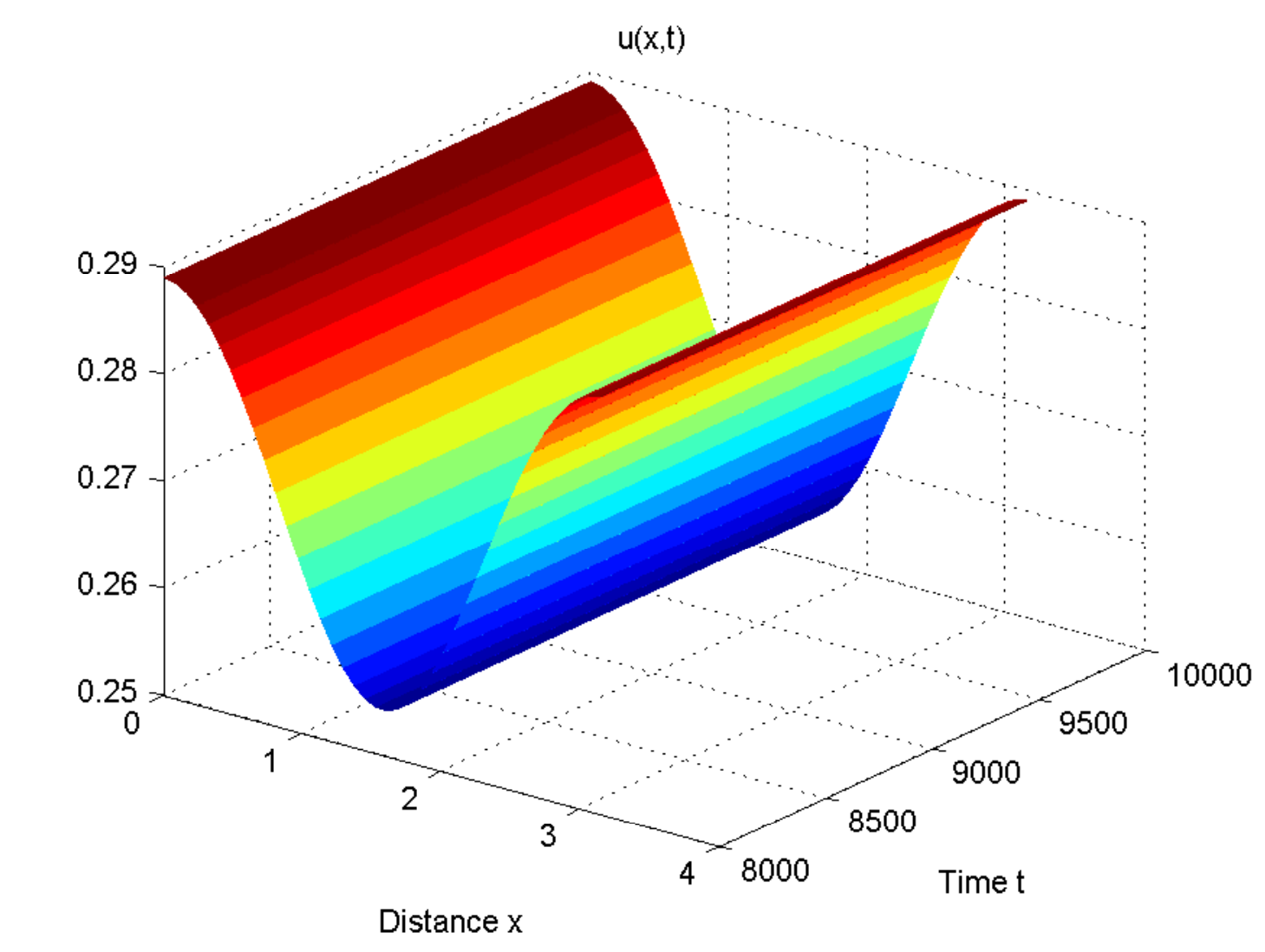}\includegraphics[width=0.49\textwidth]{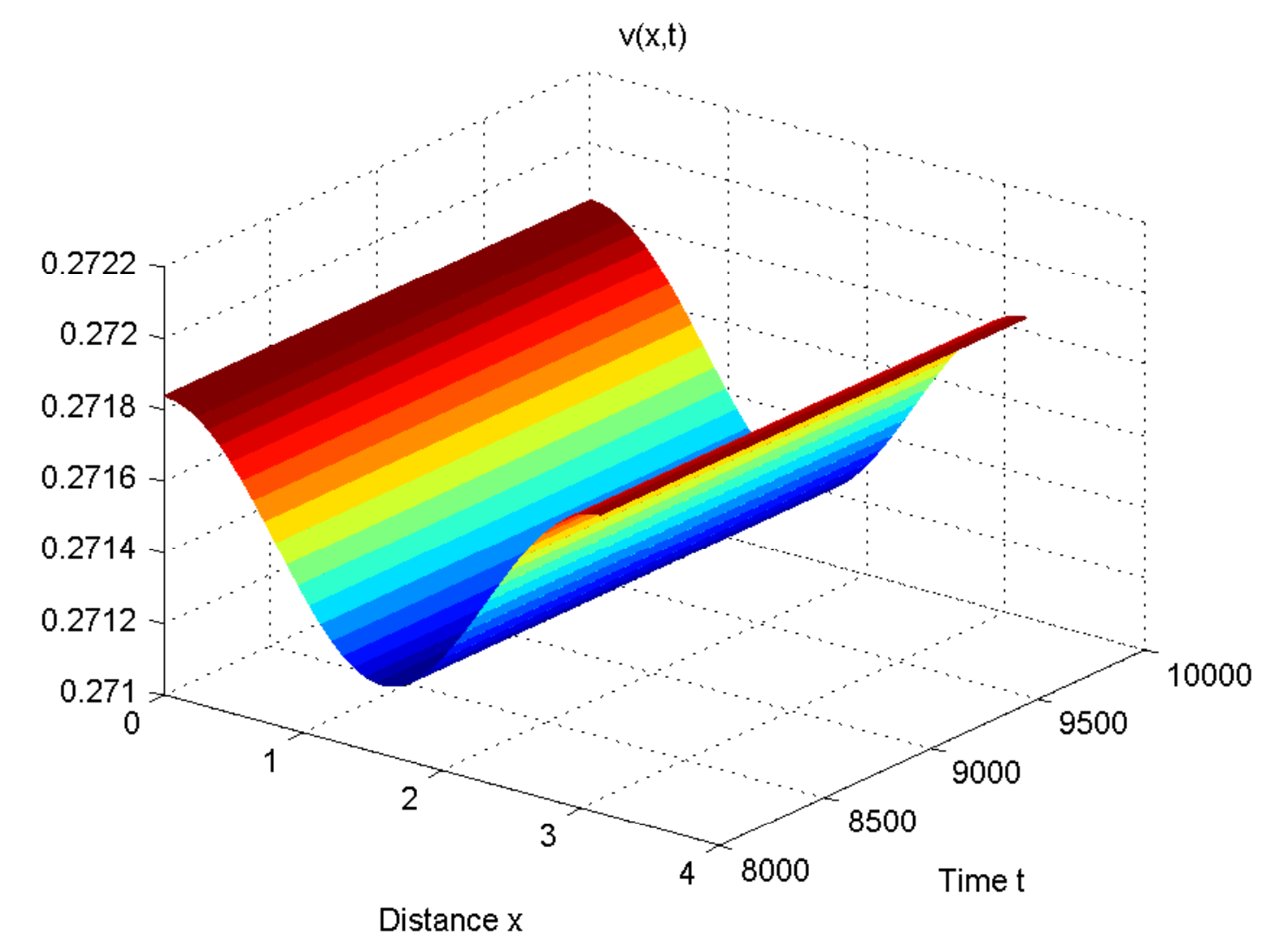}
\end{minipage}
}
\subfigure[The initial values are $u(0,x)=0.2716-0.02\cos 2x$, $v(0,x)=0.2716-0.05\cos 2x$.]{
\label{fig:7:d} 
\begin{minipage}[b]{0.48\textwidth}
  \centering\includegraphics[width=0.49\textwidth]{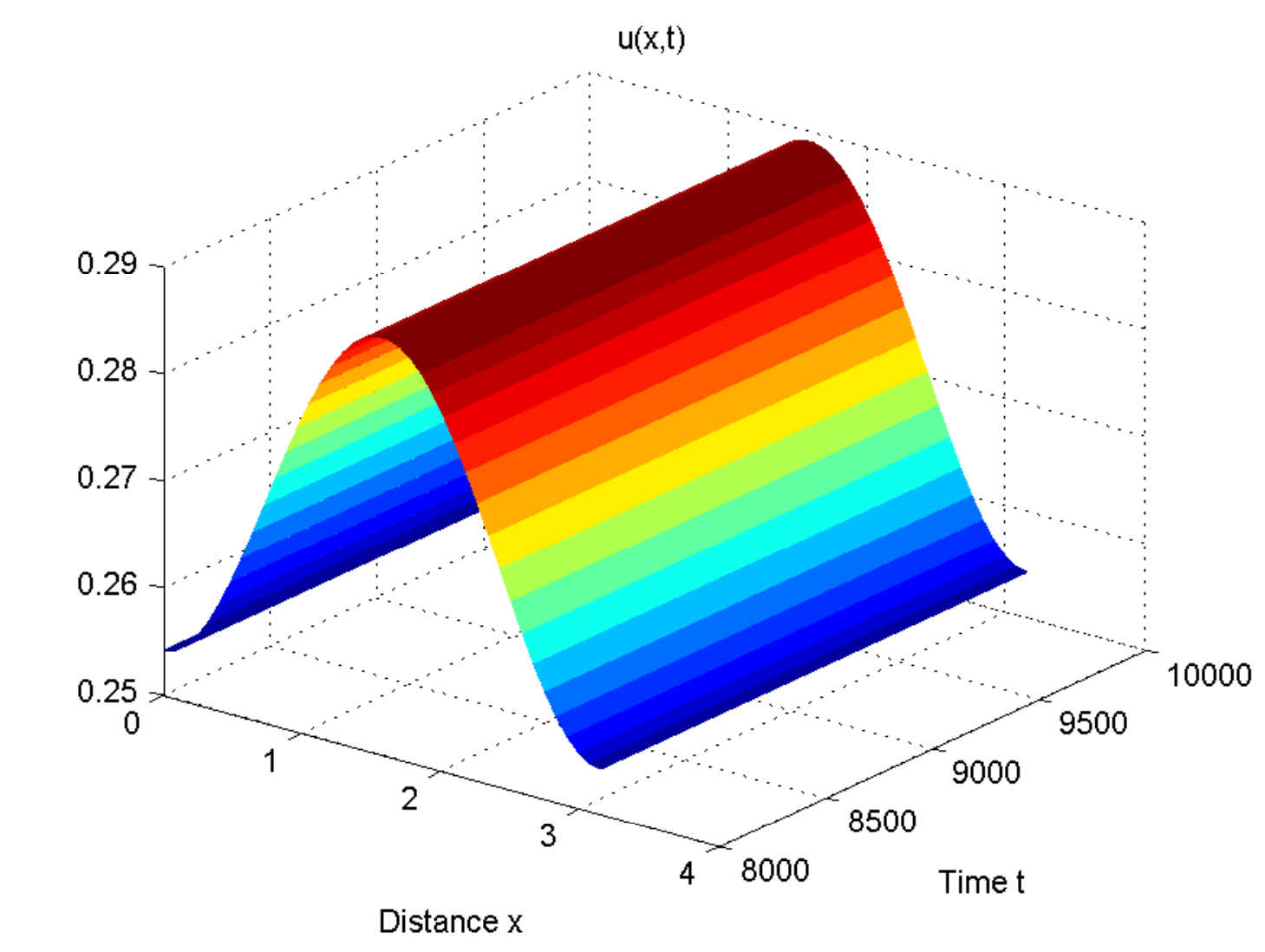}\includegraphics[width=0.49\textwidth]{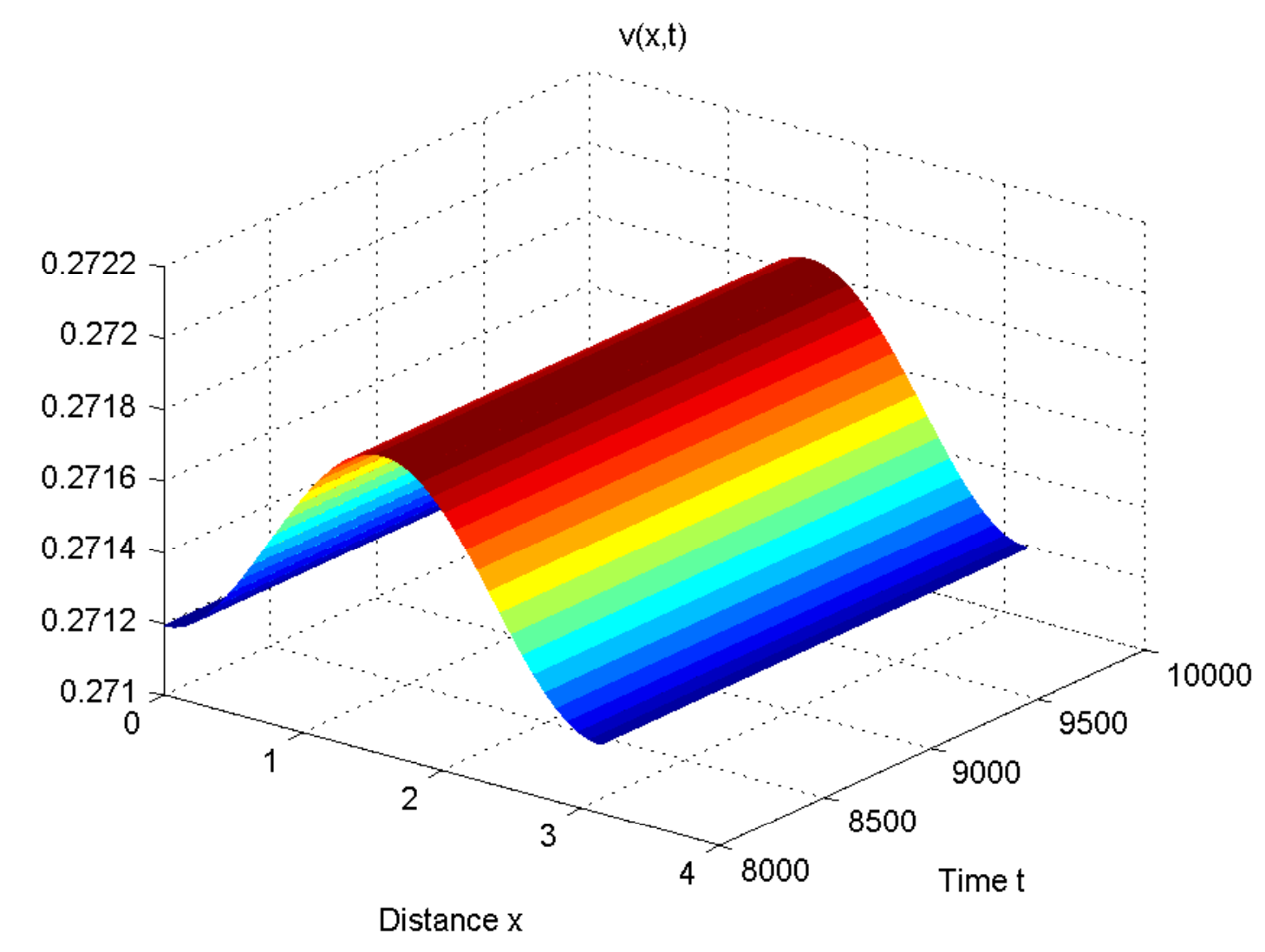}
\end{minipage}
}
\caption{For $\left(d_1,s\right)=(0.01045,0.3029)\in \mathcal{D}_2$, a pair of stable superposition steady states with the shape of $\phi_1h_1\cos x+\phi_2h_2\cos2x-$like and a pair of stable spatially inhomogeneous steady states with the shape of $\phi_2\cos2x-$like, coexist.}
\label{fig:7} 
\end{figure}

\begin{figure}[!hbt]
\centering
\subfigure[The initial values are $u(0,x)=0.2716-0.1\cos x, v(0,x)=0.2716-0.1\cos x$.]{
\label{fig:8:a} 
\begin{minipage}[b]{0.48\textwidth}
  \centering\includegraphics[width=0.49\textwidth]{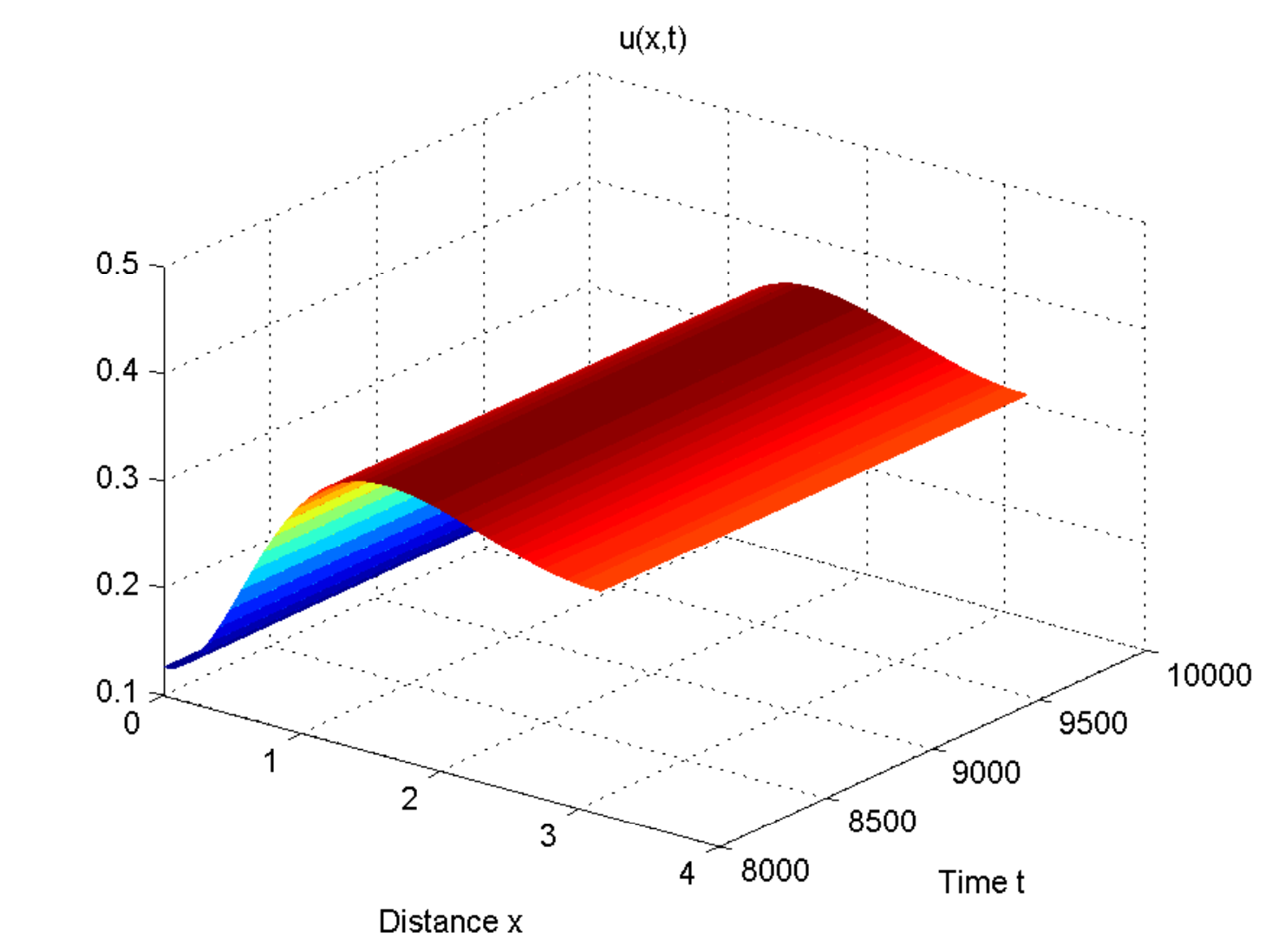}
  \includegraphics[width=0.49\textwidth]{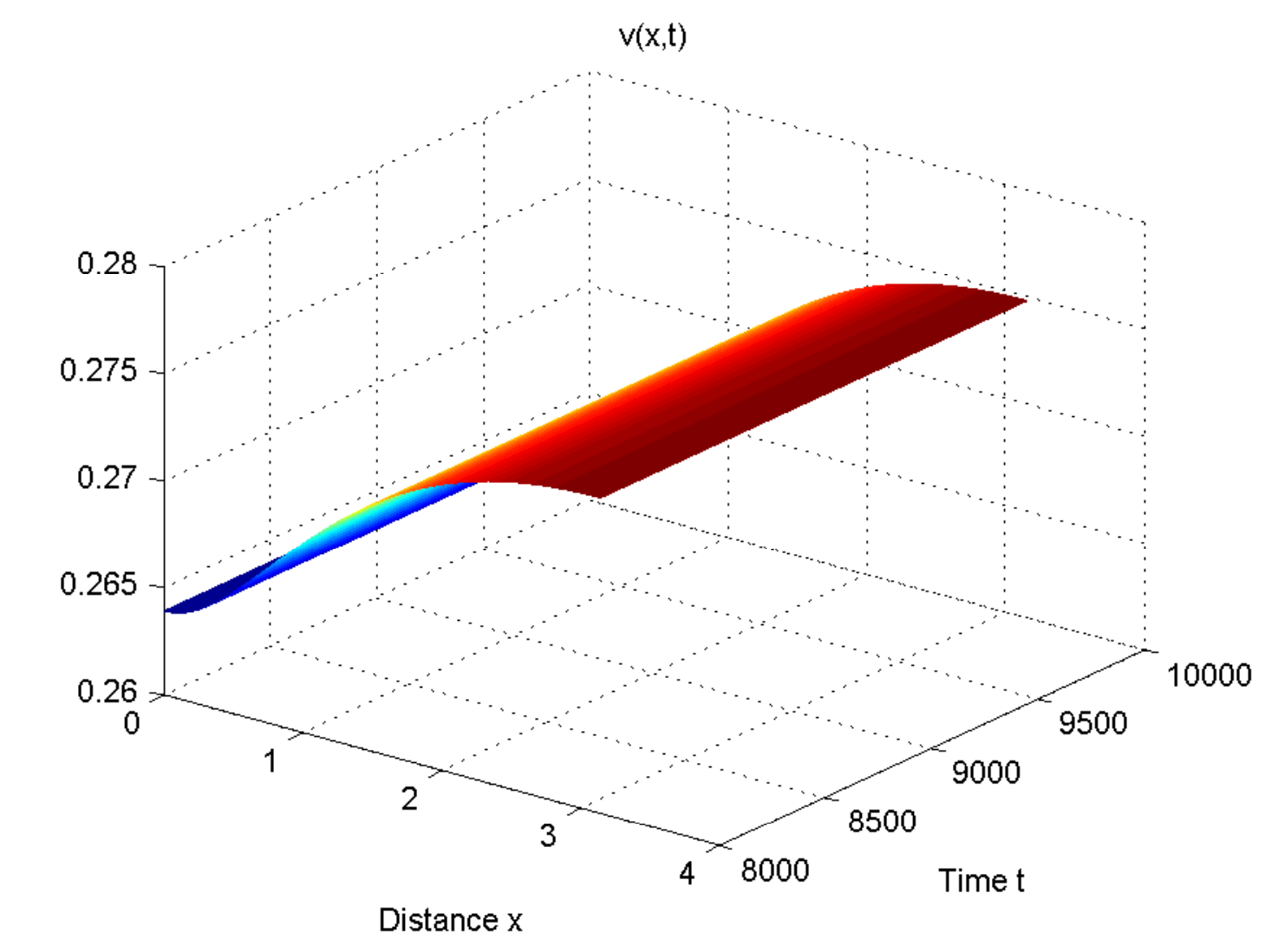}
\end{minipage}
}
\subfigure[The initial values are $u(0,x)=0.2716+0.1\cos x, v(0,x)=0.2716+0.1\cos x$.]{
\label{fig:8:b} 
\begin{minipage}[b]{0.48\textwidth}
  \centering\includegraphics[width=0.49\textwidth]{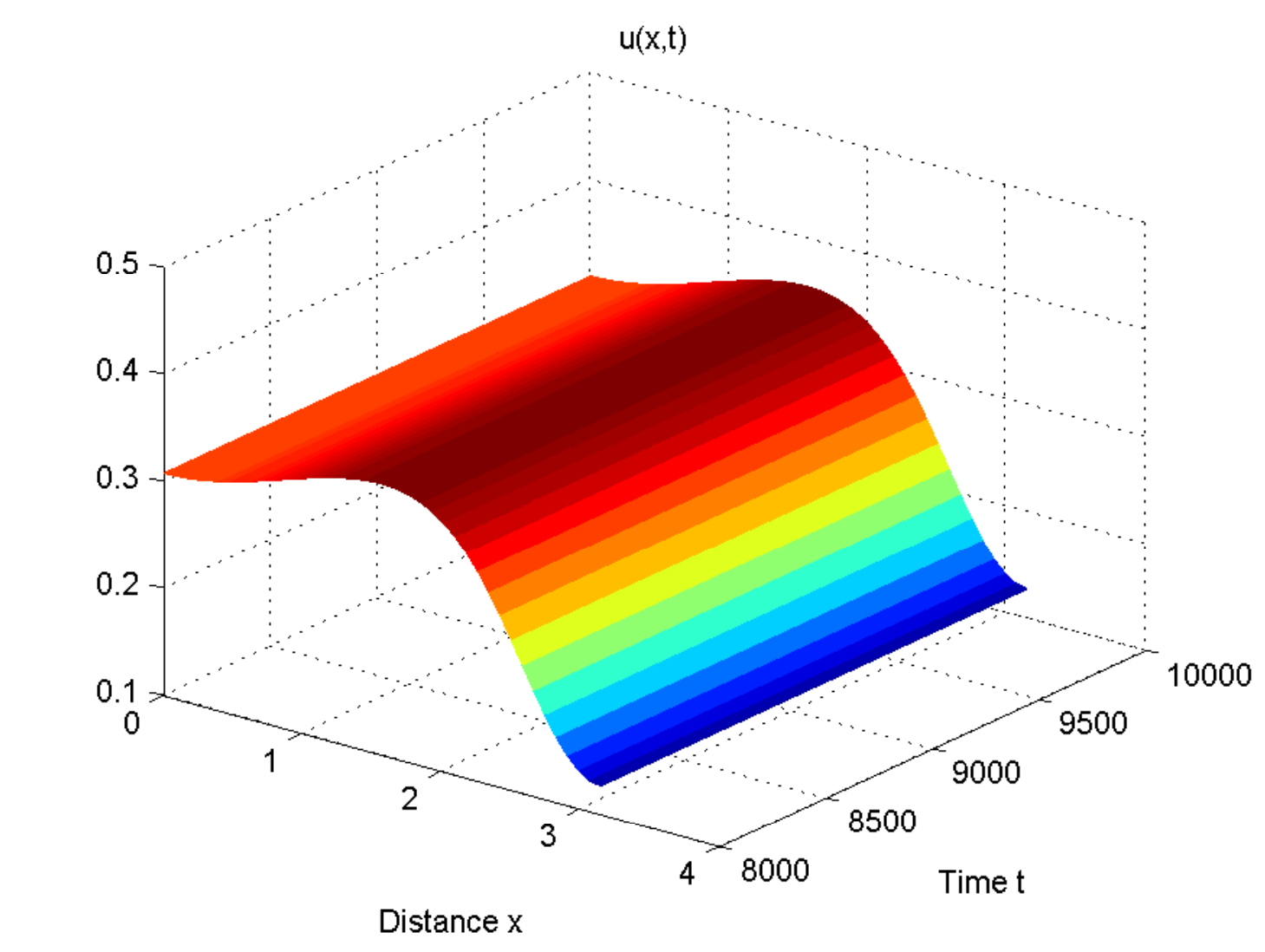}
  \includegraphics[width=0.49\textwidth]{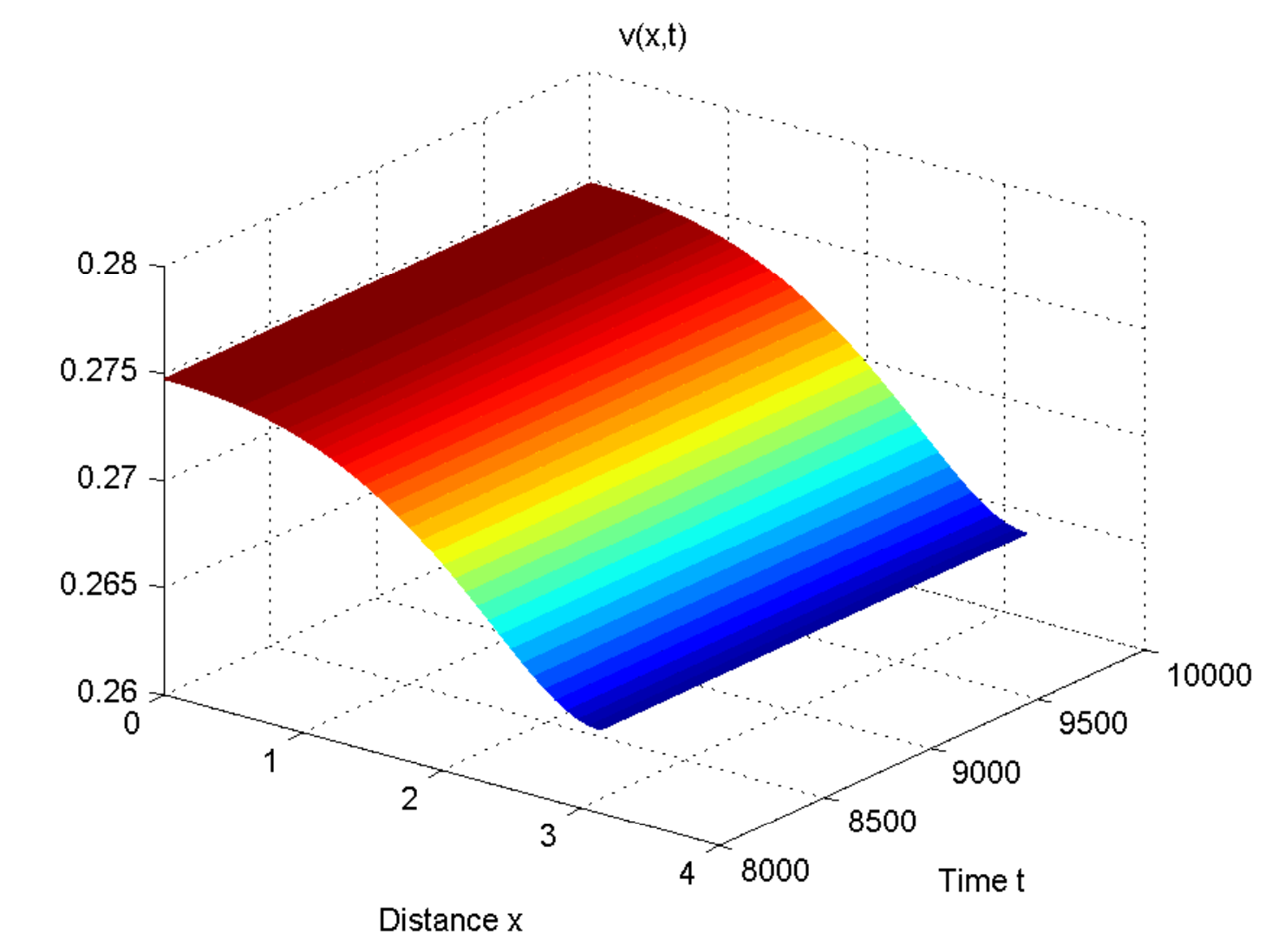}
\end{minipage}
}
\subfigure[The initial values are $u(0,x)=0.2716+0.02\cos 2x, v(0,x)=0.2716+0.05\cos 2x$.]{
\label{fig:8:c} 
\begin{minipage}[b]{0.48\textwidth}
  \centering\includegraphics[width=0.49\textwidth]{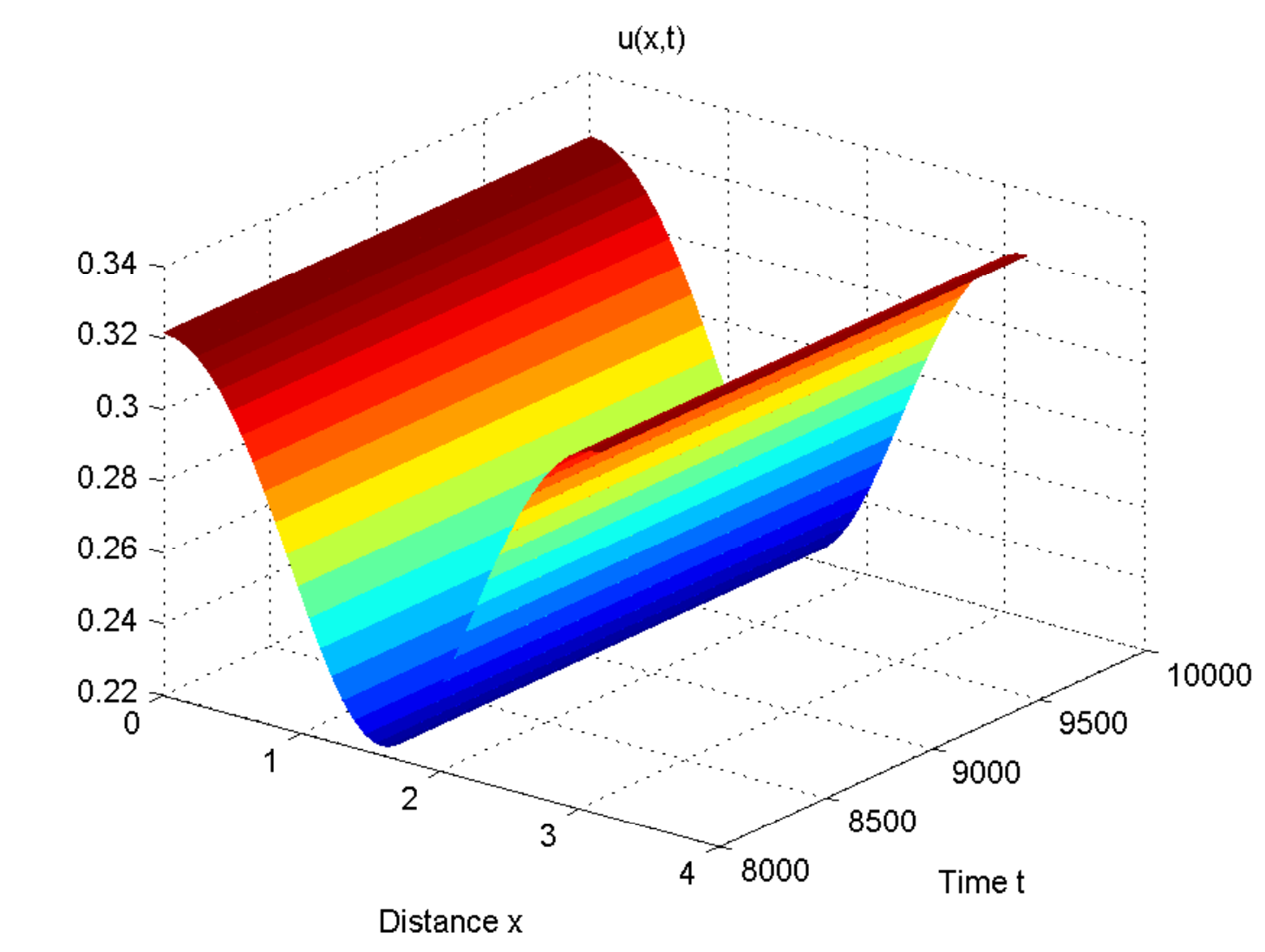}
  \includegraphics[width=0.49\textwidth]{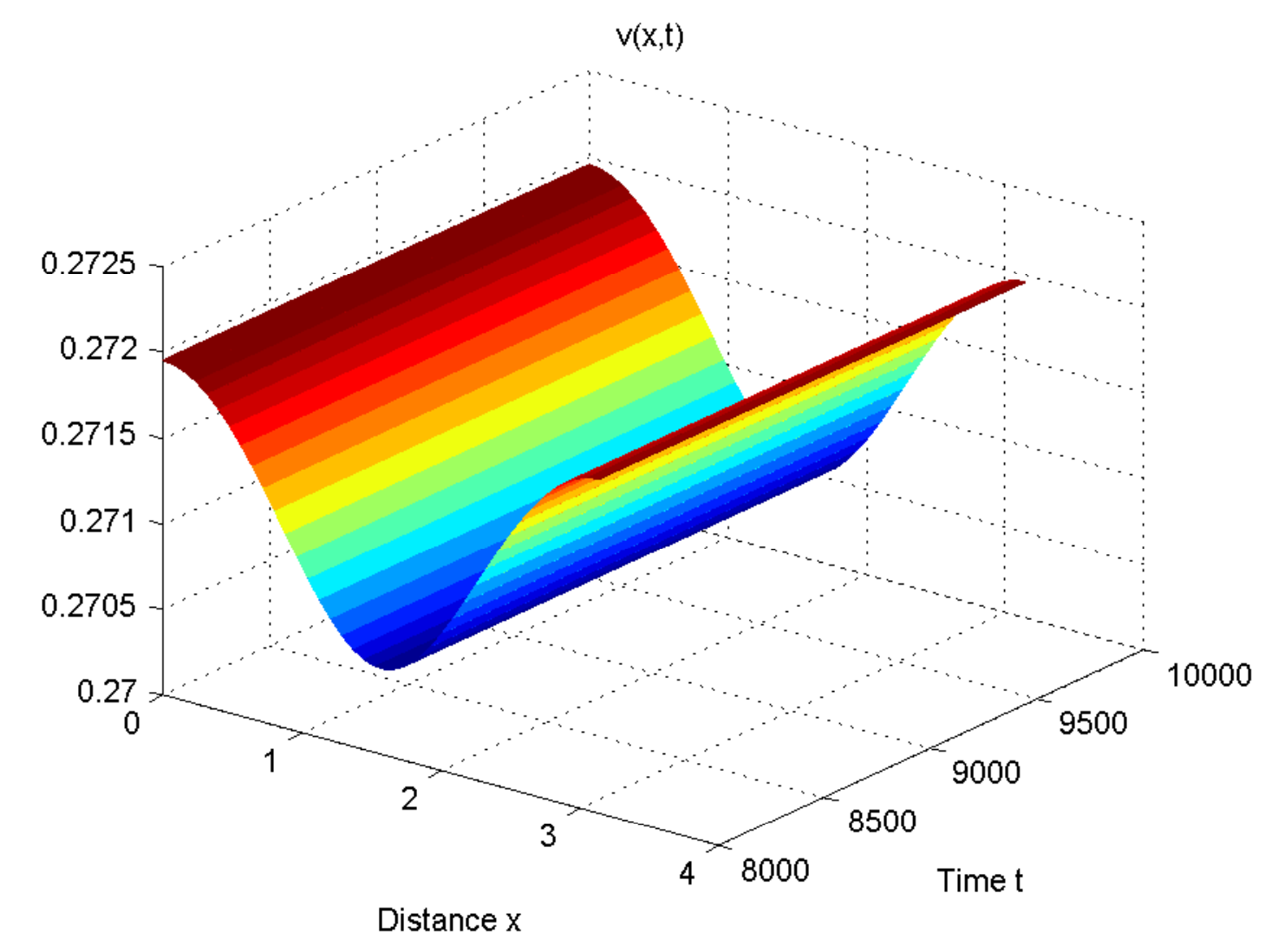}
\end{minipage}
}
\subfigure[The initial values are $u(0,x)=0.2716-0.02\cos 2x$, $v(0,x)=0.2716-0.05\cos 2x$.]{
\label{fig:8:d} 
\begin{minipage}[b]{0.48\textwidth}
  \centering\includegraphics[width=0.49\textwidth]{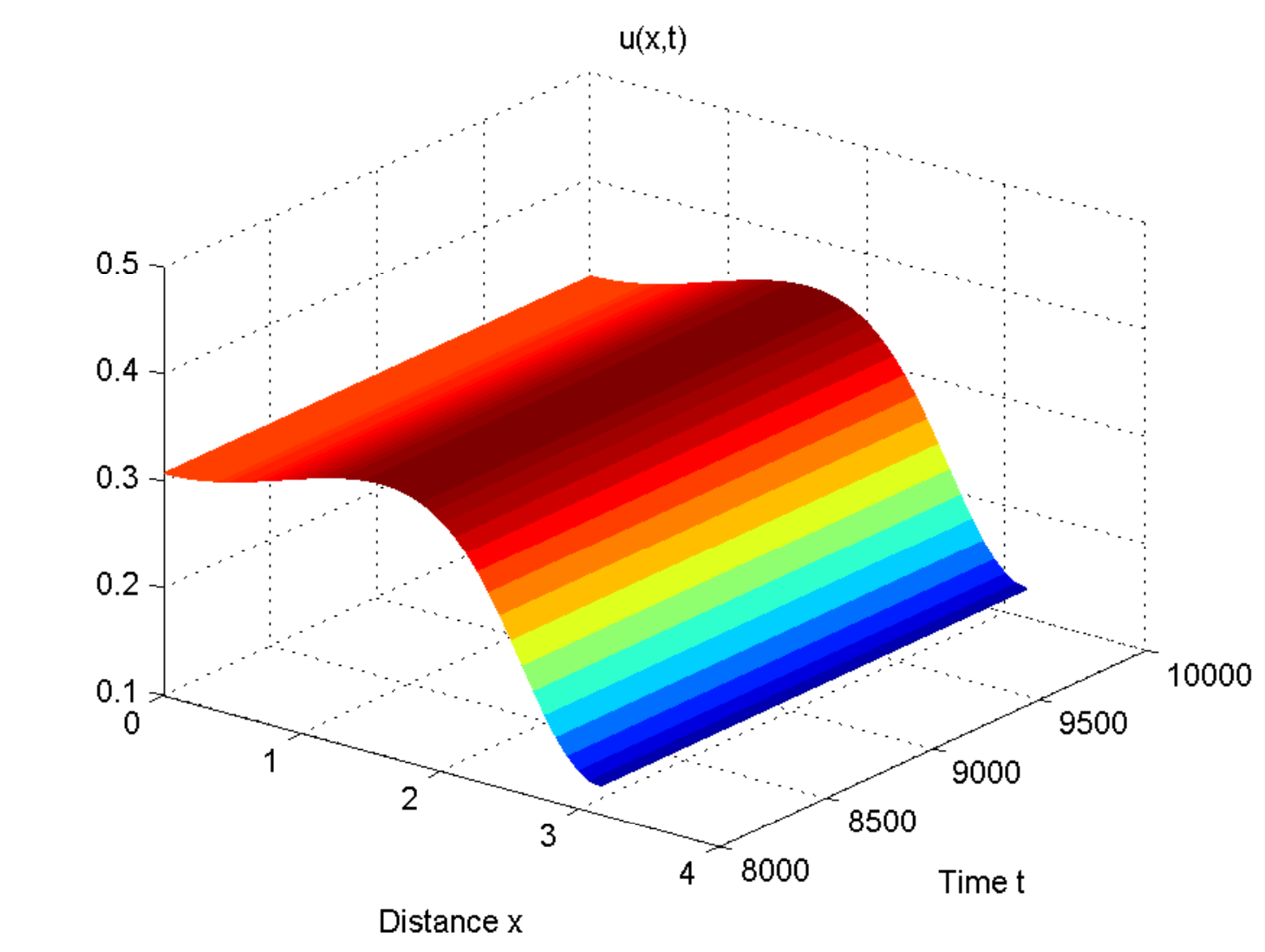}
  \includegraphics[width=0.49\textwidth]{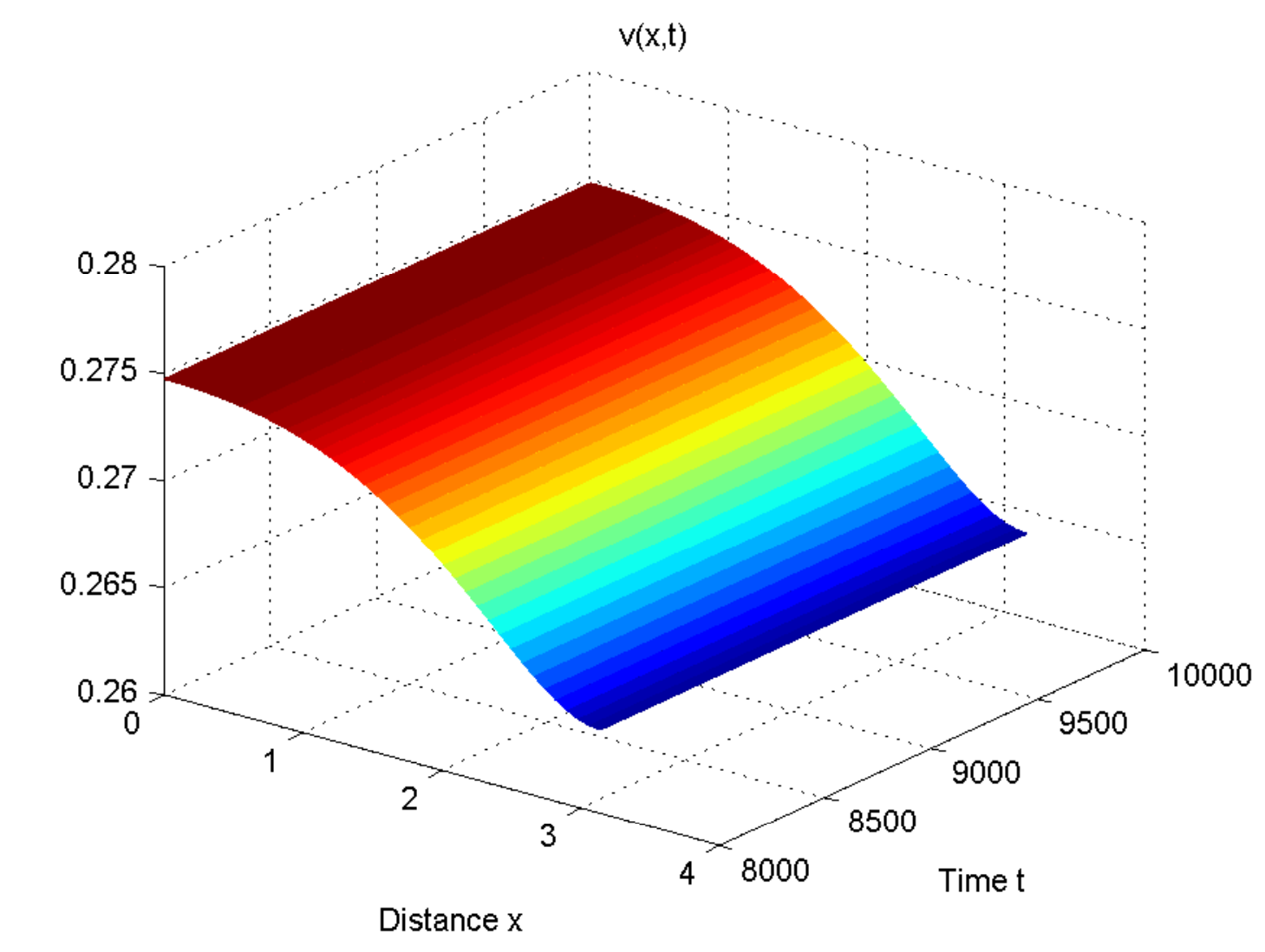}
\end{minipage}
}
\caption{For $\left(d_1,s\right)=(0.01045,0.2379)\in \mathcal{D}_4$, a pair of stable superposition steady states with the shape of $\phi_1h_1\cos x+\phi_2h_2\cos2x-$like and a stable spatially inhomogeneous steady state with the shape of $\phi_2\cos2x-$like, coexist. In (d), the solution starting from initial values $u(0,x)=0.2716-0.02\cos 2x, v(0,x)=0.2716-0.05\cos 2x$ eventually tends to one of the superposition steady states with the shape of $\phi_1h_1\cos x+\phi_2h_2\cos2x-$like, which indicates that one of the pair of spatially inhomogeneous steady states with the shape of $\phi_2\cos2x-$like is unstable.}
\label{fig:8} 
\end{figure}

Then, according to local bifurcation set and the corresponding phase portraits of normal form \eqref{e:3.8}, we summarize partial dynamical behaviors of diffusive system \eqref{e:3.1} as follows,
\begin{proposition}\label{pro:3.4}
For fixed $m=5,a=3,b=0.1,d_2=4$, diffusive predator-prey system \eqref{e:3.1} admits complex spatial patterns, when parameters $(d_1,s)$ are chosen near $(1,2)-$mode Turing-Turing bifurcation point $\left(d_{1,2}^*,s_{1,2}^*\right)=(0.01095,0.2679)$. Here are the partial results:
\begin{enumerate}
\item For $\left(d_1,s\right)\in \mathcal{D}_1$, the coexistence equilibrium $E_*$ and a pair of superposition steady states with the shape of $\phi_1h_1\cos x+\phi_2h_2\cos2x-$like are stable, which indicates that system \eqref{e:3.1} admits \emph{tri-stability} (see Fig. \ref{fig:6}). Also, system \eqref{e:3.1} has another pair of unstable superposition steady states with the shape of $\phi_1h_3\cos x+\phi_2h_4\cos2x-$like. And, $(h_1,h_2)$ and $(h_3,h_4)$ are coexistence equilibria of normal form \eqref{e:3.8}, and $\phi_1, \phi_2$ are given in \eqref{e:3.6-1}.
    Moreover, $E_*$ is unstable for $\left(d_1,s\right)\notin{\mathcal{D}}_1$, while the pair of superposition steady states with the shape of $\phi_1h_1\cos x+\phi_2h_2\cos2x-$like is always stable.
\item For $\left(d_1,s\right)\in \mathcal{D}_2$, system \eqref{e:3.1} admits a pair of stable superposition steady states with the shape of $\phi_1h_1\cos x+\phi_2h_2\cos2x-$like, and a pair of stable spatially inhomogeneous steady states with the shape of $\phi_2\cos2x-$like bifurcating from $E_*$ through Turing bifurcation. Then, system \eqref{e:3.1} supports the \emph{coexistence of four stable spatially inhomogeneous steady states} (see Fig. \ref{fig:7}). Moreover, the pair of unstable superposition steady states with the shape of $\phi_1h_3\cos x+\phi_2h_4\cos2x-$like persists.
\item When $\left(d_1,s\right)$ crosses from $\mathcal{D}_2$ into $\mathcal{D}_3$, one of the pair of  stable spatially inhomogeneous steady state with the shape of $\phi_2\cos2x-$like becomes unstable by colliding with the pair of unstable superposition steady states with the shape of $\phi_1h_3\cos x+\phi_2h_4\cos2x-$like through Turing bifurcation. Therefore, system \eqref{e:3.1} exhibits \emph{tri-stability} that the pair of stable superposition steady states with the shape of $\phi_1h_1\cos x+\phi_2h_2\cos2x-$like and another stable spatially inhomogeneous steady state with the shape of $\phi_2\cos2x-$like coexist, which is slightly different from the case $\left(d_1,s\right)\in \mathcal{D}_1$.
\item When $\left(d_1,s\right)$ crosses from $\mathcal{D}_3$ into $\mathcal{D}_4$, the pair of stable superposition steady states with the shape of $\phi_1h_1\cos x+\phi_2h_2\cos2x-$like and the stable spatially inhomogeneous steady state with the shape of $\phi_2\cos2x-$like, as well as the unstable spatially inhomogeneous steady state with the shape of $\phi_2\cos2x-$like, persist. Therefore, system \eqref{e:3.1} supports \emph{tri-stability} (see Fig. \ref{fig:8}). Moreover, a pair of new unstable superposition steady states with the shape of $\phi_1h_5\cos x+\phi_2h_6\cos2x-$like bifurcates from the unstable coexistence equilibrium $E_*$ through Turing bifurcation, where $(h_5,h_6)$ is also the coexistence equilibrium of normal form \eqref{e:3.8}.
\item When $\left(d_1,s\right)$ crosses from $\mathcal{D}_4$ into $\mathcal{D}_5$, the stable spatially inhomogeneous steady state with the shape of $\phi_2\cos2x-$like loses its stability by colliding with the pair of unstable superposition steady states with the shape of $\phi_1h_5\cos x+\phi_2h_6\cos2x-$like through Turing bifurcation.  Hence, system \eqref{e:3.1} admits \emph{bi-stability}, considering that the pair of stable superposition steady states with the shape of $\phi_1h_1\cos x+\phi_2h_2\cos2x-$like and the unstable spatially inhomogeneous steady state with the shape of $\phi_2\cos2x-$like, still remain.
\item For $\left(d_1,s\right)\in \mathcal{D}_6$, system \eqref{e:3.1} exhibits \emph{bi-stability}, due to the persistence of the pair of stable superposition steady states with the shape of $\phi_1h_1\cos x+\phi_2h_2\cos2x-$like. Moreover, the pair of unstable spatially inhomogeneous steady states with the shape of $\phi_2\cos2x-$like and the unstable coexistence equilibrium $E_*$ collide through Turing bifurcation, when $\left(d_1,s\right)$ crosses from $\mathcal{D}_5$ into $\mathcal{D}_6$.
\item When $\left(d_1,s\right)$ crosses from $\mathcal{D}_6$ back into $\mathcal{D}_1$, the pair of unstable superposition steady states with the shape of $\phi_1h_3\cos x+\phi_2h_4\cos2x-$like bifurcates from the unstable coexistence equilibrium $E_*$ through Turing bifurcation. And, the coexistence equilibrium $E_*$ becomes stable. Also, the pair of stable superposition steady states with the shape of $\phi_1h_1\cos x+\phi_2h_2\cos2x-$like persists.
\end{enumerate}
\end{proposition}

According to Proposition \ref{pro:3.4}, besides quad-stability that a pair of stable superposition steady states with the shape of $\phi_1h_1\cos x+\phi_2h_2\cos2x-$like and a pair of stable spatially inhomogeneous steady states with the shape of $\phi_2\cos2x-$like coexist, predator-prey system \eqref{e:3.1} also admits other complex spatial dynamics, like bistable superposition patterns, and tri-stable patterns that a pair of stable superposition steady states with the shape of $\phi_1h_1\cos x+\phi_2h_2\cos2x-$like coexists with the stable coexistence equilibrium $E_*$ or a stable spatially inhomogeneous steady state with the shape of $\phi_2\cos2x-$like.

\section{Proof of Theorem \ref{thm:2.1}}
In this section, we prove {\bf Theorem \ref{thm:2.1}}.
Due to $B=\mathrm{diag}(0,0)$, we have
\begin{equation}\label{e:4.1}
\left(M_j^1U_j^1\right)(z,\epsilon)=D_zU_j^1(z,\epsilon)Bz-BU_j^1(z,\epsilon)=\begin{pmatrix}0\\0\end{pmatrix},\quad U_j^1(z,\epsilon)=\begin{pmatrix}p_j^1(z,\epsilon)\\[6pt] p_j^2(z,\epsilon)\end{pmatrix},j\in \mathbb{N}, j\ge 2,
\end{equation}
where $p_j^1(z,\epsilon),p_j^2(z,\epsilon)$ are $j-$order polynomials in $(z,\epsilon)$. Thus, $\mathrm{Im}\left(M_j^1\right)=\{0\}$ for $j\ge 2,j\in\mathbb{N}$.

Considering that
\begin{align*}
V_j^{2+2}\left(\mathbb{R}^2\right)=\mathrm{span}\left\{\begin{pmatrix}z^q\epsilon^s\\0\end{pmatrix},\begin{pmatrix}0\\z^q\epsilon^s\end{pmatrix};|q|+|s|=j\in\mathbb{N},q,s\in\mathbb{N}_0^2\right\},
\end{align*}
and according to the decomposition $V_j^{2+2}\left(\mathbb{R}^2\right)=\mathrm{Im}\left(M_j^1\right)\oplus \mathrm{Im}\left(M_j^1\right)^c, j\ge2$, we have $\mathrm{Im}\left(M_j^1\right)^c=V_j^{2+2}\left(\mathbb{R}^2\right),j\ge2$.

Then, normal form of Turing-Turing bifurcation has the following form
\begin{equation}\label{e:4.2}
z=Bz+\frac{1}{2}g_2^1(z,0,\epsilon)+h.o.t.,
\end{equation}
where $g_2^1(z,0,\epsilon)$ is the quadratic function in $(z,\epsilon)$, and $h.o.t.$ stands for higher-order terms.

By $V_2^{2+2}\left(\mathbb{R}^2\right)$ and $\mathrm{Im}\left(M_2^1\right)^c$, write $f_2(z,0,\epsilon)$ and $g_2^1(z,0,\epsilon)$ as,
\begin{equation}\label{e:4.3}
\begin{aligned}
\frac{1}{2!}f_2(z,0,\epsilon)=\sum_{|q|+|s|=2} \frac{1}{\prod_{i=1}^2 q_i!\prod_{k=1}^2 s_k!}f_{qs}z^q\epsilon^s,\qquad
\frac{1}{2!}g_2^1(z,0,\epsilon)=\sum_{|q|+|s|=2} \frac{1}{\prod_{i=1}^2 q_i!\prod_{k=1}^2 s_k!}g^1_{qs}z^q\epsilon^s,
\end{aligned}
\end{equation}
where
\begin{equation*}
\begin{gathered}
q=(q_1,q_2)\in\mathbb{N}_0^2,s=(s_1,s_2)\in\mathbb{N}_0^2,z=(z_1,z_2),\epsilon=(\epsilon_1,\epsilon_2),z^q=z_1^{q_1}z_2^{q_2}, \epsilon^s=\epsilon_1^{s_1}\epsilon_2^{s_2},|q|=q_1+q_2,\\
|s|=s_1+s_2,f_{qs}\triangleq f_{q_1q_2s_1s_2}, f_{qs}=\left(f_{qs}^{1},f_{qs}^{2}\right)^T, f_{qs}^1=\left(f_{qs}^{11},f_{qs}^{12}\right)^T, g_{qs}^1\triangleq g^1_{q_1q_2s_1s_2}, g_{qs}^1=\left(g_{qs}^{11},g_{qs}^{12}\right)^T,
\end{gathered}
\end{equation*}
and polynomial coefficients $\left\{g^1_{qs}:|q|+|s|=2,\;q,s\in \mathbb{N}_0^2\right\}$ need to be determined.

Then, by $g^1_2=f^1_2-M^1_2U^1_2$ and \eqref{e:4.1}, a direct calculation yields
\begin{equation}\label{e:4.4}
g^1_{qs}=f^1_{qs},\quad |q|+|s|=2.
\end{equation}
And by Eqs.\eqref{e:4.3}, \eqref{e:2.14}, \eqref{e:2.3} and expansions \eqref{e:2.20}, we have
\begin{equation}\label{e:4.5}
\begin{aligned}
&f_{2000}^{11}=\psi_1(0)Q(\phi_1,\phi_1)\left\langle\beta_{k_1}^2,\beta_{k_1}\right\rangle;
&\qquad &f_{0200}^{11}=\psi_1(0)Q(\phi_2,\phi_2)\left\langle\beta_{k_2}^2,\beta_{k_1}\right\rangle;\\
&f_{1100}^{11}=\psi_1(0)Q(\phi_1,\phi_2)\left\langle\beta_{k_1}\beta_{k_2},\beta_{k_1}\right\rangle;
& &f_{2000}^{12}=\psi_2(0)Q(\phi_1,\phi_1)\left\langle\beta_{k_1}^2,\beta_{k_2}\right\rangle;\\
&f_{0200}^{12}=\psi_2(0)Q(\phi_2,\phi_2)\left\langle\beta_{k_2}^2,\beta_{k_2}\right\rangle;
& &f_{1100}^{12}=\psi_2(0)Q(\phi_1,\phi_2)\left\langle\beta_{k_1}\beta_{k_2},\beta_{k_2}\right\rangle;\\
&f_{1010}^{11}=\frac{1}{2}\psi_1(0)\left(L_{\epsilon_1}\phi_1-\mu_{k_1}D_{\epsilon_1}\phi_1(0)\right)\left\langle\beta_{k_1},\beta_{k_1}\right\rangle;
& &f_{1001}^{11}=\frac{1}{2}\psi_1(0)\left(L_{\epsilon_2}\phi_1-\mu_{k_1}D_{\epsilon_2}\phi_1(0)\right)\left\langle\beta_{k_1},\beta_{k_1}\right\rangle;\\
&f_{0110}^{11}=\frac{1}{2}\psi_1(0)\left(L_{\epsilon_1}\phi_2-\mu_{k_2}D_{\epsilon_1}\phi_2(0)\right)\left\langle\beta_{k_2},\beta_{k_1}\right\rangle;
& &f_{0101}^{11}=\frac{1}{2}\psi_1(0)\left(L_{\epsilon_2}\phi_2-\mu_{k_2}D_{\epsilon_2}\phi_2(0)\right)\left\langle\beta_{k_2},\beta_{k_1}\right\rangle;\\
&f_{1010}^{12}=\frac{1}{2}\psi_2(0)\left(L_{\epsilon_1}\phi_1-\mu_{k_1}D_{\epsilon_1}\phi_1(0)\right)\left\langle\beta_{k_1},\beta_{k_2}\right\rangle;
& &f_{1001}^{12}=\frac{1}{2}\psi_2(0)\left(L_{\epsilon_2}\phi_1-\mu_{k_1}D_{\epsilon_2}\phi_1(0)\right)\left\langle\beta_{k_1},\beta_{k_2}\right\rangle;\\
&f_{0110}^{12}=\frac{1}{2}\psi_2(0)\left(L_{\epsilon_1}\phi_2-\mu_{k_2}D_{\epsilon_1}\phi_2(0)\right)\left\langle\beta_{k_2},\beta_{k_2}\right\rangle;
& &f_{0101}^{12}=\frac{1}{2}\psi_2(0)\left(L_{\epsilon_2}\phi_2-\mu_{k_2}D_{\epsilon_2}\phi_2(0)\right)\left\langle\beta_{k_2},\beta_{k_2}\right\rangle.
\end{aligned}
\end{equation}

Plugging \eqref{e:4.5} into \eqref{e:4.4}, we determine polynomial coefficients $\left\{g^1_{qs}:|q|+|s|=2,\;q,s\in \mathbb{N}_0^2\right\}$, then we derive $g_2^1(z,0,\epsilon)$.
\begin{remark}
If the second-order normal form is non-degenerated, the local dynamical properties of original PFDEs could be described by the normal form truncated to order 2.
\end{remark}

According to Eq.  \eqref{e:2.19}, we have
\begin{equation*}
U_2(z,\epsilon)=\left(U_2^1,U_2^2\right)^T=\left(M_2\right)^{-1}\mathbf{P}_{\mathrm{Im}\left(M_2^1\right)\times \mathrm{Im}\left(M_2^2\right)}\circ\tilde{f}_2(z,0,\epsilon).
\end{equation*}
Then $U_2^1\in \mathrm{Ker}\left(M_2^1\right)^c=\{0\}$, by $V_2^{2+2}\left(\mathbb{R}^2\right)=\mathrm{Ker}\left(M_2^1\right)\oplus \mathrm{Ker}\left(M_2^1\right)^c$ and $\mathrm{Ker}\left(M_2^1\right)=V_2^{2+2}\left(\mathbb{R}^2\right)$.
Thus,
\begin{equation}\label{e:4.6}
U_2^1(z,\epsilon)=\begin{pmatrix}0\\0\end{pmatrix}.
\end{equation}

If the second-order normal form is degenerated, we further calculate the higher-order normal form. Firstly, we calculate center manifold $y=h(z,\epsilon)\triangleq\sum_{j\ge 2}\frac{1}{j!}U_j^2(z,\epsilon)$. Similar to \eqref{e:4.3}, $g_2^2(z,0,\epsilon)$ has the following form,
\begin{equation}\label{e:4.7}
g_2^2(z,0,\epsilon)=0.
\end{equation}
Also, write $h(z,\epsilon)$ as
\begin{equation}\label{e:4.8}
\begin{aligned}
h(z,\epsilon)=\sum_{j\ge2}\frac{1}{j!}h_j(z,\epsilon)=\sum_{j\ge2}\sum_{|q|+|s|=j} \frac{1}{\prod_{i=1}^2 q_i!\prod_{k=1}^2 s_k!}h_{qs}z^q\epsilon^s,\quad j\in\mathbb{N}.
\end{aligned}
\end{equation}
Then, according to $g^2_2=f^2_2-M^2_2U^2_2$, Eqs. \eqref{e:4.7} and \eqref{e:4.3}, a direct calculation yields
\begin{equation}\label{e:4.9}
f^2_{qs}+\mathcal{A}_1(h_{qs})=0,\quad |q|+|s|=2,
\end{equation}
where $\mathcal{A}_1\phi=\dot{\phi}+X_0\left[D_0\Delta\phi(0)+L_0^*\phi-\dot{\phi}(0)\right]$ for $\phi\in\mathcal{Q}_1$.\\
By expansions \eqref{e:2.20}, Eqs. \eqref{e:4.5} and \eqref{e:4.9}, we derive
\begin{equation}\label{e:4.10}
\begin{aligned}
h_{2000}(\theta)=&\theta\phi_1\psi_1(0)Q(\phi_1,\phi_1)\left\langle\beta_{k_1}^2,\beta_{k_1}\right\rangle\beta_{k_1}+\theta\phi_2\psi_2(0)Q(\phi_1,\phi_1)\left\langle\beta_{k_1}^2,\beta_{k_2}\right\rangle\beta_{k_2}+h_{2000}(0);\\
h_{1100}(\theta)=&\theta\phi_1\psi_1(0)Q(\phi_1,\phi_2)\left\langle\beta_{k_1}\beta_{k_2},\beta_{k_1}\right\rangle\beta_{k_1}
+\theta\phi_2\psi_2(0)Q(\phi_1,\phi_2)\left\langle\beta_{k_1}\beta_{k_2},\beta_{k_2}\right\rangle\beta_{k_2}+h_{1100}(0);\\
h_{0200}(\theta)=&\theta\phi_1\psi_1(0)Q(\phi_2,\phi_2)\left\langle\beta_{k_2}^2,\beta_{k_1}\right\rangle\beta_{k_1}
+\theta\phi_2\psi_2(0)Q(\phi_2,\phi_2)\left\langle\beta_{k_2}^2,\beta_{k_2}\right\rangle\beta_{k_2}+h_{0200}(0);\\
h_{1010}(\theta)=&\frac{1}{2}\theta\phi_1\psi_1(0)\left(L_{\epsilon_1}\phi_1-\mu_{k_1}D_{\epsilon_1}\phi_1(0)\right)\left\langle\beta_{k_1},\beta_{k_1}\right\rangle\beta_{k_1}
+\frac{1}{2}\theta\phi_2\psi_2(0)\left(L_{\epsilon_1}\phi_1-\mu_{k_1}D_{\epsilon_1}\phi_1(0)\right)\left\langle\beta_{k_1},\beta_{k_2}\right\rangle\beta_{k_2}+h_{1010}(0);\\
h_{1001}(\theta)=&\frac{1}{2}\theta\phi_1\psi_1(0)\left(L_{\epsilon_2}\phi_1-\mu_{k_1}D_{\epsilon_2}\phi_1(0)\right)\left\langle\beta_{k_1},\beta_{k_1}\right\rangle\beta_{k_1}
+\frac{1}{2}\theta\phi_2\psi_2(0)\left(L_{\epsilon_2}\phi_1-\mu_{k_1}D_{\epsilon_2}\phi_1(0)\right)\left\langle\beta_{k_1},\beta_{k_2}\right\rangle\beta_{k_2}+h_{1001}(0);\\
h_{0110}(\theta)=&\frac{1}{2}\theta\phi_1\psi_1(0)\left(L_{\epsilon_1}\phi_2-\mu_{k_2}D_{\epsilon_1}\phi_2(0)\right)\left\langle\beta_{k_2},\beta_{k_1}\right\rangle\beta_{k_1}
+\frac{1}{2}\theta\phi_2\psi_2(0)\left(L_{\epsilon_1}\phi_2-\mu_{k_2}D_{\epsilon_1}\phi_2(0)\right)\left\langle\beta_{k_2},\beta_{k_2}\right\rangle\beta_{k_2}+h_{0110}(0);\\
h_{0101}(\theta)=&\frac{1}{2}\theta\phi_1\psi_1(0)\left(L_{\epsilon_2}\phi_2-\mu_{k_2}D_{\epsilon_2}\phi_2(0)\right)\left\langle\beta_{k_2},\beta_{k_1}\right\rangle\beta_{k_1}
+\frac{1}{2}\theta\phi_2\psi_2(0)\left(L_{\epsilon_2}\phi_2-\mu_{k_2}D_{\epsilon_2}\phi_2(0)\right)\left\langle\beta_{k_2},\beta_{k_2}\right\rangle\beta_{k_2}+h_{0101}(0);
\end{aligned}
\end{equation}
with
\begin{equation*}
\begin{aligned}
L_0^*h_{2000}(\theta)+D_0\Delta h_{2000}(0)=&-Q(\phi_1,\phi_1)\beta_{k_1}^2+\phi_1(0)\psi_1(0)Q(\phi_1,\phi_1)\left\langle\beta_{k_1}^2,\beta_{k_1}\right\rangle\beta_{k_1}
+\phi_2(0)\psi_2(0)Q(\phi_1,\phi_1)\left\langle\beta_{k_1}^2,\beta_{k_2}\right\rangle\beta_{k_2};\\
L_0^*h_{1100}(\theta)+D_0\Delta h_{1100}(0)=&-Q(\phi_1,\phi_2)\beta_{k_1}\beta_{k_2}+\phi_1(0)\psi_1(0)Q(\phi_1,\phi_2)\left\langle\beta_{k_1}\beta_{k_2},\beta_{k_1}\right\rangle\beta_{k_1}\\
&+\phi_2(0)\psi_2(0)Q(\phi_1,\phi_2)\left\langle\beta_{k_1}\beta_{k_2},\beta_{k_2}\right\rangle\beta_{k_2};\\
L_0^*h_{0200}(\theta)+D_0\Delta h_{0200}(0)=&-Q(\phi_2,\phi_2)\beta_{k_2}^2+\phi_1(0)\psi_1(0)Q(\phi_2,\phi_2)\left\langle\beta_{k_2}^2,\beta_{k_1}\right\rangle\beta_{k_1}
+\phi_2(0)\psi_2(0)Q(\phi_2,\phi_2)\left\langle\beta_{k_2}^2,\beta_{k_2}\right\rangle\beta_{k_2};\\
L_0^*h_{1010}(\theta)+D_0\Delta h_{1010}(0)=&-(L_{\epsilon_1}\phi_1-\mu_{k_1}D_{\epsilon_1}\phi_1(0))\beta_{k_1}
+\phi_1(0)\psi_1(0)\left(L_{\epsilon_1}\phi_1-\mu_{k_1}D_{\epsilon_1}\phi_1(0)\right)\left\langle\beta_{k_1},\beta_{k_1}\right\rangle\beta_{k_1}\\
&+\phi_2(0)\psi_2(0)\left(L_{\epsilon_1}\phi_1-\mu_{k_1}D_{\epsilon_1}\phi_1(0)\right)\left\langle\beta_{k_1},\beta_{k_2}\right\rangle\beta_{k_2};\\
L_0^*h_{1001}(\theta)+D_0\Delta h_{1001}(0)=&-(L_{\epsilon_2}\phi_1-\mu_{k_1}D_{\epsilon_2}\phi_1(0))\beta_{k_1}
+\phi_1(0)\psi_1(0)\left(L_{\epsilon_2}\phi_1-\mu_{k_1}D_{\epsilon_2}\phi_1(0)\right)\left\langle\beta_{k_1},\beta_{k_1}\right\rangle\beta_{k_1}\\
&+\phi_2(0)\psi_2(0)\left(L_{\epsilon_2}\phi_1-\mu_{k_1}D_{\epsilon_2}\phi_1(0)\right)\left\langle\beta_{k_1},\beta_{k_2}\right\rangle\beta_{k_2};\\
L_0^*h_{0110}(\theta)+D_0\Delta h_{0110}(0)=&-(L_{\epsilon_1}\phi_2-\mu_{k_2}D_{\epsilon_1}\phi_2(0))\beta_{k_2}
+\phi_1(0)\psi_1(0)\left(L_{\epsilon_1}\phi_2-\mu_{k_2}D_{\epsilon_1}\phi_2(0)\right)\left\langle\beta_{k_2},\beta_{k_1}\right\rangle\beta_{k_1}\\
&+\phi_2(0)\psi_2(0)\left(L_{\epsilon_1}\phi_2-\mu_{k_2}D_{\epsilon_1}\phi_2(0)\right)\left\langle\beta_{k_2},\beta_{k_2}\right\rangle\beta_{k_2};\\
L_0^*h_{0101}(\theta)+D_0\Delta h_{0101}(0)=&-(L_{\epsilon_2}\phi_2-\mu_{k_2}D_{\epsilon_2}\phi_2(0))\beta_{k_2}
+\phi_1(0)\psi_1(0)\left(L_{\epsilon_2}\phi_2-\mu_{k_2}D_{\epsilon_2}\phi_2(0)\right)\left\langle\beta_{k_2},\beta_{k_1}\right\rangle\beta_{k_1}\\
&+\phi_2(0)\psi_2(0)\left(L_{\epsilon_2}\phi_2-\mu_{k_2}D_{\epsilon_2}\phi_2(0)\right)\left\langle\beta_{k_2},\beta_{k_2}\right\rangle\beta_{k_2}.
\end{aligned}
\end{equation*}
Moreover, by $\mathrm{Im}\left(M_j^1\right)^c=V_j^{2+2}\left(\mathbb{R}^2\right),j\ge2$, the third-order normal form has the form,
\begin{equation}\label{e:4.11}
\dot{z}=Bz+\frac{1}{2!}g_2^1(z,0,\epsilon)+\frac{1}{3!}g_3^1(z,0,\epsilon)+\mathrm{h.o.t.},
\end{equation}
where
\begin{equation}\label{e:4.12}
\frac{1}{3!}g_3^1(z,0,\epsilon)=\frac{1}{3!}\mathbf{P}_{\mathrm{Im}\left(M_3^1\right)^c}\circ\tilde{f}_3^1(z,0,\epsilon)
\end{equation}
with
\begin{equation}\label{e:4.13}
\begin{split}
\frac{1}{3!}\tilde{f}_3^1(z,0,\epsilon)=\sum_{|q|+|s|=3} \frac{1}{\prod_{i=1}^2 q_i!\prod_{k=1}^2 s_k!}\tilde{f}_{qs}^1z^q\epsilon^s,\qquad
\frac{1}{3!}g_3^1(z,0,\epsilon)=\sum_{|q|+|s|=3} \frac{1}{\prod_{i=1}^2 q_i!\prod_{k=1}^2 s_k!}g^1_{qs}z^q\epsilon^s,
\end{split}
\end{equation}
and polynomial coefficients $\left\{g^1_{qs}: |q|+|s|=3,\;q,s\in \mathbb{N}_0^2\right\}$ need to be determined.

By $g^1_j=\tilde{f}^1_j-M^1_jU^1_j,j\ge 3$ and Eq.  \eqref{e:4.1}, a direct calculation yields
\begin{equation}\label{e:4.14}
g^1_{qs}=\tilde{f}^1_{qs},\qquad |q|+|s|=3.
\end{equation}
And considering that $\tilde{f}_3=f_3+\frac{3}{2}\left[\left(Df_2\right)U_2-(DU_2)g_2\right]$,
let $y=0$, then $\tilde{f}_3^1(z,0,\epsilon)$ can be calculated by
\begin{equation}\label{e:4.15}
\begin{aligned}
\tilde{f}_3^1(z,0,\epsilon)=&f_3^1(z,0,\epsilon)+\frac{3}{2}\left[D_zf_2^1(z,0,\epsilon)U_2^1(z,\epsilon)+D_yf_2^1(z,y,\epsilon)|_{y=0}U_2^2(z,\epsilon)-DU_2^1(z,\epsilon)g_2^1(z,0,\epsilon)\right].
\end{aligned}
\end{equation}

According to
\begin{equation} \label{e:4.16}
\frac{1}{3!}f_3^1(z,0,\epsilon)=\sum_{|q|+|s|=3} \frac{1}{\prod_{i=1}^2 q_i!\prod_{k=1}^2 s_k!}f^1_{qs}z^q\epsilon^s,
\end{equation}
and \eqref{e:4.4}, \eqref{e:4.5}, \eqref{e:4.6}, \eqref{e:4.14}, \eqref{e:4.15}, we have
\begin{subequations}\label{e:4.17}
\begin{equation}
\begin{aligned}
g_{2100}^{11}=&f_{2100}^{11}+2\psi_1(0)\left\langle Q(\phi_1,h_{1100})\beta_{k_1},\beta_{k_1}\right\rangle+\psi_1(0)\left\langle Q(\phi_2,h_{2000})\beta_{k_2},\beta_{k_1}\right\rangle;\\
g_{2100}^{12}=&f_{2100}^{12}+2\psi_2(0)\left\langle Q(\phi_1,h_{1100})\beta_{k_1},\beta_{k_2}\right\rangle+\psi_2(0)\left\langle Q(\phi_2,h_{2000})\beta_{k_2},\beta_{k_2}\right\rangle;\\
g_{1200}^{11}=&f_{1200}^{11}+\psi_1(0)\left\langle Q(\phi_1,h_{0200})\beta_{k_1},\beta_{k_1}\right\rangle+2\psi_1(0)\left\langle Q(\phi_2,h_{1100})\beta_{k_2},\beta_{k_1}\right\rangle;\\
g_{1200}^{12}=&f_{1200}^{12}+\psi_2(0)\left\langle Q(\phi_1,h_{0200})\beta_{k_1},\beta_{k_2}\right\rangle+2\psi_2(0)\left\langle Q(\phi_2,h_{1100})\beta_{k_2},\beta_{k_2}\right\rangle;\\
g_{3000}^{11}=&f_{3000}^{11}+3\psi_1(0) \left\langle Q(\phi_1,h_{2000})\beta_{k_1},\beta_{k_1} \right\rangle;
g_{3000}^{12}=f_{3000}^{12}+3\psi_2(0) \left\langle Q(\phi_1,h_{2000})\beta_{k_1},\beta_{k_2} \right\rangle;\\
g_{0300}^{11}=&f_{0300}^{11}+3\psi_1(0) \left\langle Q(\phi_2,h_{0200})\beta_{k_2},\beta_{k_1} \right\rangle;
g_{0300}^{12}=f_{0300}^{12}+3\psi_2(0) \left\langle Q(\phi_2,h_{0200})\beta_{k_2},\beta_{k_2} \right\rangle;
\end{aligned}
\end{equation}
and
\begin{equation}
\begin{aligned}
g_{2010}^{11}=&f_{2010}^{11}+\frac{1}{2}\psi_1(0)\left\langle L_{\epsilon_1}h_{2000}+D_{\epsilon_1}\Delta h_{2000}(0),\beta_{k_1}\right\rangle
+2\psi_1(0)\left\langle Q(\phi_1,h_{1010})\beta_{k_1},\beta_{k_1}\right\rangle;\\
g_{2010}^{12}=&f_{2010}^{12}+\frac{1}{2}\psi_2(0)\left\langle L_{\epsilon_1}h_{2000}+D_{\epsilon_1}\Delta h_{2000}(0),\beta_{k_2}\right\rangle
+2\psi_2(0)\left\langle Q(\phi_1,h_{1010})\beta_{k_1},\beta_{k_2}\right\rangle;\\
g_{2001}^{11}=&f_{2001}^{11}+\frac{1}{2}\psi_1(0)\left\langle L_{\epsilon_2}h_{2000}+D_{\epsilon_2}\Delta h_{2000}(0),\beta_{k_1}\right\rangle
+2\psi_1(0)\left\langle Q(\phi_1,h_{1001})\beta_{k_1},\beta_{k_1}\right\rangle;\\
g_{2001}^{12}=&f_{2001}^{12}+\frac{1}{2}\psi_2(0)\left\langle L_{\epsilon_2}h_{2000}+D_{\epsilon_2}\Delta h_{2000}(0),\beta_{k_2}\right\rangle
+2\psi_2(0)\left\langle Q(\phi_1,h_{1001})\beta_{k_1},\beta_{k_2}\right\rangle;\\
g_{0210}^{11}=&f_{0210}^{11}+\frac{1}{2}\psi_1(0)\left\langle L_{\epsilon_1}h_{0200}+D_{\epsilon_1}\Delta h_{0200}(0),\beta_{k_1}\right\rangle
+2\psi_1(0)\left\langle Q(\phi_2,h_{0110})\beta_{k_2},\beta_{k_1}\right\rangle;\\
g_{0210}^{12}=&f_{0210}^{12}+\frac{1}{2}\psi_2(0)\left\langle L_{\epsilon_1}h_{0200}+D_{\epsilon_1}\Delta h_{0200}(0),\beta_{k_2}\right\rangle
+2\psi_2(0)\left\langle Q(\phi_2,h_{0110})\beta_{k_2},\beta_{k_2}\right\rangle;\\
g_{0201}^{11}=&f_{0201}^{11}+\frac{1}{2}\psi_1(0)\left\langle L_{\epsilon_2}h_{0200}+D_{\epsilon_2}\Delta h_{0200}(0),\beta_{k_1}\right\rangle
+2\psi_1(0)\left\langle Q(\phi_2,h_{0101})\beta_{k_2},\beta_{k_1}\right\rangle;\\
g_{0201}^{12}=&f_{0201}^{12}+\frac{1}{2}\psi_2(0)\left\langle L_{\epsilon_2}h_{0200}+D_{\epsilon_2}\Delta h_{0200}(0),\beta_{k_2}\right\rangle
+2\psi_2(0)\left\langle Q(\phi_2,h_{0101})\beta_{k_2},\beta_{k_2}\right\rangle;\\
g_{1110}^{11}=&f_{1110}^{11}+\frac{1}{2}\psi_1(0)\left\langle L_{\epsilon_1}h_{1100}+D_{\epsilon_1}\Delta h_{1100}(0),\beta_{k_1}\right\rangle
+\psi_1(0)\left\langle Q(\phi_1,h_{0110})\beta_{k_1},\beta_{k_1}\right\rangle\\
&+\psi_1(0)\left\langle Q\left(\phi_2,h_{1010}\right)\beta_{k_2},\beta_{k_1}\right\rangle;\\
g_{1110}^{12}=&f_{1110}^{12}+\frac{1}{2}\psi_2(0)\left\langle L_{\epsilon_1}h_{1100}+D_{\epsilon_1}\Delta h_{1100}(0),\beta_{k_2}\right\rangle
+\psi_2(0)\left\langle Q(\phi_1,h_{0110})\beta_{k_1},\beta_{k_2}\right\rangle\\
&+\psi_2(0)\left\langle Q\left(\phi_2,h_{1010}\right)\beta_{k_2},\beta_{k_2}\right\rangle;\\
g_{1101}^{11}=&f_{1101}^{11}+\frac{1}{2}\psi_1(0)\left\langle L_{\epsilon_2}h_{1100}+D_{\epsilon_2}\Delta h_{1100}(0),\beta_{k_1}\right\rangle
+\psi_1(0)\left\langle Q(\phi_1,h_{0101})\beta_{k_1},\beta_{k_1}\right\rangle\\
&+\psi_1(0)\left\langle Q\left(\phi_2,h_{1001}\right)\beta_{k_2},\beta_{k_1}\right\rangle;\\
g_{1101}^{12}=&f_{1101}^{12}+\frac{1}{2}\psi_2(0)\left\langle L_{\epsilon_2}h_{1100}+D_{\epsilon_2}\Delta h_{1100}(0),\beta_{k_2}\right\rangle
+\psi_2(0)\left\langle Q(\phi_1,h_{0101})\beta_{k_1},\beta_{k_2}\right\rangle\\
&+\psi_2(0)\left\langle Q\left(\phi_2,h_{1001}\right)\beta_{k_2},\beta_{k_2}\right\rangle;
\end{aligned}
\end{equation}
and
\begin{equation}
\begin{aligned}
g_{0030}^{11}=&f_{0030}^{11}+\frac{3}{2}\psi_1(0)\left\langle L_{\epsilon_1}h_{0020}+D_{\epsilon_1}\Delta h_{0020}(0),\beta_{k_1}\right\rangle;
g_{0030}^{12}=f_{0030}^{12}+\frac{3}{2}\psi_2(0)\left\langle L_{\epsilon_1}h_{0020}+\Delta D_{\epsilon_1}h_{0020}(0),\beta_{k_2}\right\rangle;\\
g_{0003}^{11}=&f_{0003}^{11}+\frac{3}{2}\psi_1(0)\left\langle L_{\epsilon_2}h_{0002}+D_{\epsilon_2}\Delta h_{0002}(0),\beta_{k_1}\right\rangle;
g_{0003}^{12}=f_{0003}^{12}+\frac{3}{2}\psi_2(0)\left\langle L_{\epsilon_2}h_{0002}+D_{\epsilon_2}\Delta h_{0002}(0),\beta_{k_2}\right\rangle;\\
g_{1020}^{11}=&f_{1020}^{11}+\psi_1(0)\left\langle L_{\epsilon_1}h_{1010}+D_{\epsilon_1}\Delta h_{1010}(0),\beta_{k_1}\right\rangle
+\psi_1(0)\left\langle Q(\phi_1,h_{0020})\beta_{k_1},\beta_{k_1}\right\rangle;\\
g_{1020}^{12}=&f_{1020}^{12}+\psi_2(0)\left\langle L_{\epsilon_1}h_{1010}+D_{\epsilon_1}\Delta h_{1010}(0),\beta_{k_2}\right\rangle
+\psi_2(0)\left\langle Q(\phi_1,h_{0020})\beta_{k_1},\beta_{k_2}\right\rangle;\\
g_{1002}^{11}=&f_{1002}^{11}+\psi_1(0)\left\langle L_{\epsilon_2}h_{1001}+D_{\epsilon_2}\Delta h_{1001}(0),\beta_{k_1}\right\rangle
+\psi_1(0)\left\langle Q(\phi_1,h_{0002})\beta_{k_1},\beta_{k_1}\right\rangle;\\
g_{1002}^{12}=&f_{1002}^{12}+\psi_2(0)\left\langle L_{\epsilon_2}h_{1001}+D_{\epsilon_2}\Delta h_{1001}(0),\beta_{k_2}\right\rangle
+\psi_2(0)\left\langle Q(\phi_1,h_{0002})\beta_{k_1},\beta_{k_2}\right\rangle;\\
g_{0120}^{11}=&f_{0120}^{11}+\psi_1(0)\left\langle L_{\epsilon_1}h_{0110}+D_{\epsilon_1}\Delta h_{0110}(0),\beta_{k_1}\right\rangle
+\psi_1(0)\left\langle Q(\phi_2,h_{0020})\beta_{k_2},\beta_{k_1}\right\rangle;\\
g_{0120}^{12}=&f_{0120}^{12}+\psi_2(0)\left\langle L_{\epsilon_1}h_{0110}+D_{\epsilon_1}\Delta h_{0110}(0),\beta_{k_2}\right\rangle
+\psi_2(0)\left\langle Q(\phi_2,h_{0020})\beta_{k_2},\beta_{k_2}\right\rangle;\\
g_{0102}^{11}=&f_{0102}^{11}+\psi_1(0)\left\langle L_{\epsilon_2}h_{0101}+D_{\epsilon_2}\Delta h_{0101}(0),\beta_{k_1}\right\rangle
+\psi_1(0)\left\langle Q(\phi_2,h_{0002})\beta_{k_2},\beta_{k_1}\right\rangle;\\
g_{0102}^{12}=&f_{0102}^{12}+\psi_2(0)\left\langle L_{\epsilon_2}h_{0101}+D_{\epsilon_2}\Delta h_{0101}(0),\beta_{k_2}\right\rangle
+\psi_2(0)\left\langle Q(\phi_2,h_{0002})\beta_{k_2},\beta_{k_2}\right\rangle;\\
g_{1011}^{11}=&f_{1011}^{11}+\frac{1}{2}\psi_1(0)\left\langle L_{\epsilon_1}h_{1001}+D_{\epsilon_1}\Delta h_{1001}(0),\beta_{k_1}\right\rangle
+\frac{1}{2}\psi_1(0)\left\langle L_{\epsilon_2}h_{1010}+D_{\epsilon_2}\Delta h_{1010}(0),\beta_{k_1}\right\rangle\\
&+\psi_1(0)\left\langle Q(\phi_1,h_{0011})\beta_{k_1},\beta_{k_1}\right\rangle;\\
g_{1011}^{12}=&f_{1011}^{12}+\frac{1}{2}\psi_2(0)\left\langle L_{\epsilon_1}h_{1001}+D_{\epsilon_1}\Delta h_{1001}(0),\beta_{k_2}\right\rangle
+\frac{1}{2}\psi_2(0)\left\langle L_{\epsilon_2}h_{1010}+D_{\epsilon_2}\Delta h_{1010}(0),\beta_{k_2}\right\rangle\\
&+\psi_2(0)\left\langle Q(\phi_1,h_{0011})\beta_{k_1},\beta_{k_2}\right\rangle;\\
g_{0111}^{11}=&f_{0111}^{11}+\frac{1}{2}\psi_1(0)\left\langle L_{\epsilon_1}h_{0101}+D_{\epsilon_1}\Delta h_{0101}(0),\beta_{k_1}\right\rangle
+\frac{1}{2}\psi_1(0)\left\langle L_{\epsilon_2}h_{0110}+D_{\epsilon_2}\Delta h_{0110}(0),\beta_{k_1}\right\rangle\\
&+\psi_1(0)\left\langle Q(\phi_2,h_{0011})\beta_{k_2},\beta_{k_1}\right\rangle;\\
g_{0111}^{12}=&f_{0111}^{12}+\frac{1}{2}\psi_2(0)\left\langle L_{\epsilon_1}h_{0101}+D_{\epsilon_1}\Delta h_{0101}(0),\beta_{k_2}\right\rangle
+\frac{1}{2}\psi_2(0)\left\langle L_{\epsilon_2}h_{0110}+D_{\epsilon_2}\Delta h_{0110}(0),\beta_{k_2}\right\rangle\\
&+\psi_2(0)\left\langle Q(\phi_2,h_{0011})\beta_{k_2},\beta_{k_2}\right\rangle.
\end{aligned}
\end{equation}
\end{subequations}
Then, applying \eqref{e:2.14}, \eqref{e:2.3} and expansions \eqref{e:2.20}, we obtain
\begin{subequations}\label{e:4.18}
\begin{equation}
\begin{aligned}
&f_{3000}^{11}=\psi_1(0)C(\phi_1,\phi_1,\phi_1)\left\langle\beta_{k_1}^3,\beta_{k_1}\right\rangle;\qquad
& &f_{2100}^{11}=\psi_1(0)C(\phi_1,\phi_1,\phi_2)\left\langle\beta_{k_1}^2\beta_{k_2},\beta_{k_1}\right\rangle;\\
&f_{0300}^{11}=\psi_1(0)C(\phi_2,\phi_2,\phi_2)\left\langle\beta_{k_2}^3,\beta_{k_1}\right\rangle;
& &f_{1200}^{11}=\psi_1(0)C(\phi_1,\phi_2,\phi_2)\left\langle\beta_{k_1}\beta_{k_2}^2,\beta_{k_1}\right\rangle;\\
&f_{3000}^{12}=\psi_2(0)C(\phi_1,\phi_1,\phi_1)\left\langle\beta_{k_1}^3,\beta_{k_2}\right\rangle;
& &f_{2100}^{12}=\psi_2(0)C(\phi_1,\phi_1,\phi_2)\left\langle\beta_{k_1}^2\beta_{k_2},\beta_{k_2}\right\rangle;\\
&f_{0300}^{12}=\psi_2(0)C(\phi_2,\phi_2,\phi_2)\left\langle\beta_{k_2}^3,\beta_{k_2}\right\rangle;
& &f_{1200}^{12}=\psi_2(0)C(\phi_1,\phi_2,\phi_2)\left\langle\beta_{k_1}\beta_{k_2}^2,\beta_{k_2}\right\rangle;
\end{aligned}
\end{equation}
and
\begin{equation}
\begin{aligned}
&f_{1011}^{11}=\frac{1}{3}\psi_1(0)\left(L_{\epsilon_1\epsilon_2}\phi_1-\mu_{k_1}D_{\epsilon_1\epsilon_2}\phi_1(0)\right)\left\langle\beta_{k_1},\beta_{k_1}\right\rangle;
& \quad &f_{1020}^{11}=\frac{1}{3}\psi_1(0)\left(L_{\epsilon_1^2}\phi_1-\mu_{k_1}D_{\epsilon_1^2}\phi_1(0)\right)\left\langle\beta_{k_1},\beta_{k_1}\right\rangle;\\
&f_{1002}^{11}=\frac{1}{3}\psi_1(0)\left(L_{\epsilon_2^2}\phi_1-\mu_{k_1}D_{\epsilon_2^2}\phi_1(0)\right)\left\langle\beta_{k_1},\beta_{k_1}\right\rangle;
& &f_{0111}^{11}=\frac{1}{3}\psi_1(0)\left(L_{\epsilon_1\epsilon_2}\phi_2-\mu_{k_2}D_{\epsilon_1\epsilon_2}\phi_2(0)\right)\left\langle\beta_{k_2},\beta_{k_1}\right\rangle;\\
&f_{0120}^{11}=\frac{1}{3}\psi_1(0)\left(L_{\epsilon_1^2}\phi_2-\mu_{k_2}D_{\epsilon_1^2}\phi_2(0)\right)\left\langle\beta_{k_2},\beta_{k_1}\right\rangle;
& &f_{0102}^{11}=\frac{1}{3}\psi_1(0)\left(L_{\epsilon_2^2}\phi_2-\mu_{k_2}D_{\epsilon_2^2}\phi_2(0)\right)\left\langle\beta_{k_2},\beta_{k_1}\right\rangle;\\
&f_{1011}^{12}=\frac{1}{3}\psi_2(0)\left(L_{\epsilon_1\epsilon_2}\phi_1-\mu_{k_1}D_{\epsilon_1\epsilon_2}\phi_1(0)\right)\left\langle\beta_{k_1},\beta_{k_2}\right\rangle;
& &f_{1020}^{12}=\frac{1}{3}\psi_2(0)\left(L_{\epsilon_1^2}\phi_1-\mu_{k_1}D_{\epsilon_1^2}\phi_1(0)\right)\left\langle\beta_{k_1},\beta_{k_2}\right\rangle;\\
&f_{1002}^{12}=\frac{1}{3}\psi_2(0)\left(L_{\epsilon_2^2}\phi_1-\mu_{k_1}D_{\epsilon_2^2}\phi_1(0)\right)\left\langle\beta_{k_1},\beta_{k_2}\right\rangle;
& &f_{0111}^{12}=\frac{1}{3}\psi_2(0)\left(L_{\epsilon_1\epsilon_2}\phi_2-\mu_{k_2}D_{\epsilon_1\epsilon_2}\phi_2(0)\right)\left\langle\beta_{k_2},\beta_{k_2}\right\rangle;\\
&f_{0120}^{12}=\frac{1}{3}\psi_2(0)\left(L_{\epsilon_1^2}\phi_2-\mu_{k_2}D_{\epsilon_1^2}\phi_2(0)\right)\left\langle\beta_{k_2},\beta_{k_2}\right\rangle;
& &f_{0102}^{12}=\frac{1}{3}\psi_2(0)\left(L_{\epsilon_2^2}\phi_2-\mu_{k_2}D_{\epsilon_2^2}\phi_2(0)\right)\left\langle\beta_{k_2},\beta_{k_2}\right\rangle.
\end{aligned}
\end{equation}
\end{subequations}
Plugging \eqref{e:4.5} and \eqref{e:4.18} into \eqref{e:4.17}, Theorem \ref{thm:2.1} follows.

\section{Conclusion}
When two Turing modes interact, there appear superposition patterns which reveal complex dynamical phenomena.
Yang et al. \cite{YL&DM&ZAM&EIR} firstly proposed a model to reproduce black-eye patterns in 2D domain, which were firstly observed in a CIMA reaction \cite{GGH&OQ&SHL}. Also, they observed a variety of other spatial patterns resulting from interactions between different Turing modes, including white-eye patterns and superposition patterns, via numerical simulations.
In the paper, to theoretically investigate superposition patterns arising from two interacting Turing modes, we explore normal form of some parameterized PFDEs at Turing-Turing singularity, then study spatial patterns of diffusive predator-prey system \eqref{e:3.1} by analyzing the obtained normal forms.
Based on Faria's work \cite{FT, FT01} and center manifold theory \cite{LXD&SJWH&WJH, FT&HWZ&WJH}, we derive the third-order normal form of Turing-Turing bifurcation, which is locally topologically equivalent to the original parameterized PFDEs. Also, an explicit algorithm of computing the third-order normal form is provided.
Moreover, the third-order normal form is simplified as three normal forms describing essentially different spatial phenomena, when considering one dimensional domain and Neumann boundary conditions.
Meanwhile, several concise formulas for computing coefficients of these three normal forms are also derived, which are expressed in explicit form of the original system parameters.
We emphasize that the process of computing coefficients of normal forms utilizing these concise formulas could be implemented by computer programs, and these formulas also apply to computing coefficients of normal forms for partial differential equations.

Then, we investigate spatial patterns of diffusive predator-prey system \eqref{e:3.1} with Crowley-Martin functional response near Turing-Turing singularity, with the aid of these three normal forms.
We find that, two of these three normal forms arise in predator-prey system \eqref{e:3.1} for two different sets of parameters. For system parameters $m=6,a=3,b=0.5,d_2=0.7,d_1=0.0051,s=0.2064$, four stable spatially inhomogeneous steady states with different single characteristic wavelengths coexist (see Fig. \ref{fig:3}), which demonstrates our conjecture in \cite{CX&JWH}. For parameters $m=5,a=3,b=0.1,d_2=4$, system \eqref{e:3.1} admits superposition steady states with the shape of $\phi_1h_1\cos x+\phi_2h_2\cos2x-$like and tri-stable patterns that a pair of stable superposition steady states with the shape of $\phi_1h_1\cos x+\phi_2h_2\cos2x-$like coexists with the stable coexistence equilibrium (see Fig. \ref{fig:6}) or a stable steady state with the shape of $\phi_2\cos2x-$like (see Fig. \ref{fig:8}), as well as quad-stable patterns that a pair of stable superposition steady states with the shape of $\phi_1h_1\cos x+\phi_2h_2\cos2x-$like and a pair of stable steady states with the shape of $\phi_2\cos2x-$like coexist (see Fig. \ref{fig:7}).
It is worth noting that, bifurcation set and phase portraits of normal form \eqref{e:3.8} might be incomplete, which indicates that system \eqref{e:3.1} might admit other complex superposition patterns.
And, complete bifurcation set and the corresponding phase portraits require further analysis.

\section*{Acknowledgements}
The authors are supported by the National Natural Science Foundation of China (No.11871176) and the Fundamental Research Funds for the Central Universities (FRFCU5630103121).

\section*{References}
\biboptions{numbers,sort&compress}


\begin{thebibliography}{41}
\providecommand{\natexlab}[1]{#1}
\providecommand{\url}[1]{\texttt{#1}}
\providecommand{\urlprefix}{URL }
\expandafter\ifx\csname urlstyle\endcsname\relax
  \providecommand{\doi}[1]{doi:\discretionary{}{}{}#1}\else
  \providecommand{\doi}[1]{doi:\discretionary{}{}{}\begingroup
  \urlstyle{rm}\url{#1}\endgroup}\fi
\providecommand{\bibinfo}[2]{#2}

\bibitem[{Turing(1952)}]{TAM}
\bibinfo{author}{A.~M. Turing}, \bibinfo{title}{The chemical basis of
  morphogenesis}, \bibinfo{journal}{Philos. Trans. Roy. Soc. London Ser. B}
  \bibinfo{volume}{237}~(\bibinfo{number}{641}) (\bibinfo{year}{1952})
  \bibinfo{pages}{37--72}.

\bibitem[{Castets et~al.(1990)Castets, Dulos, Boissonade, and
  Dekepper}]{CV&DE&BJ}
\bibinfo{author}{V.~Castets}, \bibinfo{author}{E.~Dulos},
  \bibinfo{author}{J.~Boissonade}, \bibinfo{author}{P.~Dekepper},
  \bibinfo{title}{Experimental-evidence of a sustained standing Turing-type
  nonequilibrium chemical-pattern}, \bibinfo{journal}{Phys. Rev. Lett.}
  \bibinfo{volume}{64}~(\bibinfo{number}{24}) (\bibinfo{year}{1990})
  \bibinfo{pages}{2953--2956}.

\bibitem[{Dekepper et~al.(1991)Dekepper, Castets, Dulos, and
  Boissonade}]{DP&CV&DE}
\bibinfo{author}{P.~Dekepper}, \bibinfo{author}{V.~Castets},
  \bibinfo{author}{E.~Dulos}, \bibinfo{author}{J.~Boissonade},
  \bibinfo{title}{Turing-type chemical-patterns in the chlorite-iodide-malonic
  acid reaction}, \bibinfo{journal}{Phys. D}
  \bibinfo{volume}{49}~(\bibinfo{number}{1-2}) (\bibinfo{year}{1991})
  \bibinfo{pages}{161--169}.

\bibitem[{Lengyel and Epstein(1991)}]{LI&EIR}
\bibinfo{author}{I.~Lengyel}, \bibinfo{author}{I.~R. Epstein},
  \bibinfo{title}{Modeling of Turing structures in the chlorite iodide
  malonic-acid starch reaction system}, \bibinfo{journal}{Science}
  \bibinfo{volume}{251}~(\bibinfo{number}{4994}) (\bibinfo{year}{1991})
  \bibinfo{pages}{650--652}.

\bibitem[{Gunaratne et~al.(1994)Gunaratne, Ouyang, and Swinney}]{GGH&OQ&SHL}
\bibinfo{author}{G.~H. Gunaratne}, \bibinfo{author}{Q.~Ouyang},
  \bibinfo{author}{H.~L. Swinney}, \bibinfo{title}{Pattern formation in the
  presence of symmetries}, \bibinfo{journal}{Phys. Rev. E}
  \bibinfo{volume}{50}~(\bibinfo{number}{4}) (\bibinfo{year}{1994})
  \bibinfo{pages}{2802--2820}.

\bibitem[{Yang et~al.(2002)Yang, Dolnik, Zhabotinsky, and
  Epstein}]{YL&DM&ZAM&EIR}
\bibinfo{author}{L.~Yang}, \bibinfo{author}{M.~Dolnik}, \bibinfo{author}{A.~M.
  Zhabotinsky}, \bibinfo{author}{I.~R. Epstein}, \bibinfo{title}{Spatial
  resonances and superposition patterns in a reaction-diffusion model with
  interacting Turing modes}, \bibinfo{journal}{Phys. Rev. Lett.}
  \bibinfo{volume}{88}~(\bibinfo{number}{20}) (\bibinfo{year}{2002})
  \bibinfo{pages}{208303}.

\bibitem[{Yang and Song(2016)}]{YR&SYL}
\bibinfo{author}{R.~Yang}, \bibinfo{author}{Y.~Song}, \bibinfo{title}{Spatial
  resonance and Turing-Hopf bifurcations in the Gierer-Meinhardt model},
  \bibinfo{journal}{Nonlinear Anal. Real World Appl.} \bibinfo{volume}{31}
  (\bibinfo{year}{2016}) \bibinfo{pages}{356--387}.

\bibitem[{Wei et~al.(2015)Wei, Wu, and Guo}]{WMH&WJH&GGH}
\bibinfo{author}{M.~Wei}, \bibinfo{author}{J.~Wu}, \bibinfo{author}{G.~Guo},
  \bibinfo{title}{Steady state bifurcations for a glycolysis model in
  biochemical reaction}, \bibinfo{journal}{Nonlinear Anal. Real World Appl.}
  \bibinfo{volume}{22} (\bibinfo{year}{2015}) \bibinfo{pages}{155--175}.

\bibitem[{Guo et~al.(2018)Guo, Liu, Li, and Li}]{GGH&LL&LBF&LJ}
\bibinfo{author}{G.~Guo}, \bibinfo{author}{L.~Liu}, \bibinfo{author}{B.~Li},
  \bibinfo{author}{J.~Li}, \bibinfo{title}{Qualitative analysis on positive
  steady-state solutions for an autocatalysis model with high order},
  \bibinfo{journal}{Nonlinear Anal. Real World Appl.} \bibinfo{volume}{41}
  (\bibinfo{year}{2018}) \bibinfo{pages}{665--691}.

\bibitem[{Li et~al.(2015)Li, Wu, and Doug}]{LSB&WJH&DYY}
\bibinfo{author}{S.~Li}, \bibinfo{author}{J.~Wu}, \bibinfo{author}{Y.~Doug},
  \bibinfo{title}{Turing patterns in a reaction-diffusion model with the
  Degn-Harrison reaction scheme}, \bibinfo{journal}{J. Differ. Equ.}
  \bibinfo{volume}{259}~(\bibinfo{number}{5}) (\bibinfo{year}{2015})
  \bibinfo{pages}{1990--2029}.

\bibitem[{Yang(2018)}]{YWB}
\bibinfo{author}{W.~Yang}, \bibinfo{title}{Analysis on existence of bifurcation
  solutions for a predator-prey model with herd behavior},
  \bibinfo{journal}{Appl. Math. Model.} \bibinfo{volume}{53}
  (\bibinfo{year}{2018}) \bibinfo{pages}{433--446}.

\bibitem[{Dangelmayr(1987)}]{DG}
\bibinfo{author}{G.~Dangelmayr}, \bibinfo{title}{Degenerate bifurcations near a
  double eigenvalue in the brusselator}, \bibinfo{journal}{J. Aust. Math. Soc.}
  \bibinfo{volume}{28} (\bibinfo{year}{1987}) \bibinfo{pages}{486--535}.

\bibitem[{Gambino et~al.(2013)Gambino, Lombardo, and
  Sammartino}]{GG&LMC&SM&NA13}
\bibinfo{author}{G.~Gambino}, \bibinfo{author}{M.~C. Lombardo},
  \bibinfo{author}{M.~Sammartino}, \bibinfo{title}{Pattern formation driven by
  cross-diffusion in a 2D domain}, \bibinfo{journal}{Nonlinear Anal. Real World
  Appl.} \bibinfo{volume}{14}~(\bibinfo{number}{3}) (\bibinfo{year}{2013})
  \bibinfo{pages}{1755--1779}.

\bibitem[{Gambino et~al.(2018)Gambino, Lombardo, and
  Sammartino}]{GG&LMC&SM&PRE18}
\bibinfo{author}{G.~Gambino}, \bibinfo{author}{M.~C. Lombardo},
  \bibinfo{author}{M.~Sammartino}, \bibinfo{title}{Cross-diffusion-induced
  subharmonic spatial resonances in a predator-prey system},
  \bibinfo{journal}{Phys. Rev. E} \bibinfo{volume}{97}~(\bibinfo{number}{1})
  (\bibinfo{year}{2018}) \bibinfo{pages}{012220}.

\bibitem[{Carr(1982)}]{CJ}
\bibinfo{author}{J.~Carr}, \bibinfo{title}{Applications of Centre Manifold
  Theory}, \bibinfo{publisher}{Spring -Verlag}, \bibinfo{address}{New York},
  \bibinfo{year}{1982}.

\bibitem[{Campbell and Yuan(2008)}]{CSA&YY}
\bibinfo{author}{S.~A. Campbell}, \bibinfo{author}{Y.~Yuan},
  \bibinfo{title}{Zero singularities of codimension two and three in delay
  differential equations}, \bibinfo{journal}{Nonlinearity}
  \bibinfo{volume}{21}~(\bibinfo{number}{11}) (\bibinfo{year}{2008})
  \bibinfo{pages}{2671--2691}.

\bibitem[{Chow et~al.(1994)Chow, Li, and Wang}]{CSN&LCZ&WD}
\bibinfo{author}{S.-N. Chow}, \bibinfo{author}{C.~Li},
  \bibinfo{author}{D.~Wang}, \bibinfo{title}{Normal Forms and Bifurcation of
  Planar Vector Fields}, \bibinfo{publisher}{Cambridge University Press},
  \bibinfo{year}{1994}.

\bibitem[{Faria and Magalhaes(1995)}]{FT&MLT&H}
\bibinfo{author}{T.~Faria}, \bibinfo{author}{L.~T. Magalhaes},
  \bibinfo{title}{Normal forms for retarded functional-differential equations
  with parameters and applications to Hopf-bifurcation}, \bibinfo{journal}{J.
  Differ. Equ.} \bibinfo{volume}{122}~(\bibinfo{number}{2})
  (\bibinfo{year}{1995}) \bibinfo{pages}{181--200}.

\bibitem[{Hassard et~al.(1981)Hassard, Kazarinoff, and Wan}]{HBD&KND&WYH}
\bibinfo{author}{B.~D. Hassard}, \bibinfo{author}{N.~D. Kazarinoff},
  \bibinfo{author}{Y.~H. Wan}, \bibinfo{title}{Theory and Applications of Hopf
  Bifurcation}, vol.~\bibinfo{volume}{41}, \bibinfo{publisher}{Cambridge
  University Press}, \bibinfo{address}{New York}, \bibinfo{year}{1981}.

\bibitem[{Jiang and Yuan(2007)}]{JWH&YY}
\bibinfo{author}{W.~Jiang}, \bibinfo{author}{Y.~Yuan},
  \bibinfo{title}{Bogdanov-Takens singularity in Van der Pol's oscillator with
  delayed feedback}, \bibinfo{journal}{Phys. D}
  \bibinfo{volume}{227}~(\bibinfo{number}{2}) (\bibinfo{year}{2007})
  \bibinfo{pages}{149--161}.

\bibitem[{Wang and Wei(2008)}]{WCC&WJJ}
\bibinfo{author}{C.~Wang}, \bibinfo{author}{J.~Wei}, \bibinfo{title}{Normal
  forms for NFDEs with parameters and application to the lossless transmission
  line}, \bibinfo{journal}{Nonlinear Dyn.}
  \bibinfo{volume}{52}~(\bibinfo{number}{3}) (\bibinfo{year}{2008})
  \bibinfo{pages}{199--206}.

\bibitem[{Wang and Jiang(2010)}]{WHB&JWH}
\bibinfo{author}{H.~Wang}, \bibinfo{author}{W.~Jiang},
  \bibinfo{title}{Hopf-pitchfork bifurcation in van der Pol's oscillator with
  nonlinear delayed feedback}, \bibinfo{journal}{J. Math. Anal. Appl.}
  \bibinfo{volume}{368}~(\bibinfo{number}{1}) (\bibinfo{year}{2010})
  \bibinfo{pages}{9--18}.

\bibitem[{Yuan et~al.(2015)Yuan, Jiang, and Wang}]{YR&JWH&WY}
\bibinfo{author}{R.~Yuan}, \bibinfo{author}{W.~Jiang},
  \bibinfo{author}{Y.~Wang}, \bibinfo{title}{Saddle-node-Hopf bifurcation in a
  modified Leslie-Gower predator-prey model with time-delay and prey
  harvesting}, \bibinfo{journal}{J. Math. Anal. Appl.}
  \bibinfo{volume}{422}~(\bibinfo{number}{2}) (\bibinfo{year}{2015})
  \bibinfo{pages}{1072--1090}.

\bibitem[{Faria et~al.(2002)Faria, Huang, and Wu}]{FT&HWZ&WJH}
\bibinfo{author}{T.~Faria}, \bibinfo{author}{W.~Huang},
  \bibinfo{author}{J.~Wu}, \bibinfo{title}{Smoothness of center manifolds for
  maps and formal adjoints for semilinear FDEs in general Banach spaces},
  \bibinfo{journal}{SIAM J. Math. Anal.}
  \bibinfo{volume}{34}~(\bibinfo{number}{1}) (\bibinfo{year}{2002})
  \bibinfo{pages}{173--203}.

\bibitem[{Lin et~al.(1992)Lin, So, and Wu}]{LXD&SJWH&WJH}
\bibinfo{author}{X.~Lin}, \bibinfo{author}{J.~W.~H. So},
  \bibinfo{author}{J.~Wu}, \bibinfo{title}{Center manifolds for
  partial-differential equations with delays}, \bibinfo{journal}{Proc. R. Soc.
  Edinb. Sect. A Math.} \bibinfo{volume}{122} (\bibinfo{year}{1992})
  \bibinfo{pages}{237--254}.

\bibitem[{Faria(2000)}]{FT}
\bibinfo{author}{T.~Faria}, \bibinfo{title}{Normal forms and Hopf bifurcation
  for partial differential equations with delays}, \bibinfo{journal}{Trans. Am.
  Math. Soc.} \bibinfo{volume}{352}~(\bibinfo{number}{5})
  (\bibinfo{year}{2000}) \bibinfo{pages}{2217--2238}.

\bibitem[{Faria(2001)}]{FT01}
\bibinfo{author}{T.~Faria}, \bibinfo{title}{Normal forms for semilinear
  functional differential equations in Banach spaces and applications. Part
  II}, \bibinfo{journal}{Discrete Contin. Dyn. Syst.}
  \bibinfo{volume}{7}~(\bibinfo{number}{1}) (\bibinfo{year}{2001})
  \bibinfo{pages}{155--176}.

\bibitem[{Wu(1996)}]{WJH}
\bibinfo{author}{J.~Wu}, \bibinfo{title}{Theory and Applications of Partial
  Functional Differential Equations}, vol. \bibinfo{volume}{119},
  \bibinfo{publisher}{Springer-Verlag}, \bibinfo{address}{New York},
  \bibinfo{year}{1996}.

\bibitem[{Yi et~al.(2009)Yi, Wei, and Shi}]{YFQ&WJJ&SJP}
\bibinfo{author}{F.~Yi}, \bibinfo{author}{J.~Wei}, \bibinfo{author}{J.~Shi},
  \bibinfo{title}{Bifurcation and spatiotemporal patterns in a homogeneous
  diffusive predator-prey system}, \bibinfo{journal}{J. Differ. Equ.}
  \bibinfo{volume}{246}~(\bibinfo{number}{5}) (\bibinfo{year}{2009})
  \bibinfo{pages}{1944--1977}.

\bibitem[{Su et~al.(2009)Su, Wei, and Shi}]{SY&WJJ&SJP}
\bibinfo{author}{Y.~Su}, \bibinfo{author}{J.~Wei}, \bibinfo{author}{J.~Shi},
  \bibinfo{title}{Hopf bifurcations in a reaction-diffusion population model
  with delay effect}, \bibinfo{journal}{J. Differ. Equ.}
  \bibinfo{volume}{247}~(\bibinfo{number}{4}) (\bibinfo{year}{2009})
  \bibinfo{pages}{1156--1184}.

\bibitem[{Wang(2017)}]{WJF}
\bibinfo{author}{J.~Wang}, \bibinfo{title}{Spatiotemporal patterns of a
  homogeneous diffusive predator-prey system with Holling type III functional
  response}, \bibinfo{journal}{J. Dynam. Differential Equations}
  \bibinfo{volume}{29}~(\bibinfo{number}{4}) (\bibinfo{year}{2017})
  \bibinfo{pages}{1383--1409}.

\bibitem[{Zou and Guo(2017)}]{ZR&GSJ}
\bibinfo{author}{R.~Zou}, \bibinfo{author}{S.~J. Guo},
  \bibinfo{title}{Bifurcation of reaction cross-diffusion systems},
  \bibinfo{journal}{Int. J. Bifurcation Chaos}
  \bibinfo{volume}{27}~(\bibinfo{number}{4}) (\bibinfo{year}{2017})
  \bibinfo{pages}{1750049}.

\bibitem[{Shi and Ruan(2015)}]{SHB&RSG}
\bibinfo{author}{H.~Shi}, \bibinfo{author}{S.~Ruan}, \bibinfo{title}{Spatial,
  temporal and spatiotemporal patterns of diffusive predator-prey models with
  mutual interference}, \bibinfo{journal}{IMA J. Appl. Math.}
  \bibinfo{volume}{80}~(\bibinfo{number}{5}) (\bibinfo{year}{2015})
  \bibinfo{pages}{1534--1568}.

\bibitem[{Chen and Shi(2012)}]{CSS&SJP}
\bibinfo{author}{S.~Chen}, \bibinfo{author}{J.~Shi}, \bibinfo{title}{Stability
  and Hopf bifurcation in a diffusive logistic population model with nonlocal
  delay effect}, \bibinfo{journal}{J. Differ. Equ.}
  \bibinfo{volume}{253}~(\bibinfo{number}{12}) (\bibinfo{year}{2012})
  \bibinfo{pages}{3440--3470}.

\bibitem[{Li et~al.(2013)Li, Jiang, and Shi}]{LX&JWH&SJP}
\bibinfo{author}{X.~Li}, \bibinfo{author}{W.~Jiang}, \bibinfo{author}{J.~Shi},
  \bibinfo{title}{Hopf bifurcation and Turing instability in the
  reaction-diffusion Holling-Tanner predator-prey model}, \bibinfo{journal}{IMA
  J. Appl. Math.} \bibinfo{volume}{78}~(\bibinfo{number}{2})
  (\bibinfo{year}{2013}) \bibinfo{pages}{287--306}.

\bibitem[{Jiang et~al.(2018)Jiang, An, and Shi}]{JWH&AQ&SJP}
\bibinfo{author}{W.~Jiang}, \bibinfo{author}{Q.~An}, \bibinfo{author}{J.~Shi},
  \bibinfo{title}{Formulation of the normal forms of Turing-Hopf bifurcation in
  reaction-diffusion systems with time delay}, \bibinfo{journal}{2018,
  arXiv:1802.10286}.

\bibitem[{Cao and Jiang(2018)}]{CX&JWH}
\bibinfo{author}{X.~Cao}, \bibinfo{author}{W.~Jiang},
  \bibinfo{title}{Turing-Hopf bifurcation and spatiotemporal patterns in a
  diffusive predator-prey system with Crowley-Martin functional response},
  \bibinfo{journal}{Nonlinear Anal. Real World Appl.} \bibinfo{volume}{43}
  (\bibinfo{year}{2018}) \bibinfo{pages}{428--450}.

\bibitem[{Hale and Lunel(1993)}]{HJK&LSMY}
\bibinfo{author}{J.~K. Hale}, \bibinfo{author}{S.~M.~V. Lunel},
  \bibinfo{title}{Introduction to Functional Differential Equations},
  \bibinfo{publisher}{Springer Science \& Business Media},
  \bibinfo{address}{New York}, \bibinfo{year}{1993}.

\bibitem[{May(1973)}]{MR}
\bibinfo{author}{R.~M. May}, \bibinfo{title}{Stability and Complexity in Model
  Ecosystems}, \bibinfo{publisher}{Princeton University press},
  \bibinfo{address}{Princeton, NJ}, \bibinfo{year}{1973}.

\bibitem[{Crowley and Martin(1989)}]{CPH&MEK}
\bibinfo{author}{P.~H. Crowley}, \bibinfo{author}{E.~K. Martin},
  \bibinfo{title}{Functional responses and interference within and between year
  classes of a dragonfly population}, \bibinfo{journal}{J. N. Am. Benthol.
  Soc.} \bibinfo{volume}{8}~(\bibinfo{number}{3}) (\bibinfo{year}{1989})
  \bibinfo{pages}{211--221}.

\bibitem[{Guckenheimer and Holmes(1983)}]{GJ&HP}
\bibinfo{author}{J.~Guckenheimer}, \bibinfo{author}{P.~Holmes},
  \bibinfo{title}{Nonlinear Oscillations Dynamical Systems, and Bifurcations of
  Vector Fields}, \bibinfo{publisher}{Springer-Verlag}, \bibinfo{address}{New
  York}, \bibinfo{year}{1983}.

\end{thebibliography}

\end{document}